\documentclass[a4paper,fleqn]{cas-sc}

\usepackage[authoryear,longnamesfirst]{natbib}

\PassOptionsToPackage{no-math}{fontspec}
\usepackage{kotex} 

\begin{document}
\let\WriteBookmarks\relax
\def\floatpagepagefraction{1}
\def\textpagefraction{.001}

\shorttitle{Compatible by Construction: Bottom-Up NM-ROM}    

\shortauthors{Y. Kim and Y. Choi}

\title [mode = title]{Compatible by Construction: Bottom-Up Nonlinear Model Reduction without Interface Bases}  



%

\author[1]{Youngkyu Kim}[orcid=0000-0002-4825-4072]



\ead{kim101@llnl.gov}


\credit{Writing - review \& editing, Writing - original draft, Visualization, 
Validation, Supervision, Software, Resources, Methodology, Investigation, Formal analysis, 
Data curation, Conceptualization.}

\affiliation[1]{organization={Center for Applied Scientific Computing, 
Lawrence Livermore National Laboratory},
            city={Livermore},
            postcode={94550}, 
            state={CA},
            country={USA}}

\author[1]{Youngsoo Choi}[orcid=0000-0001-8797-7970]

\cormark[1]


\ead{choi15@llnl.gov}


\credit{Writing - review \& editing, Supervision, Resources, Methodology, Formal analysis, 
Conceptualization, Project administration, Funding acquisition.}


\cortext[1]{Corresponding author}



\begin{abstract}
Nonlinear manifolds can provide compact representations of solutions to partial
differential equations (PDEs), but generating training data and learning these
representations on large domains remain costly.
We present a bottom-up framework that reuses local nonlinear representations
learned on a small training domain to solve larger PDE systems.
A shared nonlinear decoder represents the full set of degrees of freedom owned
by each subdomain, without separate interior and interface models.
Each global degree of freedom has a unique owner, and local residuals directly
use decoded state values from the corresponding neighboring subdomains.
The resulting discrete state is single-valued by construction, eliminating the
need for independent interface copies and additional interface compatibility constraints.
Physical coupling is retained through the original time-discrete PDE residual.
With the decoder held fixed, the subdomain-specific reduced coordinates are
determined jointly by minimizing this residual, and the Gauss--Newton system is
assembled from local contributions.
Experiments on parameterized viscous Burgers and incompressible Navier--Stokes
problems demonstrate reuse on larger domains without retraining.
Comparisons with local proper orthogonal decomposition models under the same
coupling construction show favorable accuracy--efficiency trade-offs in the tested regimes.
\end{abstract}




\ExplSyntaxOn
\bool_gset_true:N \g_stm_nologo_bool
\ExplSyntaxOff
\maketitle

\section{Introduction}
\label{sec:introduction}
High-fidelity partial differential equation (PDE) simulations play a central role in many areas of science and engineering, including fluid dynamics, structural analysis, and energy systems.
However, their computational cost increases with spatial resolution, domain size, and integration time, and this burden is particularly severe in many-query applications such as design optimization, uncertainty quantification, and real-time control \citep{ohlberger2015reduced-039,hesthaven2022reduced-4e1}.
State estimation, optimization, and decision-making in digital twins likewise require models that can make repeated predictions under limited computational resources \citep{willcox2024role-b0a,ferrari2024digital-655,kramer2026reduced-order-9c8}.
Reduced-order models (ROMs) alleviate these many-query computational costs by approximating high-dimensional states in low-dimensional spaces, and a broad range of linear and nonlinear approaches have been studied \citep{refId0,hesthaven2022reduced-4e1}.

Nonlinear ROMs have attracted attention because they can represent complex solution structures in low dimensions.
Linear subspaces are effective when the solution set has a rapidly decaying Kolmogorov \(n\)-width, but problems with moving structures or sharp spatial variations may require many basis vectors for an accurate approximation \citep{lee2020model-ab0,kim2022fast-c42,refId0}.
Autoencoder (AE)-based nonlinear manifold approximation offers a way to mitigate these limitations and represent full-order states with fewer reduced coordinates \citep{lee2020model-ab0,kim2022fast-c42,romor2023non-linear-012}.
Learned representations have also expanded to graph-convolutional autoencoders, smooth neural fields, and auto-decoder-based nonlinear trial manifolds \citep{pichi2024graph-30c,puri2025snf-rom:-640,ye2024meta-auto-decoder:-f5c}.

However, learning nonlinear representations requires high-fidelity simulation data, and generating these data is itself computationally expensive \citep{cocola2023hyper-reduced-15f}.
Moreover, as the computational domain grows, not only the state dimension of the full-order model (FOM) but also snapshot storage and neural network training costs may increase \citep{diaz2024fast-e0a}.
Thus, generating data directly on a large domain and learning a nonlinear representation of its full state can limit offline scalability, even when the online computation is efficient.
This calls for a construction that confines data generation and training to small problems while retaining the expressive power of nonlinear ROMs.

\begin{figure}[pos=htbp]
\centering
\includegraphics[width=0.80\textwidth]{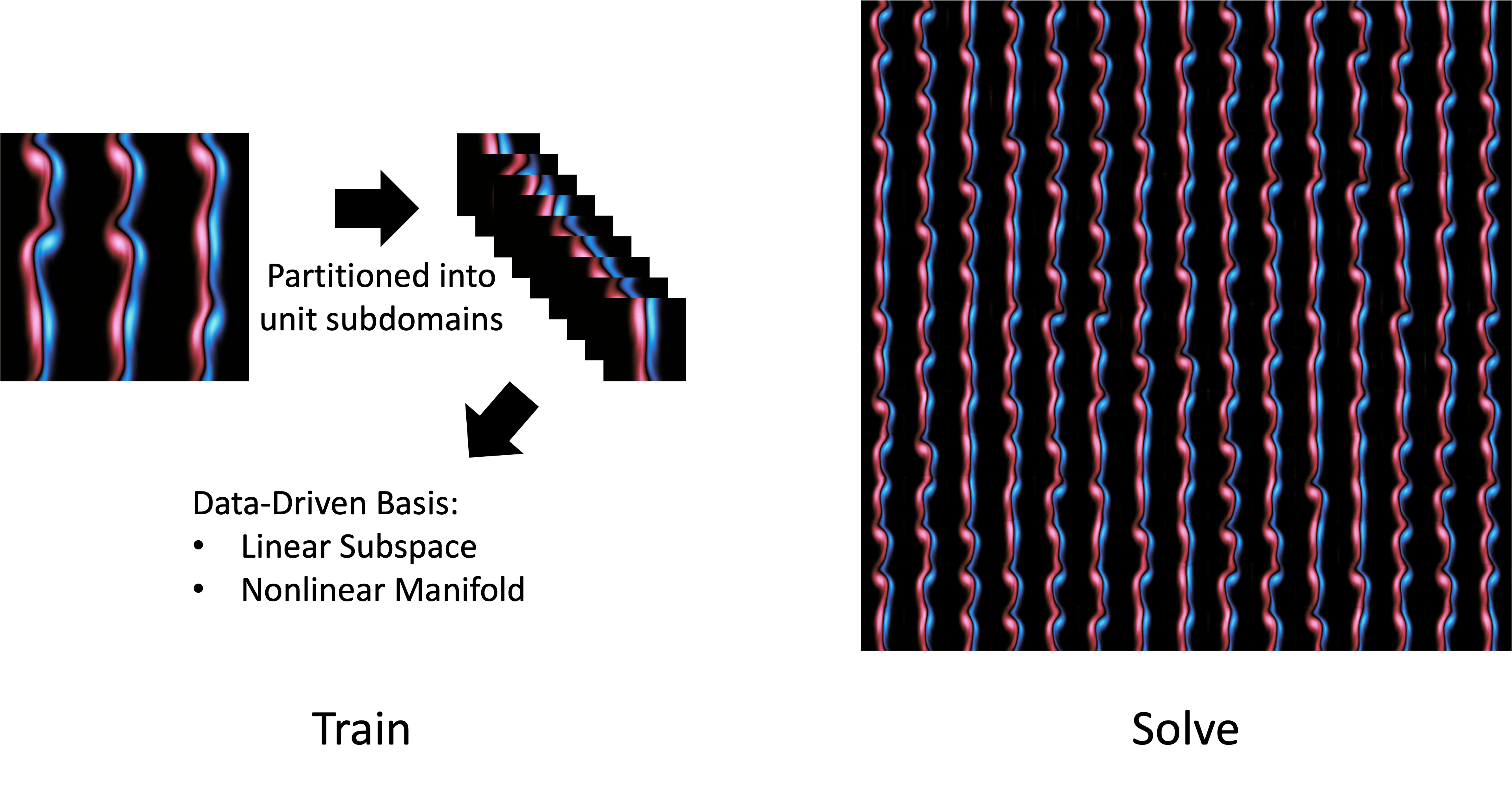}
\caption{Schematic illustration of the train-small, solve-big strategy.
A local representation learned from a small training domain is reused on the subdomains of a larger domain without retraining.
All subdomains share the decoder but have independent reduced coordinates.
Neighboring state values required by a local residual are obtained from their owning subdomains, and coupling is retained through the original time-discrete PDE residual.
The illustrated Navier--Stokes state consists of \(u,v,p\). Vorticity is computed afterward from the velocity for visualization.}
\label{fig:train_small_solve_big}
\end{figure}

This work follows a bottom-up approach in which local representations learned on a small domain are repeatedly reused on a larger computational domain.
This approach is aligned with the train-small, model-big perspective, which combines reduced representations obtained from small components or unit problems to form a larger system \citep{chung2024train-073}.
The proposed method provides a nonlinear-manifold realization of the data-driven finite element method (DD-FEM) framework described by \cite{choi2026rigor-4b6}, which combines local learning with global assembly based on the governing equations. In this realization, a shared decoder supplies the local representations, while the original time-discrete FOM residual retains inter-subdomain coupling in a global residual minimization problem.
The key idea is to partition global snapshots generated on a small domain into subdomain states, learn a local state representation, and then deploy that representation on each subdomain of a larger domain.
Rather than training a separate decoder for the larger domain, we determine the reduced coordinates of each subdomain from the governing equations.
Figure~\ref{fig:train_small_solve_big} illustrates this relationship between small-domain learning and large-domain computation.
The figure uses a Navier--Stokes flow as an example where vorticity is a visualization quantity computed from the reconstructed velocity, not a learned state variable.

Subdomain-wise reduced representations and their coupling have also been studied in DD-LSPG and nonlinear manifold DD-ROMs \citep{hoang2021domain-decomposition-838,diaz2024fast-e0a}.
In particular, \citet{zanardi2024scalable} presented a bottom-up strategy that reuses nonlinear representations learned on a small domain for a larger Burgers problem.
The focus of this work is therefore not the introduction of the bottom-up strategy itself, but the local state representation and inter-subdomain coupling used to implement it.
Existing nonlinear manifold ROMs with separate interior and port representations distinguish interior states from interface or port states and impose compatibility constraints between subdomain-specific interface copies \citep{diaz2024fast-e0a,zanardi2024scalable}.
\citet{diaz2024fast-e0a} presented both a formulation that weakly enforces full-order port compatibility and a formulation that ensures strong compatibility through constraints on reduced port coordinates.
Rather than changing the strength of these constraints, we use a unified representation of the full state owned by each subdomain and do not introduce independent interface copies.
This removes the need for separate learned models for the interior and interface and for compatibility constraints between copies.

Meanwhile, \citet{chung2024train-073} reuse proper orthogonal decomposition (POD)-based linear representations learned on small components and weakly impose inter-subdomain conditions through discontinuous Galerkin (DG) interface operators, without separate reduced bases for the interior and interface.
The absence of a separate interface reduced basis is shared with the present work, but the state representations and coupling operators differ.
Here, a shared nonlinear decoder constructs a single, non-overlapping global discrete state, and the subdomains are coupled through the dependence of the original FOM residual on neighboring states.
Thus, we retain the coupling in the existing discrete residual rather than adding a DG interface formulation to connect local representations.
The compatibility obtained by construction is the algebraic property that each global degree of freedom has the same value in every residual that uses it. It is distinct from spatial continuity between different boundary traces.
Table~\ref{tab:coupling_summary} summarizes these coupling approaches and their use of separate interface (port) representations.

\begin{table}[pos=htbp]
\centering
\small
\caption{Summary of the inter-subdomain coupling approaches discussed above.}
\label{tab:coupling_summary}
\begin{tabular}{@{}>{\raggedright\arraybackslash}p{0.27\textwidth}
                  >{\raggedright\arraybackslash}p{\dimexpr0.50\textwidth-4\tabcolsep\relax}
                  >{\raggedright\arraybackslash}p{0.23\textwidth}@{}}
\toprule
Approach & Inter-subdomain coupling & Separate interface representations \\
\midrule
DG-based ROM \citep{chung2024train-073}
& Weak coupling through DG interface operators & No \\
\addlinespace
ROMs with separate interior and port representations \citep{diaz2024fast-e0a,zanardi2024scalable}
& Compatibility constraints between subdomain-specific interface states & Yes \\
\addlinespace
Present work
& Coupling through the original FOM residual using uniquely owned state values & No \\
\bottomrule
\end{tabular}
\end{table}

Locality in the proposed approach is tied to the partition of the physical computational domain.
This differs from local and adaptive ROMs that alleviate the limitations of a single global basis by selecting or interpolating bases according to the parameter regime or region of state space \citep{amsallem2009method-a19,amsallem2011online-88b,amsallem2012nonlinear-21f}.
Unlike neural PDE solvers such as finite-basis PINNs, which use domain-decomposition-inspired network architectures \citep{moseley2023finite-63f}, the present work is an intrusive reduced-order model that directly uses the time-discrete FOM residual and its Jacobian.
After training, the decoder is held fixed, and the subdomain-specific reduced coordinates are determined jointly from a global residual minimization problem at each time step.
The network therefore provides a low-dimensional representation in which to solve the governing equations, rather than directly predicting the state at the next time step.

In this construction, the shared decoder enables reuse of the learned local representation, while non-overlapping ownership and references to the same state remove the need for separate interface compatibility constraints.
Physical coupling between subdomains is retained in the original time-discrete residual, and the local dependence of the residual is used to assemble the Gauss--Newton system from subdomain contributions.
The main contributions of this paper are as follows.
\begin{itemize}
\item We present a bottom-up ROM that combines a shared nonlinear representation of the full state owned by each subdomain with a non-overlapping state construction, ensuring algebraic compatibility by construction without separate interior/interface models or interface copies.
\item We formulate the coupling of local reduced representations through the original time-discrete PDE residual and assemble the global Gauss--Newton system using local Jacobians with respect to neighboring reduced coordinates.
\item We reuse small-domain representations on larger domains without retraining for a viscous Burgers problem with parameterized initial conditions and an incompressible Navier--Stokes problem with parameterized forcing, and compare their accuracy--efficiency trade-offs with POD models using the same residual-coupling structure.
\end{itemize}
The last comparison evaluates nonlinear and linear representations within the proposed coupling structure. It is not a direct performance comparison with existing methods that use separate interior/port representations or DG coupling.

In the following sections, we first express the full-order state and time-discrete residual in a non-overlapping domain-decomposed form, making the dependence of each local residual on its own subdomain state and neighboring states explicit.
We then introduce a nonlinear representation based on a common decoder and substitute it into the residual to derive the reduced residual minimization problem and local assembly procedure.
The numerical experiments evaluate local representation accuracy and then large-domain prediction performance for the two PDEs.
Details of data generation, autoencoder architecture and training, and numerical solvers are provided in the appendices.

\section{Domain-Decomposed Full-Order Residual Formulation}
\label{sec:dd_fom}
Consider the full-order discrete problem at a single time step:
\begin{equation}
 r(x^{(n+1)})=0
\end{equation}
For notational simplicity, we omit the time-step index and the dependence on previous-time states and parameters in what follows, and denote the current state by \(x\in\mathbb R^N\) and the residual by \(r:\mathbb R^N\to\mathbb R^{N_r}\).
Without changing this global problem, we make explicit the local structure needed for subdomain-wise residual evaluation and assembly.

The computational domain is partitioned into \(n_\Omega\) non-overlapping subdomains, and each global degree of freedom is assigned a unique owning subdomain.
Let \(x_i\) denote the state owned by subdomain \(i\). The global state is formed by stacking these blocks: 
\begin{equation}
 x=
 \begin{bmatrix}
 x_1\\x_2\\\vdots\\x_{n_\Omega}
 \end{bmatrix},
 \qquad x_i\in\mathbb R^{N_{\mathrm{loc}}},
 \qquad N=n_\Omega N_{\mathrm{loc}}
 \label{eq:full_state_partition}
\end{equation}
In the scale-up setting considered here, all subdomains within a given experiment have the same number of owned degrees of freedom, \(N_{\mathrm{loc}}\). 
The local grids are aligned with the same coordinate axes, and the degrees of freedom are ordered consistently across subdomains.
The block \(x_i\) contains the full state assigned to that subdomain, rather than two separate states for the interior and interface. 

To define the local residual, we introduce a residual row-selection matrix \(P_i\in\{0,1\}^{m_i\times N_r}\).
The \(m_i\) residual equations associated with subdomain \(i\) are given by
\begin{equation}
 r_i(x):=P_i r(x),\qquad i=1,\ldots,n_\Omega
 \label{eq:local_fom_residual}
\end{equation}
If the selection matrices partition all residual rows without overlap, selecting each row exactly once, then
\begin{equation}
 \sum_{i=1}^{n_\Omega}P_i^T P_i=I,
 \qquad
 r(x)=\sum_{i=1}^{n_\Omega}P_i^T r_i(x)
 \label{eq:assembled_fom_residual}
\end{equation}
Thus, partitioning the state and residual by subdomain can be viewed as a reordering of the original discrete problem.

Because of the locality of the spatial discretization, \(r_i\) generally does not depend on every component of the global state.
Instead, it uses only the state of its own subdomain and the neighboring states required for stencil or numerical flux evaluations near the boundary. 
Let \(\mathcal N(i)\) 
denote the indices of neighboring subdomains whose state values are required to evaluate \(r_i\), and define
\begin{equation}
 \mathcal B_i:=\{i\}\cup\mathcal N(i)
 \label{eq:neighborhood_block}
\end{equation}
Writing \(x_{\mathcal B_i}\) for the vector obtained by stacking these state blocks in a chosen order, we have
\begin{equation}
 r_i(x)=r_i(x_{\mathcal B_i})
 =r_i\bigl(x_i,\{x_j\}_{j\in\mathcal N(i)}\bigr)
 \label{eq:local_residual_dependence}
\end{equation}
Note that, although entire neighboring state blocks appear as arguments, only the components required by the boundary stencil may contribute to the residual.

Using neighboring states is different from duplicating interface states as separate unknowns.
In Eq.~\eqref{eq:local_residual_dependence}, \(x_j\) is the state owned by subdomain \(j\), and every residual that requires its values refers to that same state. Therefore,
interactions with neighboring subdomains are retained in the original residual, without compatibility constraints to match multiple copies of the same degree of freedom.
The next section introduces a nonlinear representation on top of this non-overlapping state structure.

\section{Nonlinear Reduced Representation}
\label{sec:nonlinear_rep}
The key idea is to reuse the same decoder for all subdomains while assigning independent reduced coordinates to each subdomain.
For the state \(x_i\in\mathbb R^{N_{\mathrm{loc}}}\) of subdomain \(i\), we introduce reduced coordinates \(\hat x_i\in\mathbb R^{n_{\mathrm{loc}}}\) and approximate the state as
\begin{equation}
 x_i\approx g(\hat x_i),\qquad i=1,\ldots,n_\Omega
 \label{eq:local_decoder_representation}
\end{equation}
Here, \(g:\mathbb R^{n_{\mathrm{loc}}}\to\mathbb R^{N_{\mathrm{loc}}}\) is a common nonlinear map, and \(n_{\mathrm{loc}}\ll N_{\mathrm{loc}}\) is the reduced dimension per subdomain.
Throughout this paper, we refer to this map as a decoder because it reconstructs a full subdomain state from low-dimensional coordinates.
The formulation requires a differentiable decoder. In the numerical experiments, we use the decoder of an autoencoder trained on local snapshots from a small domain.

Sharing the same function \(g\) does not mean that all subdomains have the same state.
Different inputs \(\hat x_i\) allow the same learned representation to describe different states at different spatial locations.
Moreover, all degrees of freedom owned by a subdomain are represented through a single reconstruction map, without separate models for its interior and interface.
In the bottom-up setting considered here, subdomains of the same size and resolution are repeatedly used, so the map learned on a small domain can be applied unchanged to a larger domain.

We define the global reduced state as
\begin{equation}
 \hat x=
 \begin{bmatrix}
 \hat x_1\\\hat x_2\\\vdots\\\hat x_{n_\Omega}
 \end{bmatrix}
 \in\mathbb R^{n_\Omega n_{\mathrm{loc}}}
 \label{eq:global_reduced_state}
\end{equation}
The corresponding global full-order approximation is
\begin{equation}
 x\approx\tilde x(\hat x)=
 \begin{bmatrix}
 g(\hat x_1)\\g(\hat x_2)\\\vdots\\g(\hat x_{n_\Omega})
 \end{bmatrix}
 \label{eq:global_decoder_representation}
\end{equation}
Because the subdomain states are non-overlapping, the global state is constructed by simply stacking their blocks, without introducing separate interface states.

In Eq.~\eqref{eq:global_decoder_representation}, the value of each global degree of freedom is uniquely determined by the decoder output of its owning subdomain.
Every neighboring residual that uses that degree of freedom refers to the same reconstructed value, so algebraic compatibility for the same degree of freedom holds by construction.
This is a property of the state construction and does not depend on particular reduced-coordinate values or on convergence of the iterative solver.
Whereas the shared decoder enables reuse of the learned representation, non-overlapping ownership and references to the same state ensure this compatibility.
Here, compatibility means that each global degree of freedom has a single assigned value, used consistently by every subdomain residual that depends on it. It does not mean that different boundary traces are spatially continuous for every choice of reduced coordinates.
External boundary conditions and physical coupling are handled through the discrete residual of the original problem.

The global state is thus determined by the fixed common decoder and the collection of subdomain-specific reduced coordinates.
The number of learned decoder parameters does not increase as the domain grows, although the number of global reduced coordinates increases to \(n_\Omega n_{\mathrm{loc}}\).
The next section substitutes this representation into the local full-order residuals to construct a reduced residual minimization problem that jointly determines the reduced coordinates on the larger domain.

\section{Domain-Decomposed Reduced System and Local Assembly}
\label{sec:dd_reduced_system}
The previous section represented each subdomain state using a common decoder.
We now substitute this representation into the residual in Eq.~\eqref{eq:local_residual_dependence} to define a residual in reduced coordinates and derive the Gauss--Newton system used to minimize it.

\subsection{Reduced-coordinate residual and least-squares problem}
Approximating each state \(x_j\), \(j\in\mathcal B_i\), involved in the residual of subdomain \(i\) by \(g(\hat x_j)\) gives
\begin{equation}
 \tilde r_i(\hat x_{\mathcal B_i})
 :=r_i\!\left(\{g(\hat x_j)\}_{j\in\mathcal B_i}\right)
 =r_i\!\left(g(\hat x_i),\{g(\hat x_j)\}_{j\in\mathcal N(i)}\right)
 \label{eq:local_reduced_coordinate_residual}
\end{equation}
Here, \(\hat x_{\mathcal B_i}\) is the vector of reduced coordinates for the subdomain and its required neighbors.
This expression shows that the local residual depends not only on the subdomain's own coordinates but also on neighboring coordinates.
For notational simplicity, we also write this residual as \(\tilde r_i(\hat x)\), although its actual dependence is restricted to the blocks in \(\mathcal B_i\).

At each time step, the global reduced coordinates are determined by minimizing the sum of squared subdomain residuals:
\begin{equation}
 \min_{\hat x\in\mathbb R^{n_\Omega n_{\mathrm{loc}}}}
 \mathcal F_{\mathrm{PDE}}(\hat x),
 \qquad
 \mathcal F_{\mathrm{PDE}}(\hat x)
 :=\frac12\sum_{i=1}^{n_\Omega}\|\tilde r_i(\hat x)\|_2^2.
 \label{eq:dd_reduced_ls}
\end{equation}
All local residuals are selected from the same globally reconstructed state:
\(\tilde r_i(\hat x)=P_i r(\tilde x(\hat x))\).
It follows from Eq.~\eqref{eq:assembled_fom_residual} that
\begin{equation}
 \begin{aligned}
 \sum_{i=1}^{n_\Omega}\|\tilde r_i(\hat x)\|_2^2
 &=r(\tilde x(\hat x))^T
   \left(\sum_{i=1}^{n_\Omega}P_i^TP_i\right)r(\tilde x(\hat x))\\
 &=\|r(\tilde x(\hat x))\|_2^2
 \end{aligned}
 \label{eq:global_lspg_equivalence}
\end{equation}
Thus, the sum of the subdomain objectives is equivalent to minimizing the original global time-discrete residual on the concatenated nonlinear manifold.
This is a least-squares Petrov--Galerkin (LSPG) formulation \citep{https://doi.org/10.1002/nme.3050,lee2020model-ab0}: rather than reducing the entire space--time problem at once, we apply the minimum-residual principle at each time step \citep{doi:10113717M1120531}.

The subdomain-specific reduced coordinates are therefore determined jointly through residuals that depend on neighboring states, rather than integrated independently.
The coupling through numerical fluxes and stabilization terms in the original discretization remains intact, even without additional interface compatibility constraints. 
Moreover, minimizing the residual within a given manifold is distinct from making the governing residual exactly zero.

\subsection{Gauss--Newton reduced system}
\label{sec:gn_derivation}
We apply a Gauss--Newton linearization at the current reduced state \(\hat x\) to solve Eq.~\eqref{eq:dd_reduced_ls}.
First, define the residual Jacobian for subdomain \(i\) as
\begin{equation}
 J_i(\hat x):=\frac{\partial\tilde r_i(\hat x)}{\partial\hat x}
 \in\mathbb R^{m_i\times n_\Omega n_{\mathrm{loc}}}
\end{equation}
Although this Jacobian is taken with respect to the global reduced coordinates, \(\tilde r_i\) depends only on coordinates in \(\mathcal B_i\), so nonzero entries are confined to the corresponding block columns.
Letting \(\Delta\hat x\) denote the increment at the current iterate, we have
\begin{equation}
 \tilde r_i(\hat x+\Delta\hat x)
 \approx\tilde r_i(\hat x)+J_i(\hat x)\Delta\hat x
\end{equation}
and substituting this approximation into the objective yields the linear least-squares problem
\begin{equation}
 \min_{\Delta\hat x}\frac12\sum_{i=1}^{n_\Omega}
 \|\tilde r_i(\hat x)+J_i(\hat x)\Delta\hat x\|_2^2
\end{equation}

The corresponding normal equations are
\begin{equation}
 \left(\sum_{i=1}^{n_\Omega}J_i^TJ_i\right)\Delta\hat x
 =-\sum_{i=1}^{n_\Omega}J_i^T\tilde r_i
 \label{eq:gn_normal_equations}
\end{equation}
For notational simplicity, we suppress the dependence on the current state in Gauss–Newton iterate and define the subdomain contributions as follows.
\begin{equation}
 H_i:=J_i^TJ_i,\qquad q_i:=J_i^T\tilde r_i
\end{equation}
we can write the system as
\begin{equation}
 H_{\mathrm{PDE}}\Delta\hat x=-q_{\mathrm{PDE}},
 \qquad
 H_{\mathrm{PDE}}=\sum_iH_i,\qquad q_{\mathrm{PDE}}=\sum_iq_i
 \label{eq:assembled_reduced_system}
\end{equation}
Thus, although the nonlinear representation is defined subdomain-wise, the final reduced system is a globally coupled system assembled from interactions between each subdomain and its neighbors.

In the notation below, \(\operatorname{col}\) denotes vertical concatenation, in the specified order, of vectors or matrices with the same number of columns.
Writing the global reduced residual Jacobian as \(J=\operatorname{col}_iJ_i\), we have
\begin{equation}
 H_{\mathrm{PDE}}=J^TJ,
 \qquad
 v^TH_{\mathrm{PDE}}v=\|Jv\|_2^2\ge0
\end{equation}
so the Gauss--Newton matrix is symmetric positive semidefinite.
If \(J\) has linearly independent columns, the matrix is positive definite, and preconditioned conjugate gradients (PCG) can be applied with a symmetric positive-definite preconditioner. 
If \(H_{\mathrm{PDE}}\) is only positive semidefinite, a positive diagonal shift \(\epsilon I\), with \(\epsilon>0\), can be used to obtain a positive-definite system matrix for PCG. No diagonal regularization was used in the reported experiments.
The subdomain-wise block-Jacobi preconditioner and stopping criteria used in the experiments are described in the respective solver settings and in Appendix~\ref{app:rom_solver}.

\subsection{Block-sparse local assembly}
\label{subsec:block_sparse_local_assembly}
The system in Eq.~\eqref{eq:assembled_reduced_system} can be assembled without first constructing the global reduced Jacobian.
This follows from the restricted dependence of each residual on neighboring states, not from the use of a common decoder itself.
To make this structure explicit, introduce a 
selection matrix \(S_i\in\{0,1\}^{(|\mathcal B_i|n_{\mathrm{loc}})\times(n_\Omega n_{\mathrm{loc}})}\) that selects the reduced blocks in \(\mathcal B_i\).
Here, \(|\mathcal B_i|\) counts subdomain \(i\) itself and the neighboring subdomains whose states are required to evaluate \(r_i\).
We then write
\begin{equation}
 \hat x_{\mathcal B_i}=S_i\hat x
\end{equation}
Define the Jacobian with respect to these local coordinates as
\begin{equation}
 J_i^{\mathrm{loc}}:=\frac{\partial\tilde r_i}{\partial\hat x_{\mathcal B_i}}
\end{equation}
It is related to the Jacobian with respect to global coordinates by
\begin{equation}
 J_i=J_i^{\mathrm{loc}}S_i
 \label{eq:local_global_jacobian}
\end{equation}
The derivatives that must actually be computed are therefore confined to small local Jacobians with respect to neighboring reduced coordinates.

The local Jacobian is evaluated by applying the chain rule to the full-order residual and the decoder. Define
\begin{equation}
 R_i^{\mathrm{loc}}:=\frac{\partial r_i}{\partial x_{\mathcal B_i}},
 \qquad
 G_j(\hat x_j):=\frac{\partial g(\hat x_j)}{\partial\hat x_j},
 \qquad j\in\mathcal B_i
\end{equation}
and collect the decoder Jacobians in block-diagonal form as
\begin{equation}
 G_i^{\mathrm{loc}}:=\operatorname{blkdiag}
 \bigl(G_j(\hat x_j)\bigr)_{j\in\mathcal B_i}
\end{equation}
Then the chain rule gives
\begin{equation}
 J_i^{\mathrm{loc}}=R_i^{\mathrm{loc}}G_i^{\mathrm{loc}}
 \label{eq:local_chain_rule}
\end{equation}
This requires only the local full-order residual Jacobian of subdomain \(i\) and the decoder Jacobians of the subdomains needed by that residual.

Define the local Gauss--Newton contributions as
\begin{equation}
 H_i^{\mathrm{loc}}:=(J_i^{\mathrm{loc}})^TJ_i^{\mathrm{loc}},
 \qquad
 q_i^{\mathrm{loc}}:=(J_i^{\mathrm{loc}})^T\tilde r_i
\end{equation}
The global matrix and vector can then be assembled as
\begin{equation}
 H_{\mathrm{PDE}}=\sum_{i=1}^{n_\Omega}S_i^TH_i^{\mathrm{loc}}S_i,
 \qquad
 q_{\mathrm{PDE}}=\sum_{i=1}^{n_\Omega}S_i^Tq_i^{\mathrm{loc}}
 \label{eq:scatter_assembly}
\end{equation}
This is a local-to-global scatter operation. Each local matrix and vector is computed and added to the corresponding global block locations.
Nevertheless, the computation of each contribution remains confined to the local neighborhood of that residual.

The assembly above accounts for the PDE residual contribution.
In the Navier--Stokes experiments, an additional term controls the global mean pressure. Its objective and global correction are described in Section~\ref{sec:ns_pressure_gauge} and Appendix~\ref{app:pressure_penalty}.
Local assembly does not, by itself, hyper-reduce nonlinear residual evaluation.
We evaluate the full residual in this work and leave the combination of the local residual structure with residual sampling or solution-based nonlinear subspaces for future work \citep{choi2020sns:-40d}. 

\section{Numerical Experiments}
\label{sec:numerical_experiments}
This section evaluates the proposed method on the two-dimensional viscous Burgers equations and the incompressible Navier--Stokes equations.
The two experiments have the same goals.
First, we examine how compactly a local nonlinear representation learned on a small \(3\times3\) subdomain configuration can represent snapshots not used for training.
We then hold the trained decoder fixed and apply it to the 225 subdomains of a \(15\times15\) configuration to assess whether time-dependent solutions on the larger domain can be computed without retraining.
Comparison with POD-based linear models allows us to examine both the compactness of the local representation and the accuracy--efficiency trade-off of the online computation.

\subsection{Common reduced models and evaluation setup}
\label{sec:common_evaluation}
\paragraph{Field-wise representations and the shared decoder.}
In both experiments, a separate autoencoder is trained on the local snapshots of each physical field.
The field sets are \(\mathcal Q_{\mathrm B}=\{u,v\}\) for Burgers and \(\mathcal Q_{\mathrm{NS}}=\{u,v,p\}\) for Navier--Stokes.
Denoting the encoder and decoder for field \(q\) by
\begin{equation}
 h_q:\mathbb R^m\to\mathbb R^{n_q},
 \qquad
 g_q:\mathbb R^{n_q}\to\mathbb R^m
 \label{eq:fieldwise_maps}
\end{equation}
the reconstruction of a local snapshot \(x_q\) is given by
\begin{equation}
 z_q=h_q(x_q),\qquad
 \tilde x_q=g_q(z_q)=\mathcal A_q(x_q),\qquad
 \mathcal A_q:=g_q\circ h_q
 \label{eq:experiment_field_sets}
\end{equation}
Here, \(m\) is the local number of degrees of freedom per field for the problem considered, and \(n_q\) is the field-wise reduced dimension.

These field-wise representations are related to the common decoder \(g\) in Section~\ref{sec:nonlinear_rep} through
\begin{equation}
 x_i=\operatorname{col}_{q\in\mathcal Q}x_{q,i},\qquad
 \hat x_i=\operatorname{col}_{q\in\mathcal Q}z_{q,i},\qquad
 g(\hat x_i)=\operatorname{col}_{q\in\mathcal Q}g_q(z_{q,i})
 \label{eq:fieldwise_block_decoder}
\end{equation}
The field order is \((u,v)\) for Burgers and \((u,v,p)\) for Navier--Stokes.
Thus, \(N_{\mathrm{loc}}=|\mathcal Q|m\), and the total reduced dimension per subdomain is
\begin{equation}
 n_{\mathrm{loc}}=\sum_{q\in\mathcal Q}n_q
 \label{eq:total_local_dimension}
\end{equation}
The shared decoder refers to the state-reconstruction function obtained by concatenating these field-wise maps. It does not mean that all fields are trained simultaneously in a single network.
Each field-wise decoder is reused on all subdomains within an experiment, without separate interior and interface representations.
The trained models are not shared between the Burgers and Navier--Stokes problems.

Each encoder--decoder pair is trained on local snapshots from the training split by minimizing the reconstruction loss:
\begin{equation}
 \min_{\theta_{h,q},\theta_{g,q}}
 \frac{1}{mN_{\mathrm{tr}}}\sum_{s=1}^{N_{\mathrm{tr}}}
 \left\|g_q\bigl(h_q(x_q^{[s]})\bigr)-x_q^{[s]}\right\|_2^2.
 \label{eq:fieldwise_ae_training_loss}
\end{equation}
Here, \(x_q^{[s]}\) denotes the \(s\)-th training sample of field \(q\), pooled over subdomains, time instances, and training cases, and \(N_{\mathrm{tr}}\) is the number of training samples per field.
After training, the decoder parameters are fixed, and the reduced coordinates are treated directly as unknowns in the online residual minimization problem.
The encoder is used for representation learning and local reconstruction evaluation, and can also be used to estimate the initial reduced coordinates.
Details of the network architecture, normalization, sparse output layer, and training settings are provided in Appendix~\ref{sec:model_architecture}.

The POD baseline uses an orthonormal basis \(\Phi_q\) constructed from the same local training data and has the form 
\begin{equation}
 x_q\approx\bar x_q+\Phi_q a_q
 \label{eq:pod_local_model}
\end{equation}
Here, \(\bar x_q\) is the reference state and \(a_q\) contains the reduced coordinates.
We denote the AE-based model by DD-NM-ROM and the POD-based model by DD-LS-ROM.
Both models use the same non-overlapping state construction and residual coupling, so the comparisons below assess nonlinear and linear representations within the same reduced-system construction.

\paragraph{Local representation accuracy.}
Before time integration on the larger domain, we evaluate representation accuracy using local snapshots extracted from separate small-domain test trajectories.
For AE models, we use the encoder--decoder reconstruction \(\tilde x_q^{\mathrm{AE}}=g_q(h_q(x_q))\). 
For POD models, we use the orthogonal projection \(\tilde x_q^{\mathrm{POD}}=\bar x_q+\Phi_q\Phi_q^T(x_q-\bar x_q)\). 
The AE reconstruction is not assumed to be an orthogonal projection obtained by directly finding the nearest point on the nonlinear manifold. 
We use the term representation error to refer to the two cases together.

Let \(X_q\in\mathbb R^{m\times N_s}\) contain the \(N_s\) test snapshots of field \(q\) as columns, and let \(\tilde X_q\) be the corresponding reconstructed or projected snapshot matrix.
For both problems, we use the relative Frobenius error
\begin{equation}
 \varepsilon_{\mathrm{rep},q}
 =\frac{\|\tilde X_q-X_q\|_F}{\|X_q\|_F}
 \label{eq:representation_error}
\end{equation}
This evaluation measures the compactness and reconstruction capability of the learned representation and is distinct from the ROM prediction error obtained by solving the residual problem on the larger domain.

\paragraph{Large-domain prediction accuracy and computational time.}
Let \(t_\ell\), \(\ell=1,\ldots,N_{\mathrm{snap}}\), denote the saved time instances in the scale-up experiment, and let \(q_{\mathrm{ROM}}^\ell\) and \(q_{\mathrm{FOM}}^\ell\) be the corresponding spatial state vectors.
For a single test case, the field-wise space--time relative error is defined as
\begin{equation}
 E_q^{L_2}
 =\frac{\left(\sum_{\ell=1}^{N_{\mathrm{snap}}}
 \|q_{\mathrm{ROM}}^\ell-q_{\mathrm{FOM}}^\ell\|_2^2\right)^{1/2}}
 {\left(\sum_{\ell=1}^{N_{\mathrm{snap}}}
 \|q_{\mathrm{FOM}}^\ell\|_2^2\right)^{1/2}}
 \label{eq:field_space_time_error}
\end{equation}
This is the relative Frobenius error of the snapshot matrix for the complete space--time trajectory and measures the error accumulated over all saved time instances.
To assess the largest instantaneous deviation, we also use the global-normalized maximum error:
\begin{equation}
 E_q^{\mathrm{Gmax}}
 =\frac{\max_{1\le\ell\le N_{\mathrm{snap}}}
 \|q_{\mathrm{ROM}}^\ell-q_{\mathrm{FOM}}^\ell\|_2}
 {\max_{1\le\ell\le N_{\mathrm{snap}}}
 \|q_{\mathrm{FOM}}^\ell\|_2}.
 \label{eq:field_global_max_error}
\end{equation}
The maxima in the numerator and denominator are taken independently over all saved time instances.
Thus, this metric is the maximum deviation normalized by a single global trajectory scale, not the maximum of the instantaneous relative errors.

The all-field metrics \(E_{\mathrm{all}}^{L_2}\) and \(E_{\mathrm{all}}^{\mathrm{Gmax}}\) are computed in the same way by concatenating \([u;v]\) for Burgers and \([u;v;p]\) for Navier--Stokes.
The all-field error is not a simple average of the field-wise relative errors and depends on the magnitude of each field.
This relationship and the interpretation of the two error metrics are discussed further in Appendix~\ref{app:metric_interpretation}.

Letting \(T_{\mathrm{FOM}}\) and \(T_{\mathrm{ROM}}\) denote the FOM and ROM wall-clock times, respectively, we define the speed-up as
\begin{equation}
 S=\frac{T_{\mathrm{FOM}}}{T_{\mathrm{ROM}}}
 \label{eq:speedup}
\end{equation}
All timings are measured in the same single-threaded environment on an Apple M3 Max with 36 GB of unified memory.
The scale-up experiment for each problem uses 10 distinct test cases.
Errors and speed-ups are computed for each case and then arithmetically averaged for the tables and results discussion; representative state plots show one of these cases.
The timing scope for each problem is described in Appendix~\ref{app:timing_details}.

\subsection{Two-dimensional viscous Burgers equations}
\label{sec:burgers_experiment}
\subsubsection{Problem setup}
The first example is the two-dimensional viscous Burgers problem for the velocity field \(\mathbf u=(u,v)\):
\begin{equation}
 \frac{\partial\mathbf u}{\partial t}
 +(\mathbf u\cdot\nabla)\mathbf u
 =\nu\nabla^2\mathbf u.
 \label{eq:burgers_momentum}
\end{equation}
The viscosity is fixed at \(\nu=5\times10^{-3}\), and the equations are integrated with \(\Delta t=0.02\) to the final time \(T=1\).
To generate diverse initial states, the field-wise initial conditions on each subdomain \(\Omega_i\) are defined as
\begin{equation}
 q(x,y,0;\boldsymbol\mu)
 =\mu_i^q\sin(2\pi x)\sin(2\pi y),
 \qquad q\in\{u,v\},\quad (x,y)\in\Omega_i
 \label{eq:burgers_ic}
\end{equation}
Here, \(x,y\) are global physical coordinates, and the amplitude, sign, and activation of the coefficient \(\mu_i^q\) are sampled independently for each subdomain and field.
The coefficient distributions and an example of the initial state are provided in Appendix~\ref{app:burgers_data_generation}.

\subsubsection{Training and scale-up configurations}
Training data are generated on a \(3\times3\) subdomain configuration over \([0,1.5]^2\).
Each subdomain has size \(0.5\times0.5\) and uses a \(50\times50\) cell-centered grid, giving a global grid of \(150\times150\) and \(m=2500\) local degrees of freedom per field.
Periodic boundary conditions are imposed in both directions on the small domain.
Independent trajectories are generated by initial-condition case: 1600 for training, 400 for validation, and 200 for local reconstruction testing, with 51 snapshots stored per trajectory, including the initial state.
All time and subdomain samples from a given case belong to the same data split.

\begin{figure}[pos=htbp]
\centering
\includegraphics[width=0.48\textwidth]{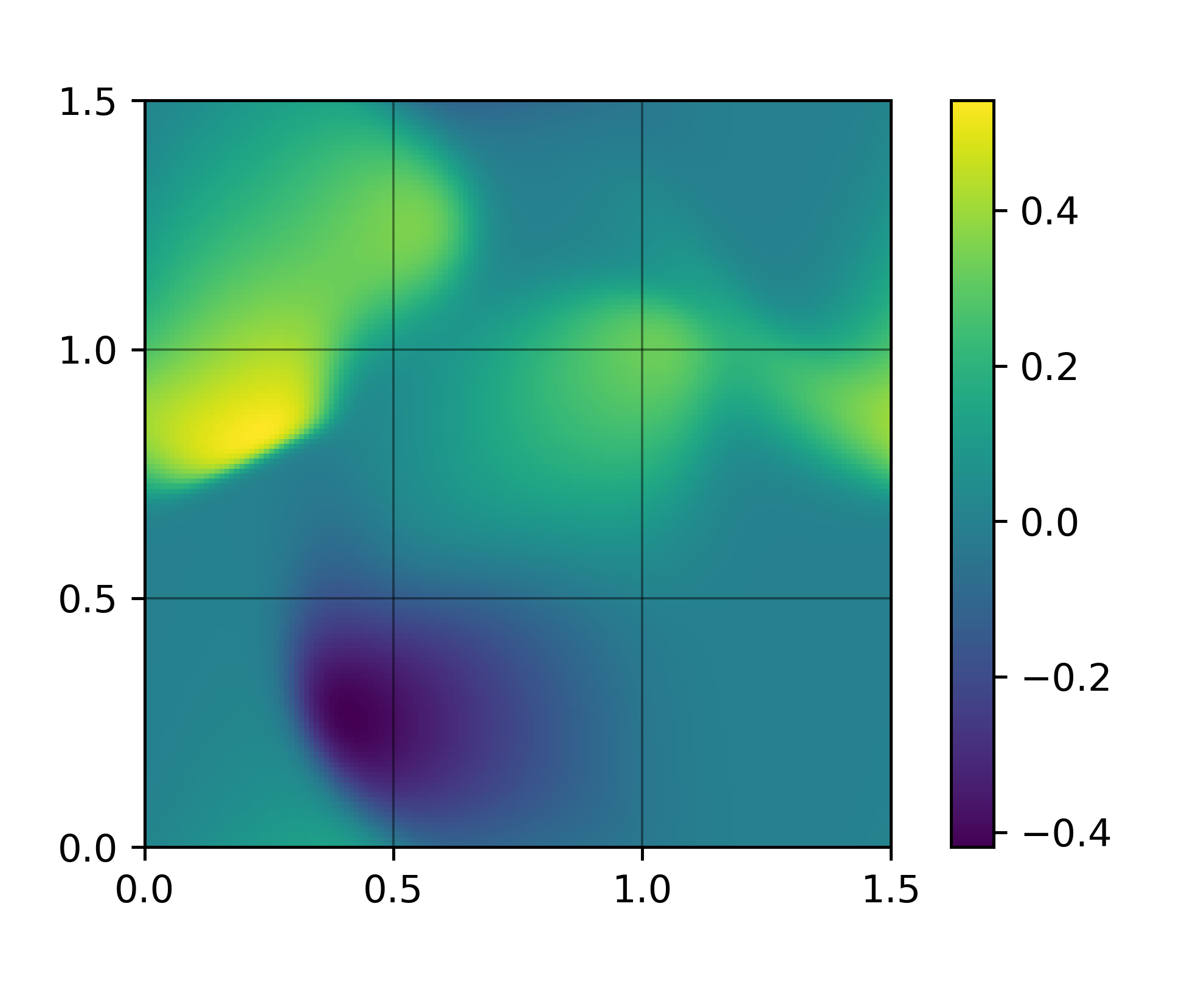}\hfill
\includegraphics[width=0.48\textwidth]{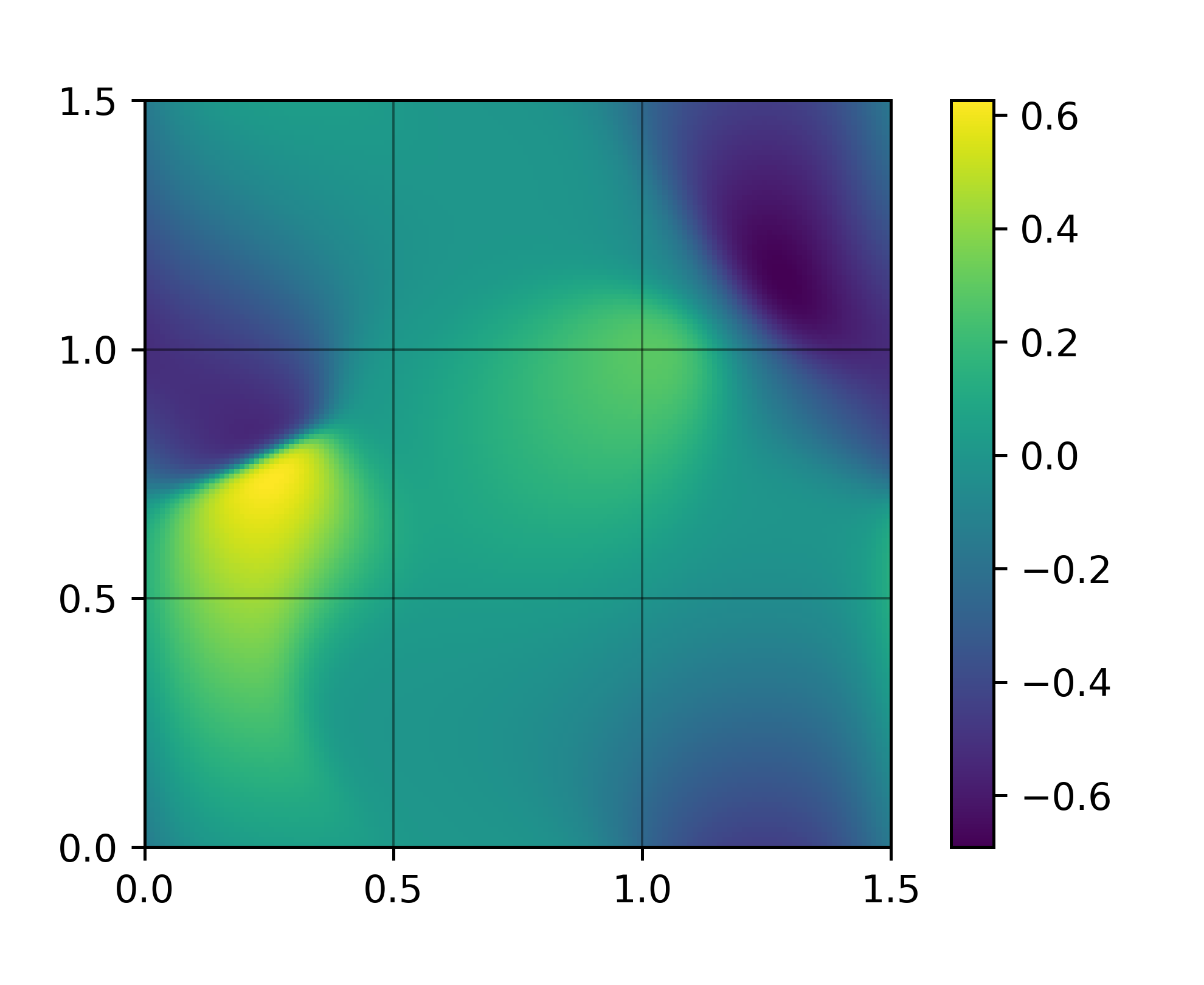}
\caption{Representative FOM training snapshot for the \(3\times3\) Burgers configuration at \(t=1\).
The left and right panels show \(u\) and \(v\), respectively.
Each \(50\times50\) subdomain block provides one local snapshot per field.
These local snapshots are pooled over subdomains, time instances, and training cases to construct the AE and POD representations.}
\label{fig:burgers_training_data_example}
\end{figure}

Nine local states per field are extracted from each global snapshot and pooled over time and cases to train the AEs and construct the POD bases.
Figure~\ref{fig:burgers_training_data_example} illustrates the relationship between a global training snapshot and its local samples.
The field-wise sample counts for training, validation, and testing, together with the snapshot collection procedure, are described in Appendices~\ref{app:data_split} and~\ref{sec:burgers_data_split}.

For the scale-up evaluation, the domain is enlarged to \([0,7.5]^2\) with a \(750\times750\) global grid, while the local grid spacing is kept unchanged.
The number of subdomains is then \(15\times15=225\), and homogeneous Neumann conditions are imposed on the outer boundary.
Thus, this experiment changes not only the number of subdomains but also the external boundary conditions relative to the small training domain.
The evaluation uses 10 separately generated large-domain initial-condition cases, with the local AEs and POD bases from the small domain applied unchanged and without retraining.

\subsubsection{Reduced models and solver settings}
Table~\ref{tab:burgers_reduced_model_configs} summarizes the reduced dimensions compared in the Burgers experiments.
The two velocity components use the same dimension, \(n_u=n_v=n_q\), so the total dimension per subdomain is \(n_{\mathrm{loc}}=2n_q\).
The AEs use the seven specified field-wise dimensions from 20 to 50.
POD dimensions are selected using cumulative energy-retention levels of \(99\%\), \(99.9\%\), \(99.95\%\), and \(99.99\%\), with additional field-wise dimensions of 53 and 61.

\begin{table}[pos=htbp]
\centering
\small
\caption{Reduced model configurations for the Burgers experiments.
The two velocity components use the same dimension, \(n_u=n_v=n_q\).
The total local dimension is the sum of the two component dimensions on one subdomain.}
\label{tab:burgers_reduced_model_configs}
\begin{tabular}{lcc}
\toprule
Model & Field-wise dimension \(n_q\) & Total local dimension \\
\midrule
AE & \(\{20,25,30,35,40,45,50\}\) & \(2n_q\) \\
POD & \(\{16,49,53,61,65,120\}\) & \(2n_q\) \\
\bottomrule
\end{tabular}
\end{table}

The Burgers autoencoder uses an encoder with hidden widths \((300,200,100)\) and a decoder with one hidden layer and a sparse output layer.
The encoder widths and decoder-mask parameters are fixed across latent dimensions, and no data normalization is applied.
The detailed architecture and training settings are provided in Appendices~\ref{app:architecture_settings} and~\ref{app:training_hparams}.

The FOM uses fully implicit backward Euler time integration and a cell-centered finite volume spatial discretization.
Both velocity components are solved jointly by Newton's method. The linear systems are solved by a sparse direct method on the small training domain and by ILU-preconditioned GMRES for the large-domain reference FOM.
The ROM uses Gauss--Newton iterations with subdomain-wise block-Jacobi-preconditioned PCG.
The FOM discretization and linear solver settings are given in Appendices~\ref{app:burgers_spatial_discretization} and~\ref{app:burgers_fom_solver}, and the ROM stopping criteria are given in Appendix~\ref{app:burgers_rom_solver}.

\subsubsection{Representation accuracy on held-out local snapshots}
\label{sec:burgers_representation}
Before evaluating the reduced dynamics on the larger domain, we compare representation accuracy using local test snapshots extracted from 200 separate small-domain cases.
The AE values in Table~\ref{tab:burgers_representation_accuracy} are reconstruction errors induced by the learned encoder--decoder pairs, whereas the POD values are orthogonal projection errors onto the corresponding bases.
All values are computed using the relative Frobenius error in Eq.~\eqref{eq:representation_error}.

\begin{table}[pos=htbp]
\centering
\small
\caption{Field-wise representation errors on held-out local Burgers snapshots, measured by the relative Frobenius error in Eq.~\eqref{eq:representation_error}.
AE values are encoder--decoder reconstruction errors, and POD values are orthogonal projection errors.
The dimension \(n_q\) is specified per velocity component. The reported errors are dimensionless, not percentages.}
\label{tab:burgers_representation_accuracy}
\begin{tabular}{lccc}
\toprule
Model & \(n_q\) & \(\varepsilon_{\mathrm{rep},u}\) & \(\varepsilon_{\mathrm{rep},v}\) \\
\midrule
AE & 20 & \(3.189\times10^{-2}\) & \(3.111\times10^{-2}\) \\
AE & 25 & \(2.230\times10^{-2}\) & \(2.328\times10^{-2}\) \\
AE & 30 & \(1.809\times10^{-2}\) & \(1.906\times10^{-2}\) \\
AE & 35 & \(1.485\times10^{-2}\) & \(1.548\times10^{-2}\) \\
AE & 40 & \(1.252\times10^{-2}\) & \(1.294\times10^{-2}\) \\
AE & 45 & \(1.072\times10^{-2}\) & \(1.106\times10^{-2}\) \\
AE & 50 & \(9.596\times10^{-3}\) & \(9.935\times10^{-3}\) \\
\midrule
POD & 16  & \(9.848\times10^{-2}\) & \(9.819\times10^{-2}\) \\
POD & 49  & \(3.039\times10^{-2}\) & \(3.081\times10^{-2}\) \\
POD & 53  & \(2.774\times10^{-2}\) & \(2.812\times10^{-2}\) \\
POD & 61  & \(2.327\times10^{-2}\) & \(2.373\times10^{-2}\) \\
POD & 65  & \(2.156\times10^{-2}\) & \(2.203\times10^{-2}\) \\
POD & 120 & \(9.420\times10^{-3}\) & \(9.842\times10^{-3}\) \\
\bottomrule
\end{tabular}
\end{table}

Increasing the AE dimension reduces the reconstruction errors for both velocity components.
At a field-wise dimension of 50, the errors in \(u\) and \(v\) are \(9.596\times10^{-3}\) and \(9.935\times10^{-3}\), respectively, or approximately \(1\%\).
These values are close to the POD projection errors of \(9.420\times10^{-3}\) and \(9.842\times10^{-3}\) at a field-wise dimension of 120.
Thus, at this error level, the nonlinear representation reconstructs the test snapshots with a smaller local dimension.
Local snapshot reconstruction accuracy alone does not determine large-domain time-integration performance, however, we next compare predictions obtained by using the same representations in the residual-minimization problem.

\subsubsection{Scale-up prediction on a \(15\times15\) configuration}
\label{sec:burgers_scaleup_results}
Table~\ref{tab:burgers_scaleup_results} reports errors in the full velocity state \([u;v]\) and speed-ups over the FOM for the 10 large-domain test cases.
All models reuse the local representations constructed from the \(3\times3\) training data, and the reported errors and speed-ups are arithmetic means of the case-wise metrics.

\begin{table}[pos=htbp]
\centering
\small
\caption{Scale-up performance on the \(15\times15\) Burgers configuration.
Errors and speed-ups are arithmetic means over 10 test cases.
Field-wise dimensions are listed in the order \((u,v)\), and \(n_{\mathrm{loc}}\) is their sum on one subdomain.
Both error columns are percentages. \(S\) is the speed-up over the FOM.
All local representations are constructed from the \(3\times3\) training data and reused without retraining.}
\label{tab:burgers_scaleup_results}
\small
\setlength{\tabcolsep}{5pt}
\begin{tabular}{lrrrrr}
\toprule
Model & \((u,v)\) dimensions & \(n_{\mathrm{loc}}\) & \(E_{\mathrm{all}}^{L_2}\)(\%) & \(E_{\mathrm{all}}^{\mathrm{Gmax}}\)(\%) & \(S\) \\
\midrule
DD-NM-ROM & \((20,20)\) & 40  & 7.249 & 6.927 & 33.12 \\
DD-NM-ROM & \((25,25)\) & 50  & 5.436 & 5.103 & 26.71 \\
DD-NM-ROM & \((30,30)\) & 60  & 4.239 & 3.916 & 21.42 \\
DD-NM-ROM & \((35,35)\) & 70  & 3.317 & 3.067 & 18.01 \\
DD-NM-ROM & \((40,40)\) & 80  & 2.752 & 2.576 & 16.89 \\
DD-NM-ROM & \((45,45)\) & 90  & 2.370 & 2.263 & 13.40 \\
DD-NM-ROM & \((50,50)\) & 100 & 2.023 & 1.939 & 12.33 \\
\midrule
DD-LS-ROM & \((16,16)\)   & 32  & 13.849 & 12.058 & 61.24 \\
DD-LS-ROM & \((49,49)\)   & 98  & 4.558  & 4.303  & 19.58 \\
DD-LS-ROM & \((53,53)\)   & 106 & 4.236  & 4.033  & 16.38 \\
DD-LS-ROM & \((61,61)\)   & 122 & 3.879  & 3.737  & 12.71 \\
DD-LS-ROM & \((65,65)\)   & 130 & 3.716  & 3.592  & 11.32 \\
DD-LS-ROM & \((120,120)\) & 240 & 2.538  & 2.575  & 4.21 \\
\bottomrule
\end{tabular}
\end{table}

As the AE field-wise dimension increases from 20 to 50, the mean \(E_{\mathrm{all}}^{L_2}\) decreases from \(7.249\%\) to \(2.023\%\), and the mean \(E_{\mathrm{all}}^{\mathrm{Gmax}}\) decreases from \(6.927\%\) to \(1.939\%\).
The mean speed-up correspondingly decreases from \(33.12\) to \(12.33\).
Increasing the POD dimension also improves accuracy, but incurs the computational cost of a larger reduced system.
It is therefore necessary to compare accuracy and computational speed together, rather than considering each model's error alone.

For representative state comparisons, we use AE \((50,50)\) and POD \((65,65)\), which have similar mean speed-ups.
The AE has a total local dimension of 100 and a mean speed-up of \(12.33\), whereas the POD model has a total local dimension of 130 and a speed-up of \(11.32\).
Their mean space--time errors are \(2.023\%\) and \(3.716\%\), respectively, and their mean global-normalized maximum errors are \(1.939\%\) and \(3.592\%\).
These two models therefore provide a comparison of prediction accuracy at similar computational speeds.

\begin{figure}[pos=htbp]
\centering
\includegraphics[width=0.95\textwidth]{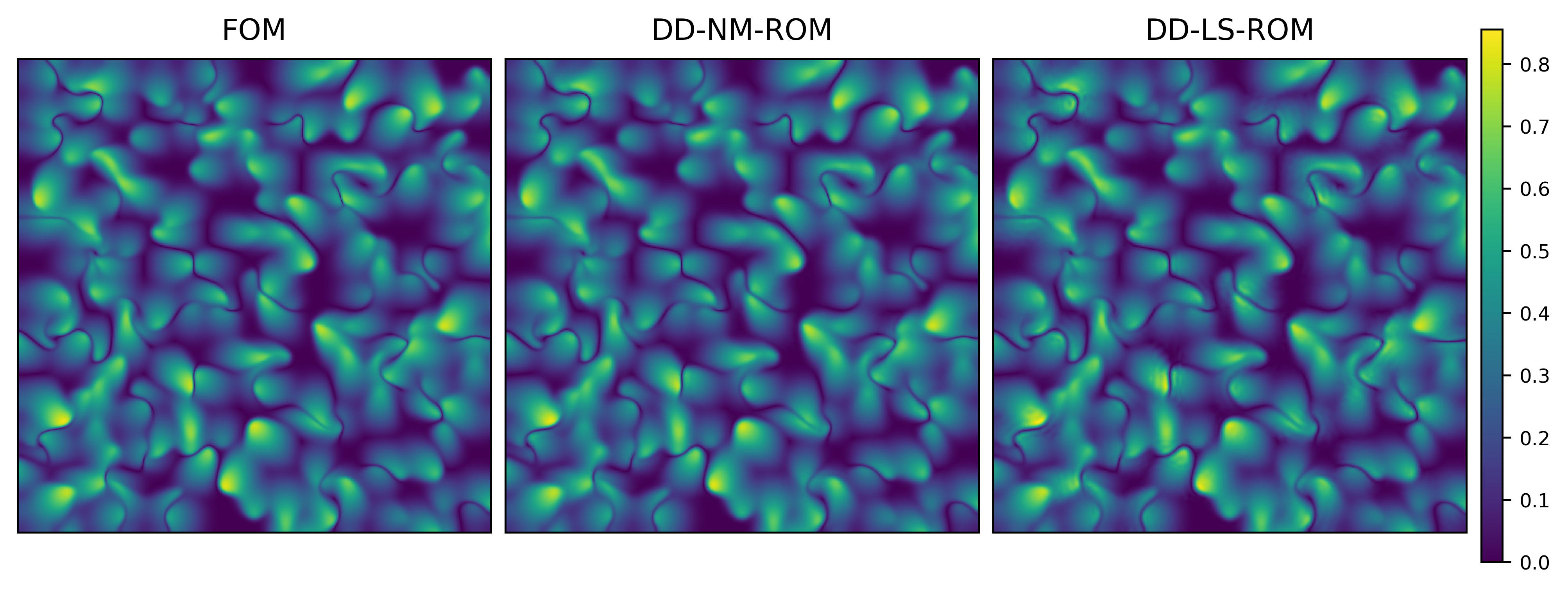}\\
\includegraphics[width=0.69\textwidth]{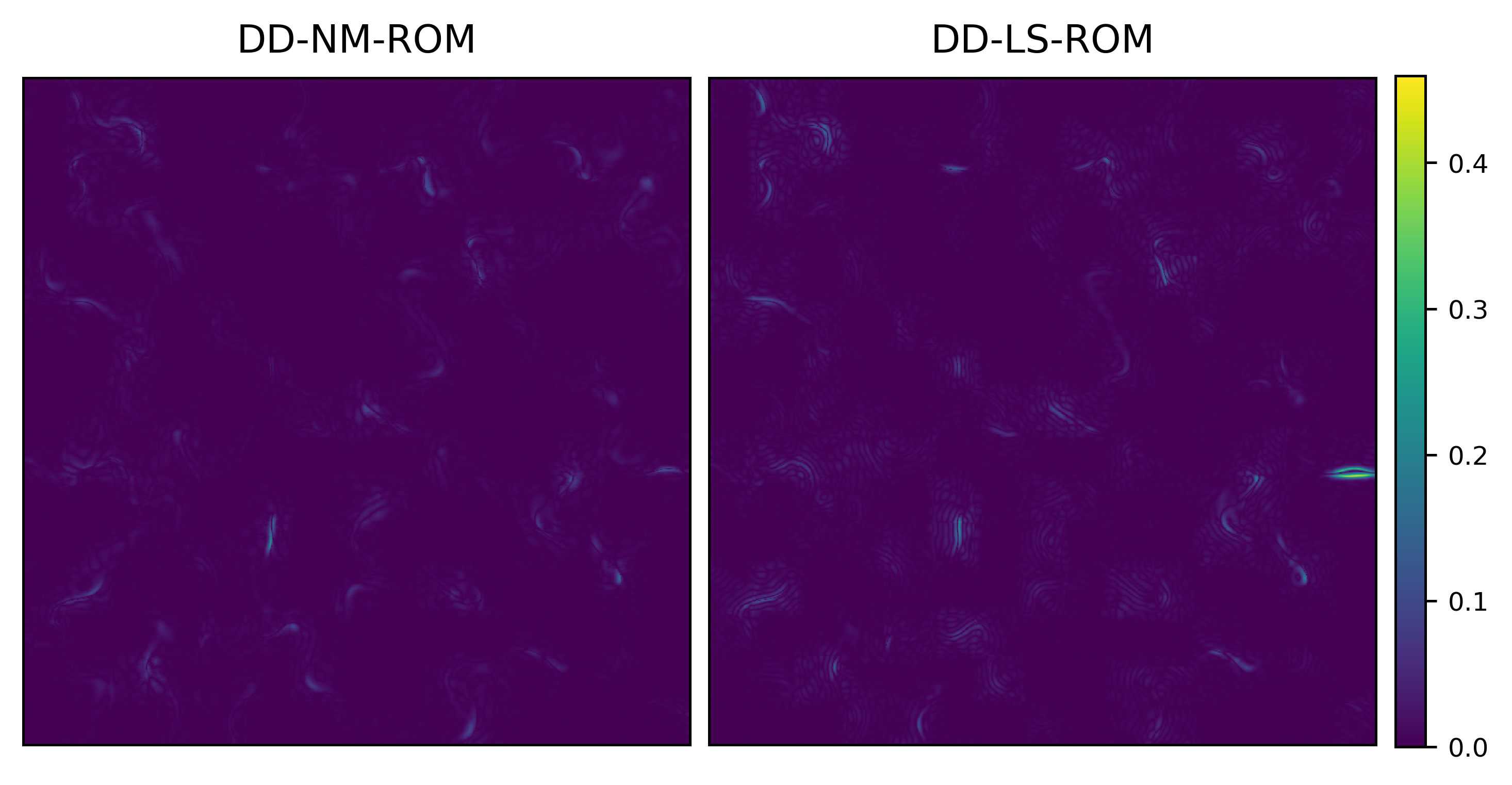}
\caption{Speed field and pointwise absolute speed error at \(t=1\) for one representative \(15\times15\) Burgers test case.
The top row shows the FOM, DD-NM-ROM with AE \((50,50)\), and DD-LS-ROM with POD \((65,65)\), from left to right.
The bottom row shows \(|s_{\mathrm{ROM}}-s_{\mathrm{FOM}}|\) for the two ROMs.
Panels in the same row use a common color scale.
The figure represents one case, whereas Table~\ref{tab:burgers_scaleup_results} reports averages over 10 cases.}
\label{fig:burgers_speed_solution}
\label{fig:burgers_speed_error}
\end{figure}

Figure~\ref{fig:burgers_speed_solution} shows the speed at the final time \(t=1\) for one test case,
\begin{equation}
 s=\sqrt{u^2+v^2}
 \label{eq:speed_definition}
\end{equation}
together with the pointwise absolute error \(|s_{\mathrm{ROM}}-s_{\mathrm{FOM}}|\).
Speed is used as a visualization quantity to examine the two velocity components together. The final-time error in this figure is distinct from the full-trajectory errors reported in the table.
The corresponding \(u\) and \(v\) solutions and errors for the same case are provided in Appendix~\ref{app:burgers_scaleup_fields}.

\subsubsection{Accuracy--efficiency trade-off}
\begin{figure}[pos=htbp]
\centering
\begin{minipage}[t]{0.49\textwidth}
\centering
\includegraphics[width=\linewidth]{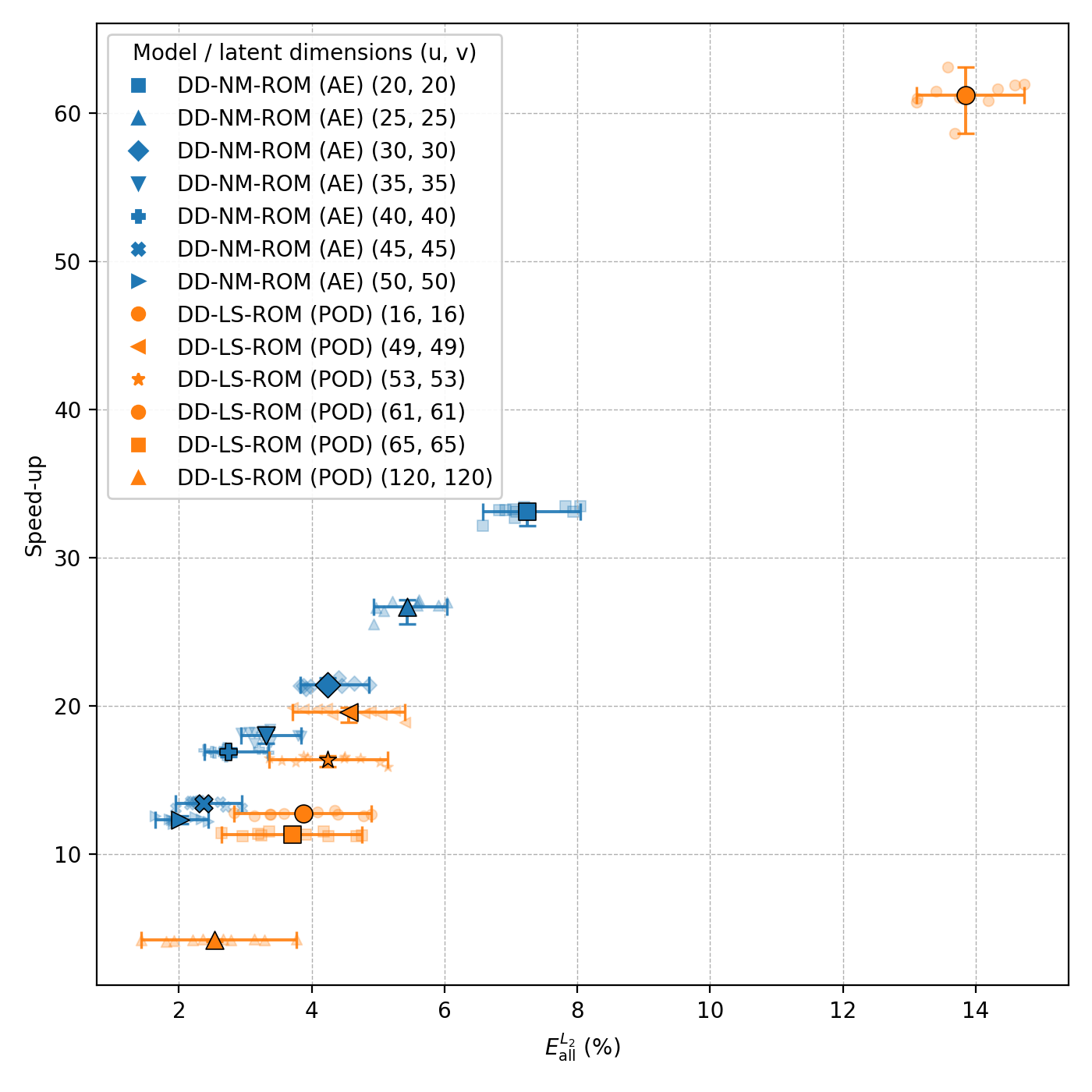}\\
{\small (a) Space--time relative error}
\end{minipage}\hfill
\begin{minipage}[t]{0.49\textwidth}
\centering
\includegraphics[width=\linewidth]{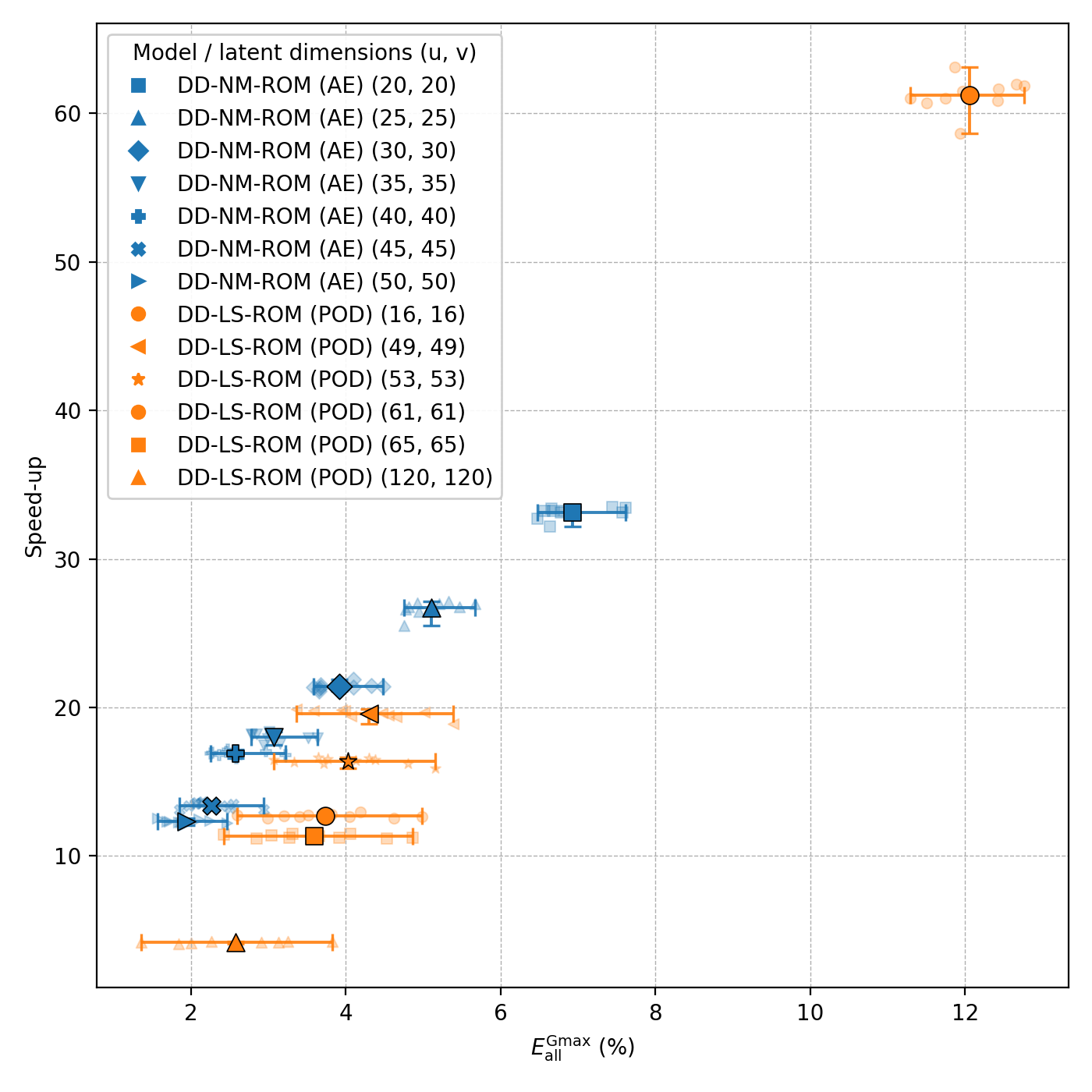}\\
{\small (b) Global-normalized maximum error}
\end{minipage}
\caption{Accuracy--efficiency trade-offs for the \(15\times15\) Burgers experiments.
The horizontal axes show \(E_{\mathrm{all}}^{L_2}\) in panel (a) and \(E_{\mathrm{all}}^{\mathrm{Gmax}}\) in panel (b). Both vertical axes show the speed-up over the FOM.
Small translucent markers represent individual test cases, large opaque markers represent arithmetic means over 10 cases, and horizontal and vertical bars span the observed minimum-to-maximum ranges.
Legend dimensions are listed in the order \((u,v)\).
All local representations are reused from the \(3\times3\) training configuration without retraining.}
\label{fig:burgers_speedup_error_tradeoff}
\label{fig:burgers_gmax_tradeoff}
\end{figure}

Figure~\ref{fig:burgers_speedup_error_tradeoff} presents the comparison in Table~\ref{tab:burgers_scaleup_results} as accuracy--efficiency trade-offs for the two error metrics.
Large markers are placed at the mean error and mean speed-up over 10 cases, while small markers and range bars show the case-wise distribution and minimum-to-maximum ranges.
Together, the two panels allow model-selection trends to be compared in terms of both accumulated trajectory error and the largest instantaneous deviation.

The smallest POD model, \((16,16)\), is the fastest model in the comparison, with a mean speed-up of \(61.24\), but has mean errors of \(E_{\mathrm{all}}^{L_2}=13.849\%\) and \(E_{\mathrm{all}}^{\mathrm{Gmax}}=12.058\%\).
In contrast, the smallest AE model, \((20,20)\), has a lower speed-up of \(33.12\) but smaller errors of \(7.249\%\) and \(6.927\%\), respectively.
Thus, even in the most strongly compressed settings, the nonlinear representation offers a lower-error option, whereas the linear representation offers a higher-speed but less accurate option.

The efficiency advantage of the nonlinear representation is more pronounced when higher accuracy is required.
For the AE \((50,50)\) and POD \((65,65)\) models compared above, the AE has smaller errors in both metrics at a similar speed-up.
Meanwhile, AE \((40,40)\) and POD \((120,120)\) have similar mean global-normalized maximum errors of \(2.576\%\) and \(2.575\%\), respectively, but speed-ups of \(16.89\) and \(4.21\).
This comparison illustrates the need to assess cost at a required accuracy or accuracy at an allowable cost, rather than merely matching the dimensions of different representations.

Taken together, the local reconstruction and scale-up prediction results show that a nonlinear representation learned on a small domain can be repeatedly applied to a larger Burgers problem and, at an appropriate reduced dimension, provides a more favorable accuracy--efficiency trade-off than POD.
These results are obtained by changing the domain and external boundary conditions while holding the small-domain local representation fixed.

\subsection{Incompressible Navier--Stokes equations}
\label{sec:navier_stokes_experiment}
\subsubsection{Problem setup}
The second example is the incompressible Navier--Stokes problem for velocity \(\mathbf u=(u,v)\) and pressure \(p\):
\begin{align}
 \frac{\partial\mathbf u}{\partial t}
 +(\mathbf u\cdot\nabla)\mathbf u
 &=-\frac1\rho\nabla p+\nu\nabla^2\mathbf u+\mathbf f,
 \label{eq:ns_momentum}\\
 \nabla\cdot\mathbf u&=0.
 \label{eq:ns_incompressibility}
\end{align}
The density is set to \(\rho=1\) and the kinematic viscosity to \(\nu=10^{-4}\), with periodic boundary conditions imposed on velocity and pressure in both directions.
In this example, we parameterize the spatial direction of the forcing rather than the initial conditions.
A localized Gaussian forcing source is placed at the center of each \(1\times1\) subdomain, and its direction is varied within \(\pm10^\circ\) of the upward direction.
Specifically, the direction of source \(s\) is defined as
\begin{equation}
 \theta_s=\frac\pi2+\mu_s,
 \qquad \mu_s\in[-\pi/18,\pi/18]
 \label{eq:forcing_angle_param}
\end{equation}
with combinations of forcing directions sampled using a Sobol sequence.
The source centers and widths are fixed, and only the directions are varied to generate diverse flow responses.
The spatial profile and time-dependent amplitude are described in Appendix~\ref{app:ns_forcing}.

\subsubsection{Training and scale-up configurations}
\begin{figure}[pos=htbp]
\centering
\includegraphics[width=0.32\textwidth]{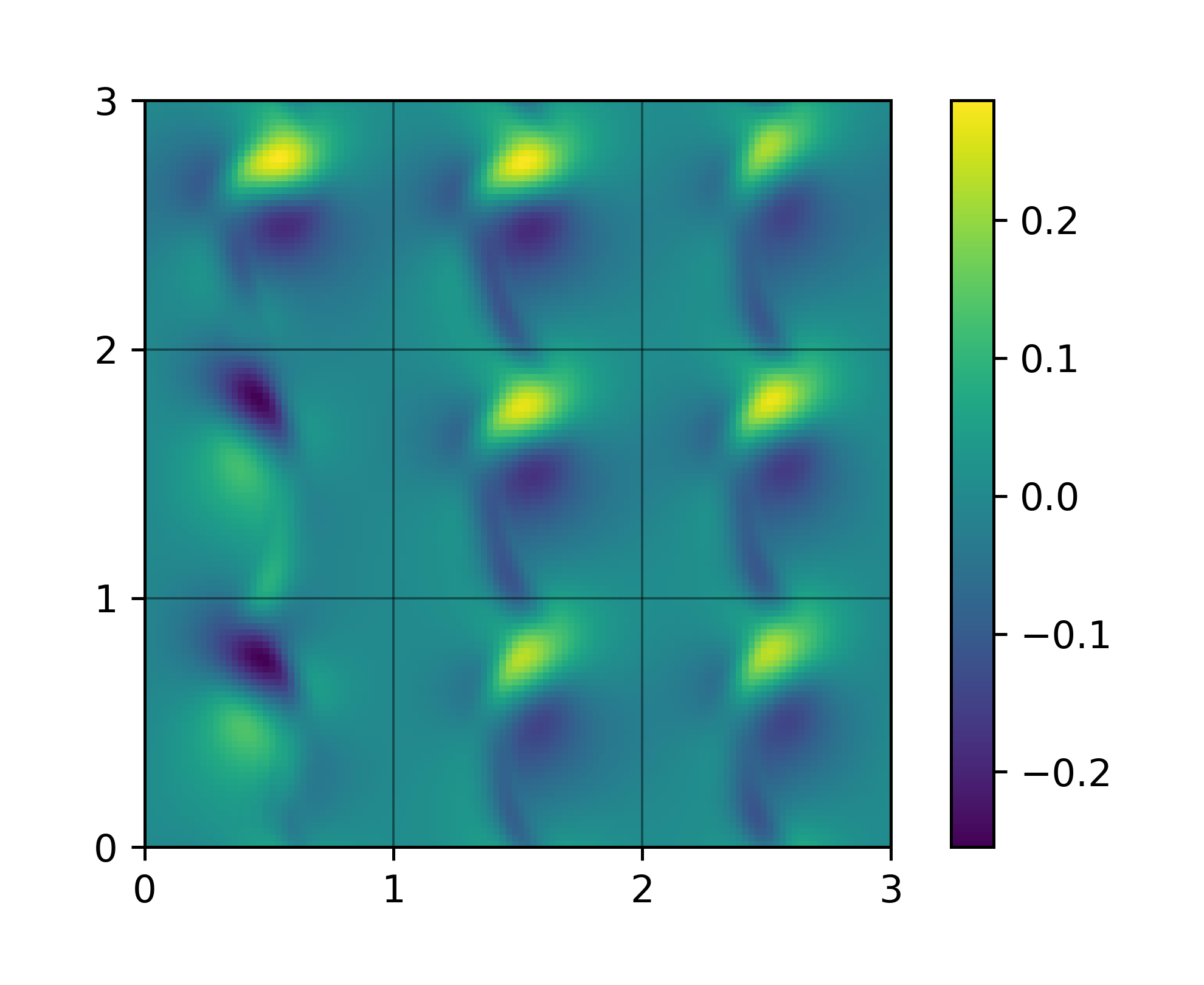}\hfill
\includegraphics[width=0.32\textwidth]{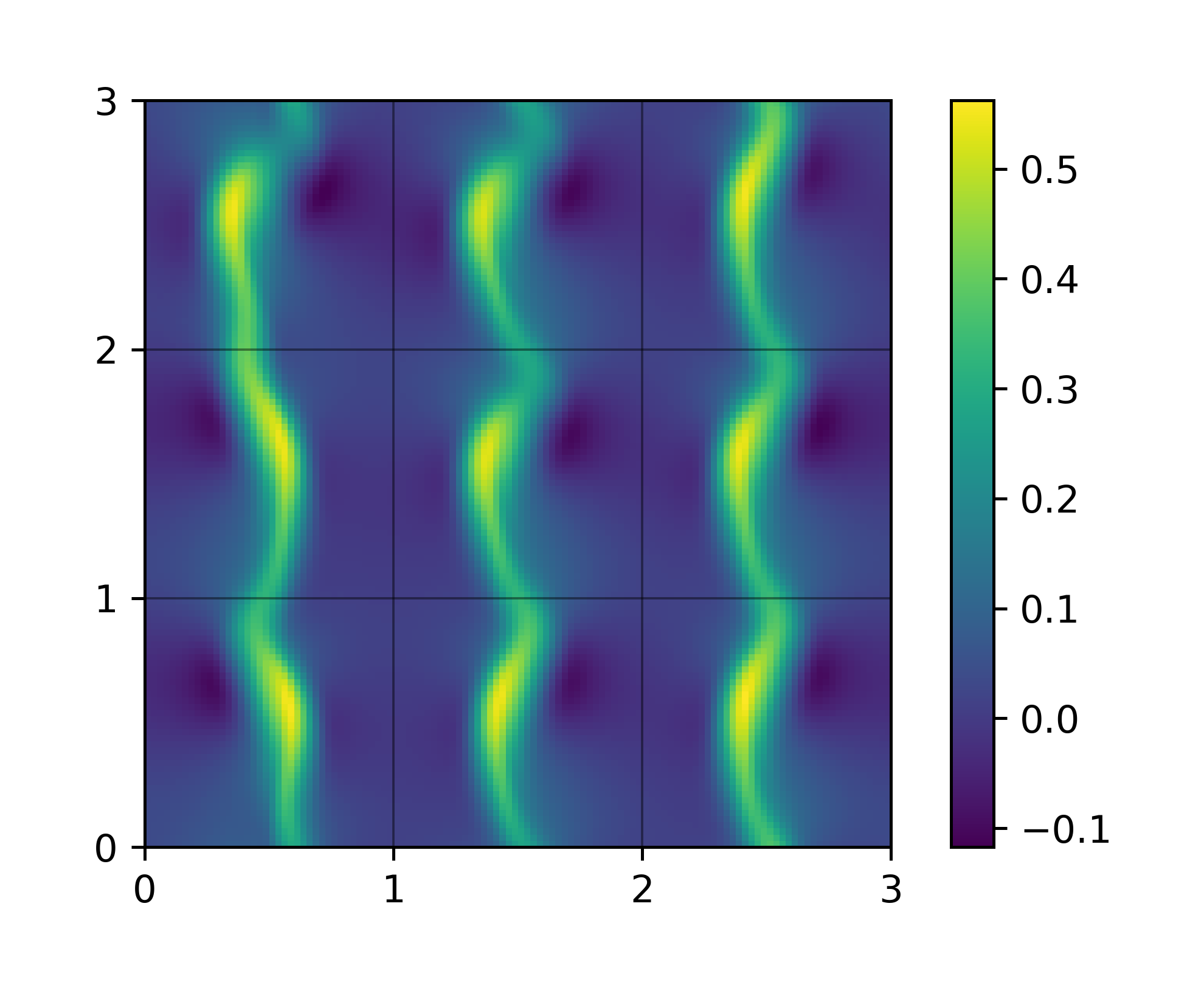}\hfill
\includegraphics[width=0.32\textwidth]{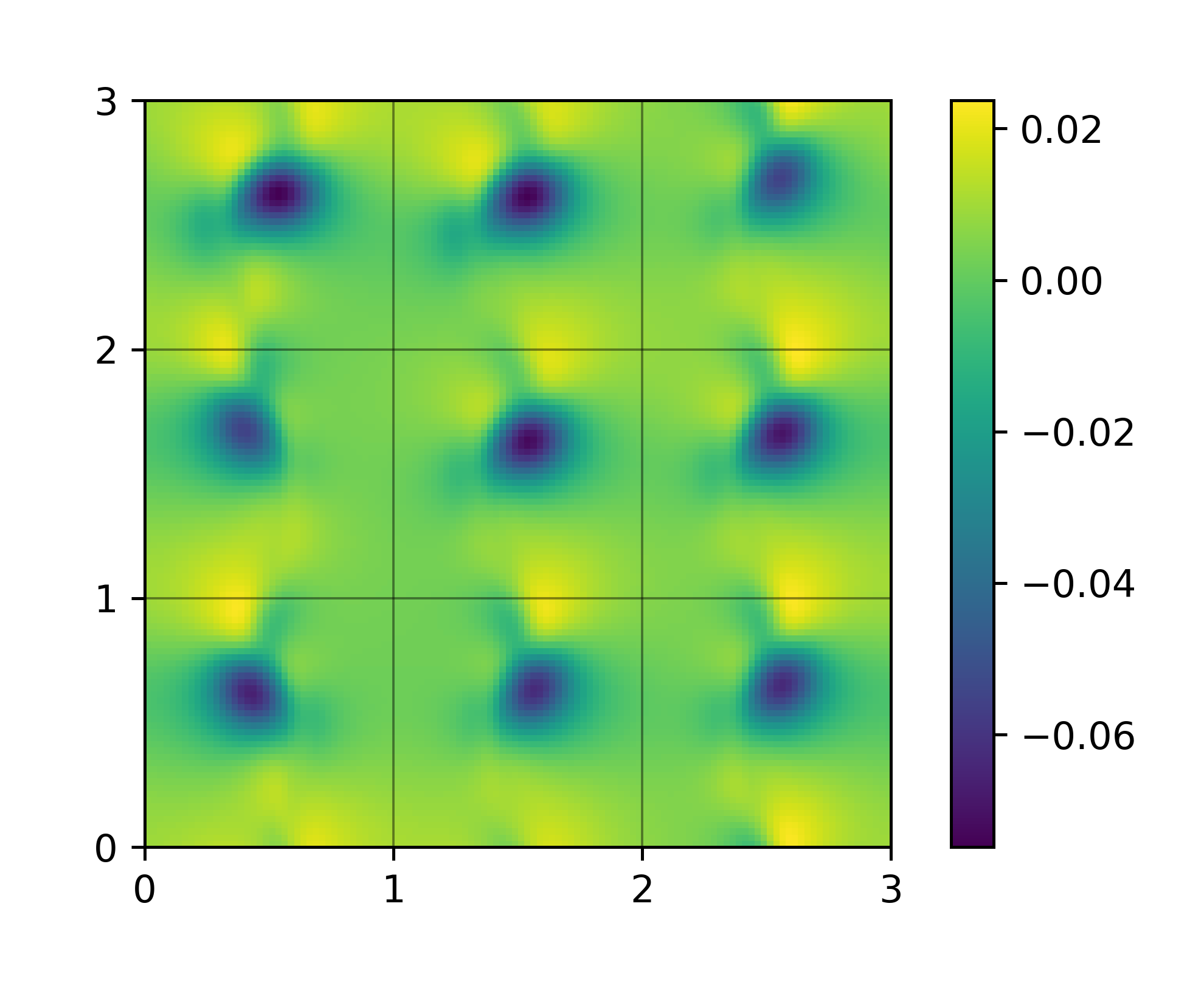}
\caption{Representative FOM training snapshot on the \(3\times3\) Navier--Stokes configuration.
The panels show \(u\), \(v\), and \(p\), from left to right, at the final time for one training realization.
Subdomain boundaries indicate the nine non-overlapping unit subdomains.
The local snapshots are pooled over subdomains, saved time instances, and training realizations to construct the field-wise AE and POD representations.}
\label{fig:training_data_example}
\end{figure}

Training data are generated on a small \(3\times3\) subdomain configuration.
Each subdomain is assigned \(40\times40\) degrees of freedom per field, giving a \(120\times120\) global grid and a field-wise local dimension of \(m=1600\).
The \(u,v,p\) fields are extracted subdomain-wise from each FOM time snapshot to construct the field-wise datasets.
Figure~\ref{fig:training_data_example} shows how nine local samples are obtained from a representative global training snapshot.

Local snapshots from 1000 small-domain forcing cases are randomly split into training and validation sets with a training fraction of 0.8.
Unlike the case-wise split used for Burgers, this split is performed at the local-snapshot level.
A separate set of 200 small-domain forcing cases, excluded from training and model selection, is generated to evaluate local reconstruction accuracy.
Details of the data split are provided in Appendix~\ref{sec:ns_data_split}.

In the scale-up experiment, the learned local representations are reused on a \(15\times15\) subdomain configuration.
The global grid is \(600\times600\), and the field-wise state size on each subdomain is the same as in the training configuration.
The large-domain reference FOM is computed for 10 separate forcing cases.
Neither the AE weights nor the POD bases are updated using large-domain data, so this comparison tests the application of a small-domain local representation to a coupled problem on a larger domain.

\subsubsection{Reduced models and solver settings}
\label{sec:solver_settings}
Table~\ref{tab:reduced_model_configs} summarizes the reduced dimensions used in the Navier--Stokes experiments.
The AE configurations are \((5,5,5)\), \((7,7,7)\), and \((10,10,10)\), with equal dimensions for all three fields.
The POD bases are selected field-wise according to cumulative energy-retention levels of \(99\%\), \(99.9\%\), and \(99.99\%\). The corresponding dimensions are listed in the table.
The AEs for the three fields are trained independently, but each learned field-wise decoder is shared across all subdomains.

\begin{table}[pos=htbp]
\centering
\small
\caption{Reduced model configurations for the Navier--Stokes experiments.
Field-wise dimensions are AE latent dimensions or POD ranks, listed in the order \((u,v,p)\).
The total local dimension is the sum of the three component dimensions on one subdomain.}
\label{tab:reduced_model_configs}
\begin{tabular}{lcc}
\toprule
Model & Field-wise dimension & Total local dimension \\
\midrule
AE (5,5,5) & \([n_u,n_v,n_p]=[5,5,5]\) & 15 \\
AE (7,7,7) & \([n_u,n_v,n_p]=[7,7,7]\)  & 21 \\
AE (10,10,10) & \([n_u,n_v,n_p]=[10,10,10]\)  & 30 \\
\midrule
POD (20,17,17) & \([n_u,n_v,n_p]=[20,17,17]\) & 54 \\
POD (46,42,37) & \([n_u,n_v,n_p]=[46,42,37]\) & 125 \\
POD (87,79,67) & \([n_u,n_v,n_p]=[87,79,67]\) & 233 \\
\bottomrule
\end{tabular}
\end{table}

The Navier--Stokes encoder has one hidden layer and uses standardized inputs.
The decoder has one hidden layer and a sparse output layer that follows spatial connectivity.
In this experiment, the encoder width and decoder-mask parameters vary together with the latent dimension. The specific values are given in Appendix~\ref{app:architecture_settings}.

Time integration uses fully implicit backward Euler.
With \(\Delta t=0.1\), 50 time steps are computed to reach the final time \(T=5\).
The ROM starts from \(u=v=p=0\), and the initial reduced coordinates are obtained using the field-wise encoders.
Prediction errors are evaluated using states saved at \(t=0.1,0.2,\ldots,5\).
The FOM is solved by Newton's method with velocity and pressure treated jointly as unknowns.
Because the Newton linear systems are generally not symmetric positive definite due to velocity--pressure coupling and convection, we use GMRES.
The time-discrete residual, including pressure stabilization, and the grid configuration are described in Appendix~\ref{app:ns_discretization}.

\paragraph{Treatment of the global mean pressure.}
\label{sec:ns_pressure_gauge}
To set a pressure reference under periodic boundary conditions, we add a global mean-pressure penalty to the ROM's PDE residual objective.
Let \(N_p\) denote the number of global pressure degrees of freedom and \(\bar p(\hat x)\) the mean of the reconstructed global pressure. The reduced objective actually minimized is
\begin{equation}
 \mathcal F(\hat x)=\mathcal F_{\mathrm{PDE}}(\hat x)
 +\frac{\gamma N_p}{2}\bigl(\bar p(\hat x)-p_\star\bigr)^2,
 \qquad \gamma=1,\quad p_\star=0
 \label{eq:pressure_penalty_main}
\end{equation}
This term weakly controls the mean pressure over the entire domain, rather than subdomain-wise means, and is distinct from a constraint or penalty for interface compatibility.
The definition of the global mean and the additional Gauss--Newton contributions are given in Appendix~\ref{app:pressure_penalty}.

This objective is minimized by Gauss--Newton iterations, with each linear system solved using subdomain-wise block-Jacobi-preconditioned PCG.
The stopping criteria, linear solver tolerances, and preconditioning details are provided in Appendix~\ref{app:ns_rom_solver}.

\subsubsection{Representation accuracy on held-out local snapshots}
\label{sec:ns_representation}
Table~\ref{tab:representation_accuracy} compares local representation errors in \(u,v,p\) obtained from separate small-domain test cases.
As in the Burgers experiments, we use encoder--decoder reconstruction for AE models and orthogonal projection for POD models, measured by the relative Frobenius error in Eq.~\eqref{eq:representation_error}.
This comparison evaluates how compactly the nonlinear representation describes each field relative to a linear basis.

\begin{table}[pos=htbp]
\centering
\small
\caption{Field-wise representation errors on held-out local Navier--Stokes snapshots, measured by the relative Frobenius error in Eq.~\eqref{eq:representation_error}.
AE values are encoder--decoder reconstruction errors, and POD values are orthogonal projection errors.
The test snapshots are extracted from separate \(3\times3\) trajectories not used for training or model selection.
The dimension \(r\) is specified per field; the reported errors are dimensionless, not percentages.}
\label{tab:representation_accuracy}
\begin{tabular}{llcc}
\toprule
Model & Field & \(r\) & $\varepsilon_{\mathrm{rep},q}$ \\
\midrule
\(\mathrm{AE}\) & \(u\) & 5  & \(4.227\times10^{-2}\) \\
\(\mathrm{AE}\) & \(u\) & 7  & \(2.385\times10^{-2}\) \\
\(\mathrm{AE}\) & \(u\) & 10 & \(1.103\times10^{-2}\) \\
\(\mathrm{AE}\) & \(v\) & 5  & \(1.891\times10^{-2}\) \\
\(\mathrm{AE}\) & \(v\) & 7  & \(1.062\times10^{-2}\) \\
\(\mathrm{AE}\) & \(v\) & 10 & \(4.285\times10^{-3}\) \\
\(\mathrm{AE}\) & \(p\) & 5  & \(3.716\times10^{-2}\) \\
\(\mathrm{AE}\) & \(p\) & 7  & \(2.257\times10^{-2}\) \\
\(\mathrm{AE}\) & \(p\) & 10 & \(1.109\times10^{-2}\) \\
\midrule
\(\mathrm{POD}\) & \(u\) & 20 & \(9.4765\times10^{-2}\) \\
\(\mathrm{POD}\) & \(u\) & 46 & \(2.9032\times10^{-2}\) \\
\(\mathrm{POD}\) & \(u\) & 87 & \(9.3822\times10^{-3}\) \\
\(\mathrm{POD}\) & \(v\) & 17 & \(5.4155\times10^{-2}\) \\
\(\mathrm{POD}\) & \(v\) & 42 & \(1.6980\times10^{-2}\) \\
\(\mathrm{POD}\) & \(v\) & 79 & \(5.4596\times10^{-3}\) \\
\(\mathrm{POD}\) & \(p\) & 17 & \(9.7261\times10^{-2}\) \\
\(\mathrm{POD}\) & \(p\) & 37 & \(3.0629\times10^{-2}\) \\
\(\mathrm{POD}\) & \(p\) & 67 & \(9.7189\times10^{-3}\) \\
\bottomrule
\end{tabular}
\end{table}

At an AE field-wise dimension of 10, the errors in \(u\) and \(p\) are \(1.103\times10^{-2}\) and \(1.109\times10^{-2}\), respectively, or approximately \(1\%\).
These are close to the POD errors of \(9.3822\times10^{-3}\) and \(9.7189\times10^{-3}\), obtained using 87 basis vectors for \(u\) and 67 for \(p\).
For \(v\), the AE error at dimension 10 is \(4.285\times10^{-3}\), smaller than the POD error of \(5.4596\times10^{-3}\) at dimension 79.
Thus, the nonlinear representation effectively describes the test snapshots with a small reduced dimension for all three fields.
We next examine whether this compactness also translates into accuracy and computational efficiency in residual-based time integration on the larger domain.

\subsubsection{Scale-up prediction on a \(15\times15\) configuration}
\label{sec:ns_scaleup_results}
Table~\ref{tab:scaleup_results} summarizes errors in the full state \([u;v;p]\) and mean speed-ups for the 10 large-domain test cases.
Each model solves the time-discrete residual on the larger domain while holding fixed the local representation learned or constructed on the small \(3\times3\) configuration.

\begin{table}[pos=htbp]
\centering
\small
\caption{Scale-up performance on the \(15\times15\) Navier--Stokes configuration.
Errors and speed-ups are arithmetic means over 10 test cases.
All representations are constructed from the \(3\times3\) data and reused without retraining.
The total local dimension is \(n_{\mathrm{loc}}\), both error columns are percentages, and \(S\) is the speed-up over the FOM.}
\label{tab:scaleup_results}
\begin{tabular}{lcccc}
\toprule
Model & \(n_{\mathrm{loc}}\) &
\(E_{\mathrm{all}}^{L_2}\)(\%) &
\(E_{\mathrm{all}}^{\mathrm{Gmax}}\)(\%) &
\(S\) \\
\midrule
AE (5,5,5)      & 15 & \(7.881\) & \(13.90\) & 33.99 \\
AE (7,7,7)      & 21 & \(4.702\) & \(7.934\) & 23.87 \\
AE (10,10,10)   & 30 & \(3.484\) & \(5.520\) & 21.13 \\
\midrule
POD (20,17,17)     & 54 & \(20.94\) & \(32.72\) & 79.45 \\
POD (46,42,37)     & 125 & \(7.298\) & \(12.07\) & 19.20 \\
POD (87,79,67)     & 233 & \(3.229\) & \(5.299\) & 3.99 \\
\bottomrule
\end{tabular}
\end{table}

With a total local dimension of 30, AE \((10,10,10)\) achieves a mean space--time error of \(3.484\%\), a mean global-normalized maximum error of \(5.520\%\), and a mean speed-up of \(21.13\).
POD \((46,42,37)\), which has a comparable mean speed-up of \(19.20\), uses a local dimension of 125 and has errors of \(7.298\%\) and \(12.070\%\), respectively.
Thus, at a comparable computational speed, the nonlinear representation provides lower mean errors with a much smaller local dimension.
These two models are used for the qualitative comparison of representative states.

\begin{figure}[pos=htbp]
\centering
\includegraphics[width=0.95\textwidth]{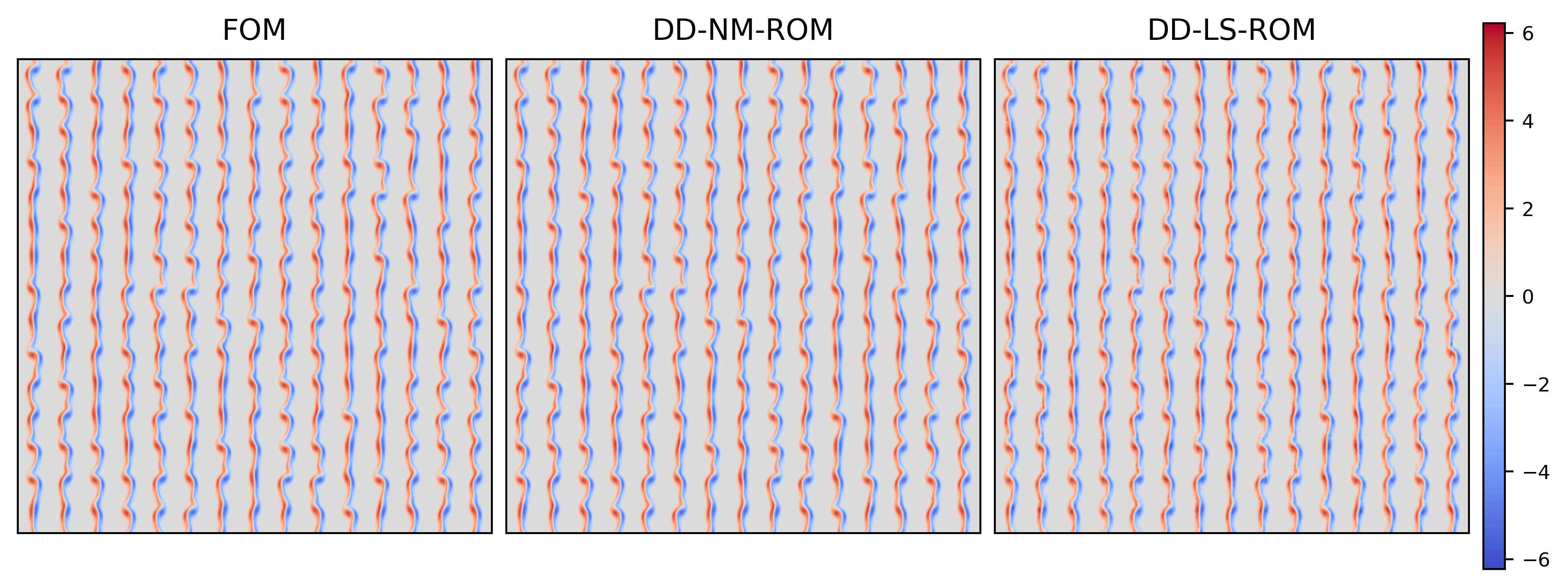}\\
\includegraphics[width=0.69\textwidth]{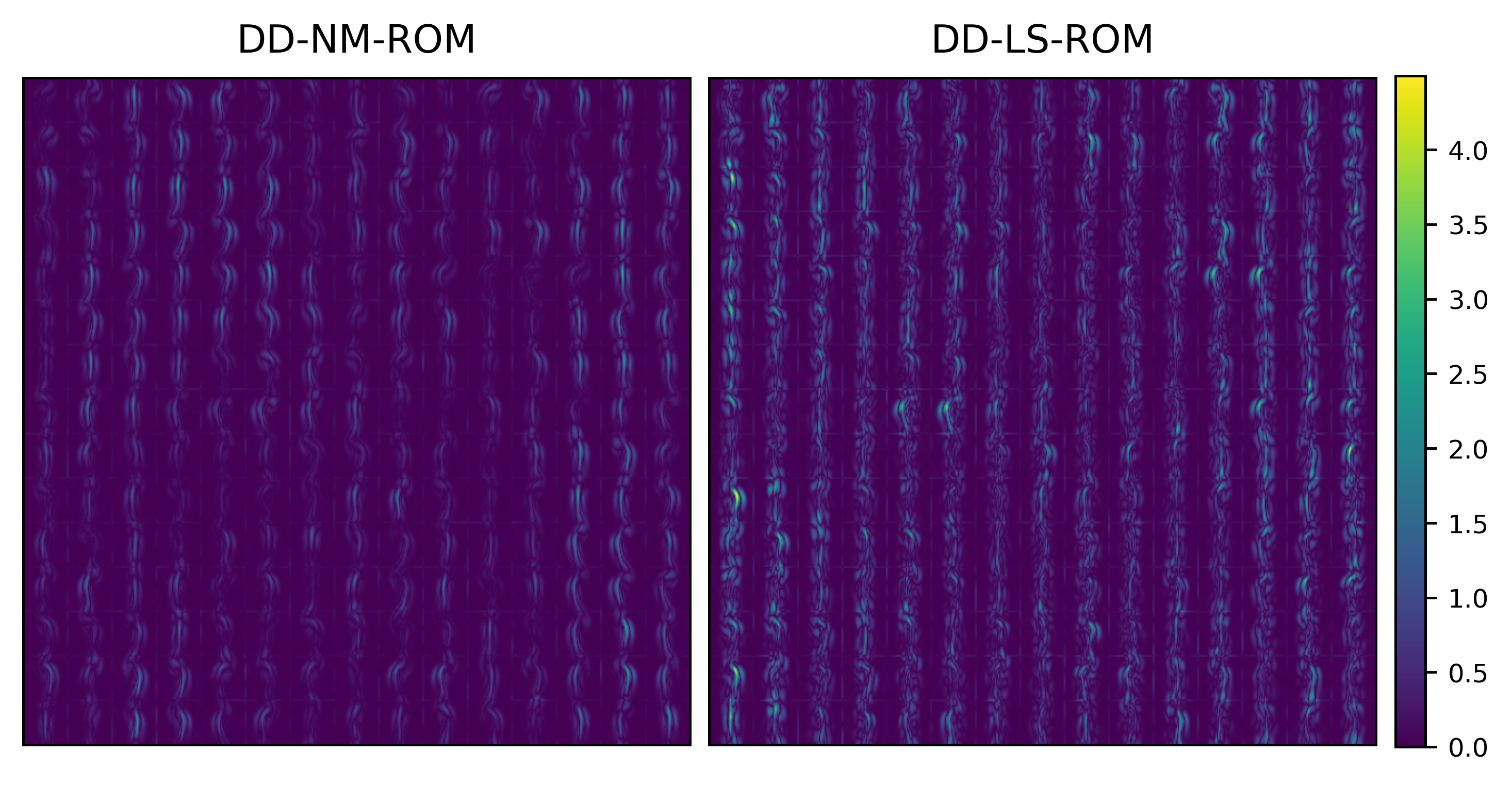}
\caption{Vorticity field and pointwise absolute vorticity error at \(t=5\) for one representative \(15\times15\) Navier--Stokes test case.
The top row shows the FOM, DD-NM-ROM with AE \((10,10,10)\), and DD-LS-ROM with POD \((46,42,37)\), from left to right.
The bottom row shows \(|\omega_{\mathrm{ROM}}-\omega_{\mathrm{FOM}}|\) for the two ROMs.
The same derivative operator is applied to the FOM and ROM velocities, and panels in each row use a common color scale.
The figure represents one case, whereas Table~\ref{tab:scaleup_results} reports averages over 10 cases.}
\label{fig:scaleup_vorticity_solution}
\label{fig:scaleup_vorticity_error}
\end{figure}

Figure~\ref{fig:scaleup_vorticity_solution} compares vorticity and pointwise absolute vorticity errors at \(t=5\) for one test case.
Vorticity is defined as
\begin{equation}
 \omega=\frac{\partial v}{\partial x}-\frac{\partial u}{\partial y}
 \label{eq:vorticity_definition}
\end{equation}
and is computed by applying the same discrete derivative operator to the FOM and ROM velocity fields.
It is a postprocessed quantity obtained from the velocity components of the learned state \(u,v,p\), used to examine flow structures and spatial errors.
The pointwise visualization error is
\begin{equation}
 e_\omega^{\mathrm{ROM}}(\mathbf x,t)
 =|\omega_{\mathrm{ROM}}(\mathbf x,t)-\omega_{\mathrm{FOM}}(\mathbf x,t)|
 \label{eq:vorticity_visual_error}
\end{equation}
which is distinct from the full-trajectory errors in the primitive variables reported in the table.
The \(u,v,p\) states and errors for the same representative case are provided in Appendix~\ref{app:ns_scaleup_fields}.

\subsubsection{Accuracy--efficiency trade-off}
\begin{figure}[pos=htbp]
\centering
\begin{minipage}[t]{0.49\textwidth}
\centering
\includegraphics[width=\linewidth]{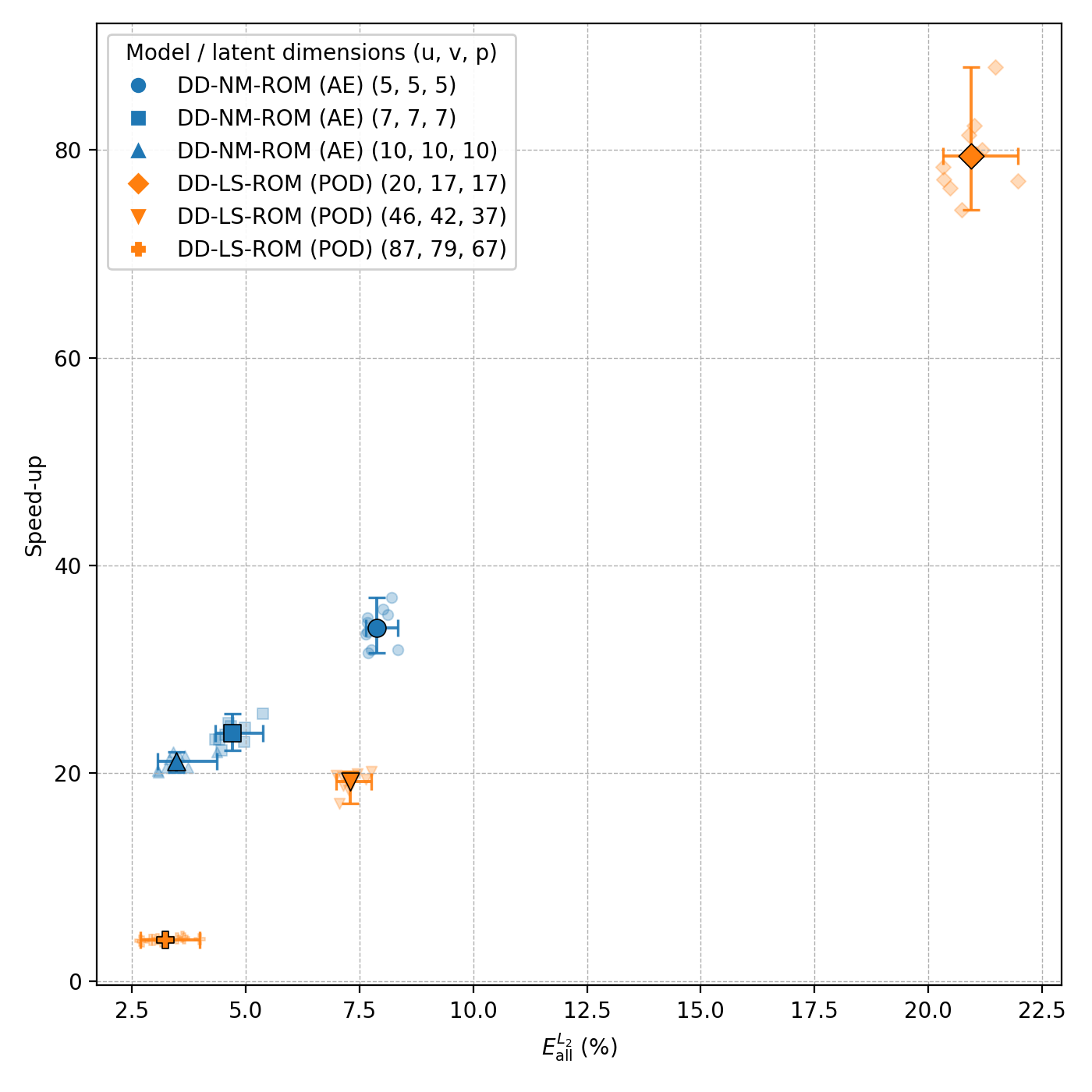}\\
{\small (a) Space--time relative error}
\end{minipage}\hfill
\begin{minipage}[t]{0.49\textwidth}
\centering
\includegraphics[width=\linewidth]{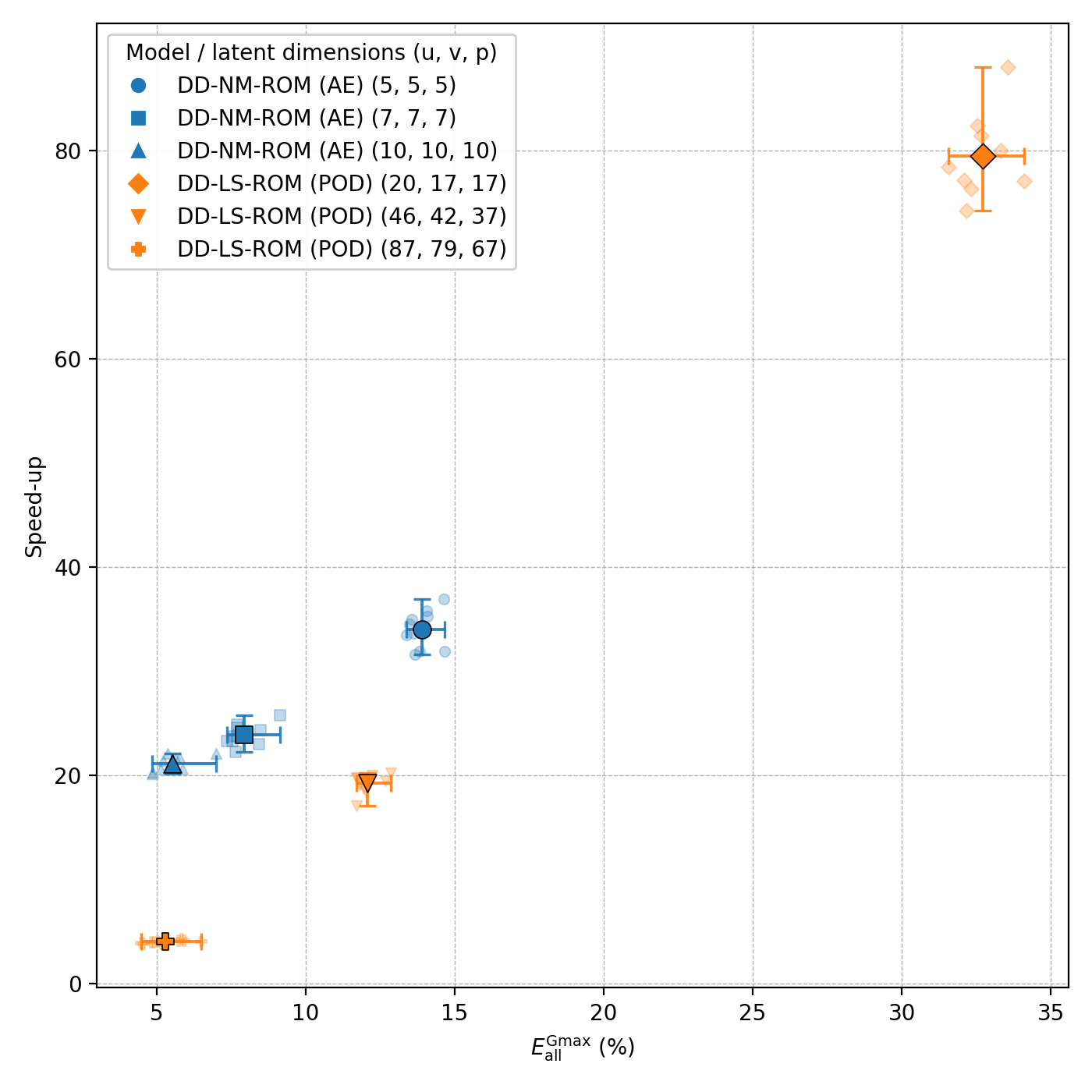}\\
{\small (b) Global-normalized maximum error}
\end{minipage}
\caption{Accuracy--efficiency trade-offs for the \(15\times15\) Navier--Stokes experiments.
The horizontal axes show \(E_{\mathrm{all}}^{L_2}\) in panel (a) and \(E_{\mathrm{all}}^{\mathrm{Gmax}}\) in panel (b). Both vertical axes show the speed-up over the FOM.
Small translucent markers represent individual test cases, large opaque markers represent arithmetic means over 10 cases, and horizontal and vertical bars span the observed minimum-to-maximum ranges.
Legend dimensions are listed in the order \((u,v,p)\).
All local representations are reused from the \(3\times3\) training configuration without retraining.}
\label{fig:pareto_scaleup}
\label{fig:ns_gmax_tradeoff}
\end{figure}

Figure~\ref{fig:pareto_scaleup} shows the accuracy--efficiency trade-offs for the two trajectory error metrics.
As in the Burgers results, large markers denote arithmetic means over 10 cases, small markers denote individual cases, and horizontal and vertical bars span the observed minimum-to-maximum ranges.
The plots therefore show both the mean performance of each model and the variation in error and speed-up across test cases.

Among the smallest-dimension configurations, AE \((5,5,5)\) has smaller errors than POD \((20,17,17)\) in both metrics.
Their mean space--time errors are \(7.881\%\) and \(20.94\%\), respectively, and their mean global-normalized maximum errors are \(13.90\%\) and \(32.72\%\).
The POD model, however, has a higher speed-up of \(79.45\), compared with \(33.99\) for the AE.
This shows that under strong compression, the fastest model and the more accurate model need not be the same.

At intermediate computational speeds, AE \((10,10,10)\) provides lower errors than POD \((46,42,37)\).
Comparing AE \((10,10,10)\) with the largest POD model, \((87,79,67)\), gives mean space--time errors of \(3.484\%\) and \(3.229\%\), respectively, which are in a similar range.
The POD model attains the smallest values in both error metrics, but uses a total local dimension of 233 and achieves a mean speed-up of only \(3.99\).
In contrast, the AE maintains a mean speed-up of \(21.13\) with a local dimension of 30.
The advantage of the nonlinear representation is therefore not that it achieves the smallest error in every setting, but that it offers more favorable computational costs in a practical accuracy range.

Taken together, the representation-accuracy results, scale-up predictions, and two accuracy--efficiency comparisons show that a local nonlinear representation learned on a small \(3\times3\) domain can be applied to a \(15\times15\) Navier--Stokes problem without retraining.
This demonstrates that the same non-overlapping residual-coupling construction can be used for a problem with different state variables and parameterization from Burgers.
In both experiments, the nonlinear representation provides a more favorable accuracy--efficiency trade-off than POD at an appropriate local dimension.

\clearpage
\section{Discussion and Conclusions}
\label{sec:conclusion}
This work presented a bottom-up reduced-order modeling framework that reuses local nonlinear representations learned on a small training domain for larger PDE problems.
The full state owned by each subdomain is represented using a shared decoder and independent reduced coordinates, and this representation is substituted into the original time-discrete residual to construct the global reduced problem.
Because the number of subdomains can be increased while holding the trained decoder fixed, there is no need to relearn a nonlinear representation on the full large domain.

In the proposed state construction, each global degree of freedom is represented by only one subdomain, and neighboring residuals directly use its reconstructed value.
Algebraic compatibility for the same degree of freedom therefore holds without separate reduced models for the interior and interface or compatibility constraints between independent interface copies.
This distinguishes the approach from constraint-based methods with separate interior/port representations. It also retains the coupling in the original FOM residual rather than introducing a separate DG interface formulation to couple local representations.
The central feature of this work is a single formulation that combines representation reuse through a shared decoder with compatibility through a non-overlapping state construction.

The reduced coordinates are determined jointly from a reduced global residual minimization problem at each time step.
Using the dependence of each local residual only on its own subdomain state and required neighboring states, we compute the reduced Jacobian as a product of the local full-order Jacobian and decoder Jacobians, and assemble the Gauss--Newton system from local contributions.
The global mean pressure in Navier--Stokes is handled by adding a separate mean-pressure term to this common PDE residual formulation.

In the numerical experiments, representations constructed from \(3\times3\) data were applied to \(15\times15\) domains without retraining for both the Burgers problem with parameterized initial conditions and the Navier--Stokes problem with parameterized forcing.
For Burgers, AE \((50,50)\) achieved a mean space--time error of \(2.023\%\) and a mean speed-up of \(12.33\) over 10 cases, with a lower error than POD \((65,65)\) at a similar speed.
For Navier--Stokes, AE \((10,10,10)\) achieved a mean space--time error of \(3.484\%\) and a mean speed-up of \(21.13\).
It was more accurate than POD \((46,42,37)\) at a comparable speed and used a smaller local dimension with greater computational efficiency than POD \((87,79,67)\) in a similar error range.
These results demonstrate the utility of the nonlinear representation within the same residual-coupling structure. They do not establish a direct performance advantage over existing constraint-based or DG-based ROMs.

The validation in this work is limited to domain scale-up with uniform local state configurations and within the specified ranges of initial conditions and forcing.
Algebraic compatibility of the discrete state and approximation accuracy on a larger domain are distinct properties. Separate evaluation is needed when domain enlargement produces global structures absent from the small training domain.
The experiments also do not quantify savings in offline data-generation and training costs, parallel scalability with the number of subdomains, or sensitivity to nonlinear stopping tolerances.

An important direction for future work is to combine the local residual structure with hyper-reduction to further reduce the cost of nonlinear residual and Jacobian evaluations.
DEIM \citep{chaturantabut2010nonlinear-a21}, solution-based nonlinear subspaces, and residual-sampling strategies \citep{choi2020sns:-40d} offer possible means of making this extension.
Hyper-reduction was not used in the present experiments.
The scope can also be extended to connecting representations across nonuniform subdomain sizes and resolutions, broader parameter ranges and more complex geometries, and parallel computation that exploits local assembly.
These extensions should be used to evaluate the accuracy and computational efficiency of bottom-up nonlinear ROMs that reuse small-domain representations across a wider range of problems.

\clearpage
\appendix
\counterwithin{equation}{section}
\counterwithin{figure}{section}
\counterwithin{table}{section}

\section{Full-Order Problems and Data Generation}
\label{sec:data_generation}
This appendix describes the generation of full-order trajectories on small training domains and their conversion into subdomain-wise training samples.
For both problems, each \(3\times3\) global snapshot is partitioned into nine non-overlapping local states for AE training and POD basis construction.
We first describe the domain, parameterization, and temporal and spatial discretizations for each problem, followed by the common snapshot collection procedure and the problem-specific data splits.

\subsection{Domain and initial conditions for the Burgers problem}
\label{app:burgers_data_generation}
The Burgers problem uses Eq.~\eqref{eq:burgers_momentum} from the main text with the viscosity fixed at \(\nu=0.005\).
The training domain \([0,1.5]^2\) is discretized on a \(150\times150\) cell-centered grid and partitioned into \(3\times3\) subdomains.
Each subdomain has size \(0.5\times0.5\) and contains \(50\times50\) cells, so \(\Delta x=\Delta y=0.01\).
Periodic boundary conditions are imposed in both directions on the small domain.

For large-domain evaluation, the domain is enlarged to \([0,7.5]^2\) with a \(750\times750\) grid, and homogeneous Neumann conditions are imposed on the outer boundary:
\begin{equation}
 \frac{\partial q}{\partial n}=0,\qquad q\in\{u,v\}
 \label{eq:burgers_neumann_bc}
\end{equation}
Here, \(n\) denotes the outward normal direction.
The local cell count and grid spacing are the same as on the small domain, and the total number of global degrees of freedom for the two velocity components is \(2\times750^2=1{,}125{,}000\).

\begin{figure}[pos=htbp]
\centering
\includegraphics[width=0.48\textwidth]{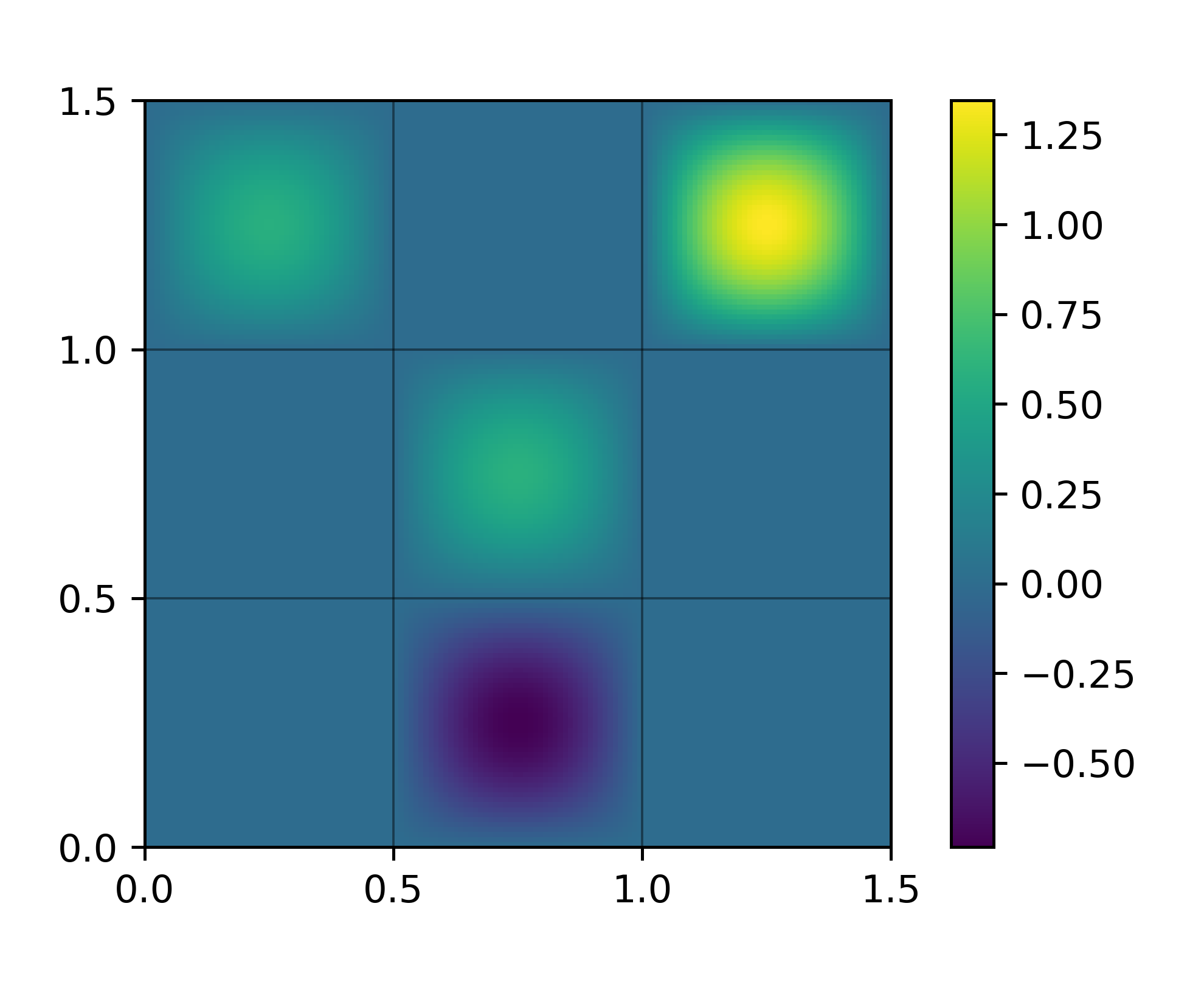}\hfill
\includegraphics[width=0.48\textwidth]{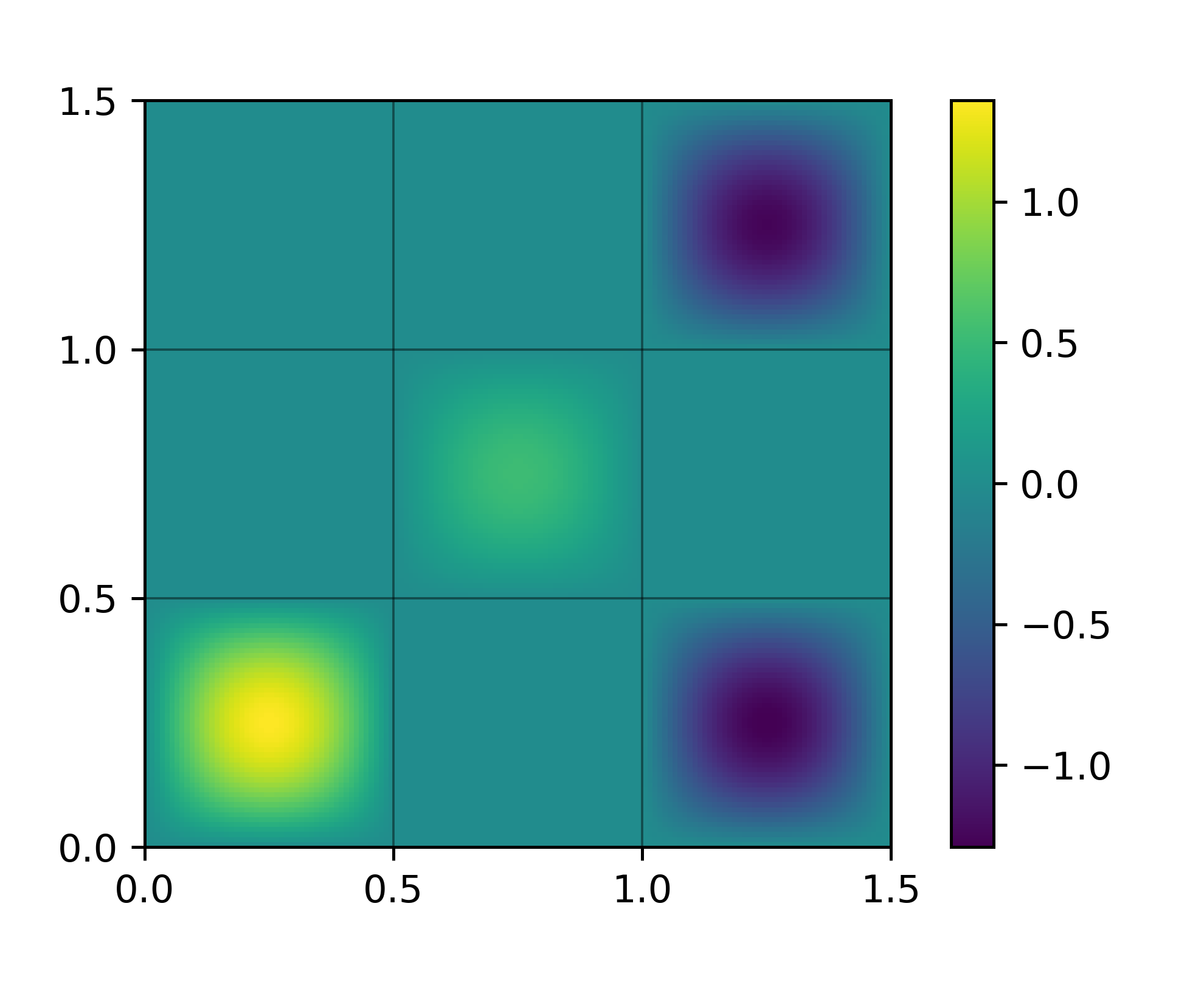}
\caption{Initial condition at \(t=0\) for the Burgers training case shown at \(t=1\) in Figure~\ref{fig:burgers_training_data_example}.
The left and right panels show \(u\) and \(v\), respectively.
The \(150\times150\) grid is divided into a \(3\times3\) configuration, and the activation, sign, and amplitude are sampled independently for each subdomain and velocity component.}
\label{fig:burgers_initial_condition}
\end{figure}

The initial conditions are constructed by multiplying the global sinusoidal profile in Eq.~\eqref{eq:burgers_ic} by subdomain-specific coefficients.
The sine functions are evaluated at the global cell-center coordinates and vanish at subdomain boundaries because the subdomain side length is \(0.5\).
The coefficient for field \(q\) and subdomain \(i\) is sampled as
\begin{align}
 \mu_i^q&=(-1)^{\eta_i^q}\gamma_i^q\xi_i^q,
 \label{eq:burgers_ic_coefficient}\\
 \gamma_i^q&\sim\mathcal U(0.5,1.5),\qquad
 \xi_i^q\sim\operatorname{Bernoulli}(0.5),\qquad
 \eta_i^q\sim\operatorname{Bernoulli}(0.5)
 \label{eq:burgers_ic_sampling}
\end{align}
The amplitude, activation, and sign are determined by \(\gamma_i^q\), \(\xi_i^q\), and \(\eta_i^q\), respectively, with independent sampling across all subdomains and fields.
This construction assigns different initial velocity magnitudes and directions to different subdomains.
The parameter dimension is 18 for the \(3\times3\) configuration and 450 for the \(15\times15\) configuration.
The base random seed is \(20260721\), with separate random streams for data splits and individual cases.
Figure~\ref{fig:burgers_initial_condition} shows the initial state of the same training case as Figure~\ref{fig:burgers_training_data_example} in the main text.

\subsection{Discretization and full-order solver for Burgers}
\subsubsection{Spatial discretization}
\label{app:burgers_spatial_discretization}
On a uniform cell-centered finite volume grid, the face-normal velocity is computed as the arithmetic mean of the velocities in the two adjacent cells.
A smooth upwind splitting is used for the convective flux.
For a face-crossing velocity \(a\), the positive and negative parts are defined as
\begin{equation}
 a^\pm=\frac12\left(a\pm\sqrt{a^2+\epsilon_{\mathrm{up}}^2}\right),
 \qquad \epsilon_{\mathrm{up}}=1.5\times10^{-6}
 \label{eq:burgers_smooth_upwind}
\end{equation}
For a face oriented in the positive coordinate direction, let \(q_L\) and \(q_R\) denote the field values in the two adjacent cells. The numerical flux is defined as
\begin{equation}
 F_q=a^+q_L+a^-q_R
\end{equation}

Using the face velocities \(a_x\) and \(a_y\) and face-flux divergence operators \(D_x\) and \(D_y\), the nonconservative convection term is discretized as
\begin{equation}
 C_h(q)=D_xF_{x,q}+D_yF_{y,q}
 -q\odot(D_xa_x+D_ya_y)
 \label{eq:burgers_discrete_convection}
\end{equation}
Here, \(\odot\) denotes element-wise multiplication.
The last term corrects for the velocity-divergence contribution to represent the nonconservative convection term in Eq.~\eqref{eq:burgers_momentum}.
Diffusion is represented by a discrete Laplacian \(L_h\) constructed from centered face gradients.

\subsubsection{Time integration and Newton solver}
\label{app:burgers_fom_solver}
Fully implicit backward Euler integration is performed with time step \(\Delta t=0.02\) to the final time \(T=1\).
The residual for field \(q\) at the new time step is
\begin{equation}
 R_q^{n+1}
 =\frac{q^{n+1}-q^n}{\Delta t}
 +C_h(q^{n+1})-\nu L_hq^{n+1}=0,
 \qquad q\in\{u,v\}
 \label{eq:burgers_backward_euler}
\end{equation}
The face velocities in the convection term are also evaluated using \((u,v)\) at the new time step, so both components are unknowns in a coupled nonlinear system.
The previous-time solution is used as the initial guess, and Newton updates are computed using an analytically constructed Jacobian.
Let \(w=[u;v]\) denote the coupled state and \(R(w)\) its residual. The linear system and update at iteration \(k\) are
\begin{equation}
 \frac{\partial R}{\partial w}(w_k)\delta w_k=-R(w_k),
 \qquad w_{k+1}=w_k+\delta w_k
 \label{eq:fom_newton_update}
\end{equation}
The stopping criterion is a residual norm at or below \(10^{-10}\), with a maximum of 20 iterations.

The Newton linear systems on the \(3\times3\) training domain are solved by a sparse direct method.
The \(15\times15\) reference FOM uses ILU-preconditioned GMRES with a relative tolerance of \(10^{-10}\), a restart length of 100, and a maximum of 1000 iterations.
The ILU drop tolerance is \(10^{-4}\), and the fill factor is 10.
The ROM minimizes the same time-discrete residual in reduced coordinates. Its iterative solver settings are given in Appendix~\ref{app:burgers_rom_solver}.

\subsection{Navier--Stokes problem and parameterized forcing}
\label{app:ns_data_generation}
The Navier--Stokes FOM uses Eqs.~\eqref{eq:ns_momentum}--\eqref{eq:ns_incompressibility}.
On the domain \(\Omega=[0,L_x]\times[0,L_y]\), velocity and pressure are periodic in both coordinate directions:
\begin{align}
 \mathbf u(0,y,t)&=\mathbf u(L_x,y,t),&
 \mathbf u(x,0,t)&=\mathbf u(x,L_y,t),\\
 p(0,y,t)&=p(L_x,y,t),&p(x,0,t)&=p(x,L_y,t).
 \label{eq:ns_periodic_bc}
\end{align}
The FOM pressure reference is set by requiring a zero global mean:
\begin{equation}
 \int_\Omega p(\mathbf x,t)\,d\mathbf x=0
 \label{eq:pressure_mean_zero}
\end{equation}
In the ROM, this mean condition is imposed weakly through the penalty in Eq.~\eqref{eq:pressure_penalty_main} of the main text. The formulation is described in Appendix~\ref{app:pressure_penalty}.

\subsubsection{Spatial distribution and time dependence of the forcing}
\label{app:ns_forcing}
There is one forcing source per unit subdomain, so the number of sources \(n_s\) equals the number of subdomains \(n_\Omega\).
The combination of forcing directions is denoted by \(\boldsymbol\mu=(\mu_1,\ldots,\mu_{n_s})\) and sampled using a Sobol sequence.
For a given forcing realization, the total forcing is
\begin{equation}
 \mathbf f(\mathbf x,t;\boldsymbol\mu)
 =\sum_{s=1}^{n_s}a(t)
 \begin{bmatrix}\cos\theta_s(\mu_s)\\\sin\theta_s(\mu_s)\end{bmatrix}
 \chi_s(\mathbf x)
 \label{eq:forcing_param}
\end{equation}
Here, \(\chi_s\) is the spatial weight of source \(s\), and \(a(t)\) is a common time-decay function for all sources.
We choose the decay function
\begin{equation}
 a(t)=\max\{2-t/2,0\}
 \label{eq:forcing_time_decay}
\end{equation}
so that it decreases linearly in time while remaining nonnegative.

\begin{figure}[pos=htbp]
\centering
\includegraphics[width=0.25\textwidth]{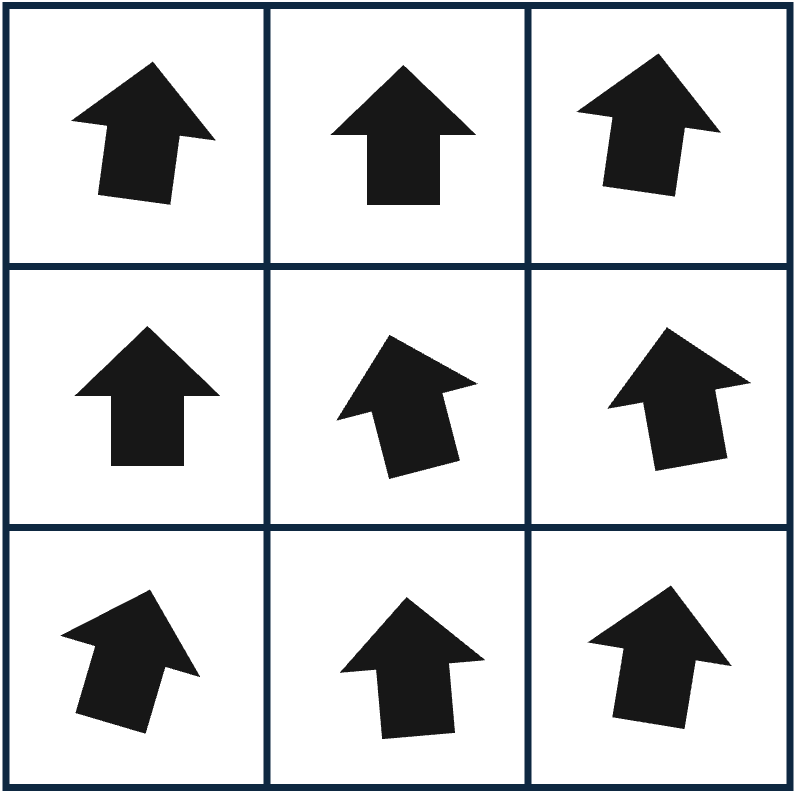}
\caption{Direction parameterization of the localized Navier--Stokes forcing.
Each unit subdomain contains one Gaussian source with a fixed center and width.
Its direction is \(\theta_s=\pi/2+\mu_s\), with \(\mu_s\in[-10^\circ,10^\circ]\) sampled using a Sobol sequence.
The spatial profile is given by \(\chi_s\) in Eq.~\eqref{eq:forcing_spatial_weight}.}
\label{fig:force_param}
\end{figure}

Let \(\boldsymbol\xi_s(\mathbf x)\in[0,1]^2\) denote the normalized local coordinates within subdomain \(\Omega_s\). The spatial weight is defined as
\begin{equation}
 \chi_s(\mathbf x)=
 \begin{cases}
 \displaystyle\exp\!\left[-\frac{\|\boldsymbol\xi_s(\mathbf x)-\mathbf c_s\|_2^2}
 {2\sigma_f^2}\right],&\mathbf x\in\Omega_s,\\[1ex]
 0,&\mathbf x\notin\Omega_s
 \end{cases}
 \label{eq:forcing_spatial_weight}
\end{equation}
The source center and width are fixed at \(\mathbf c_s=(1/2,1/2)\) and \(\sigma_f=0.06\), respectively.
Each source therefore has a Gaussian profile that is largest at the subdomain center and decreases away from it.
The direction is varied within \(\pm10^\circ\) of the upward direction according to Eq.~\eqref{eq:forcing_angle_param}.
Figure~\ref{fig:force_param} illustrates this direction parameterization.

\subsection{Discretization and full-order solver for Navier--Stokes}
\label{app:ns_discretization}
Letting \(n_{\mathrm{sub},x}\) and \(n_{\mathrm{sub},y}\) denote the numbers of subdomains, we set the grid-point counts as
\begin{equation}
 n_x=40n_{\mathrm{sub},x},\qquad n_y=40n_{\mathrm{sub},y},\qquad
 \label{eq:ns_grid_configuration}
\end{equation}
and global domain dimensions as
\begin{equation}
 L_x=n_{\mathrm{sub},x},\qquad L_y=n_{\mathrm{sub},y}.
 \label{eq:ns_domain_configuration}
\end{equation}
Thus, the \(3\times3\) training configuration uses a \(120\times120\) grid, and the \(15\times15\) scale-up configuration uses a \(600\times600\) grid, with \(m=1600\) local degrees of freedom per field.
The grid spacings are \(h_x=L_x/(n_x-1)\) and \(h_y=L_y/(n_y-1)\).

Time integration uses fully implicit backward Euler with \(\Delta t=0.1\).
The viscosity is \(\nu=10^{-4}\) and the density is \(\rho=1\). 50 time steps are taken to reach \(T=5\).
Let \(D_x\) and \(D_y\) denote the discrete derivative operators in the two directions and \(L\) the discrete Laplacian.
Pressure stabilization is introduced by adding \(-\tau Lp\) to the continuity residual, with
\begin{equation}
 \tau=10^{-4}h_x^2
 \label{eq:ns_pressure_stabilization}
\end{equation}

Let \((u,v,p)\) denote the unknowns at the new time step and \((u^n,v^n)\) the previous-time velocities. We define
\begin{equation}
 A=I/\Delta t-\nu L,\qquad
 \mathcal C(u,v)=\operatorname{diag}(u)D_x+\operatorname{diag}(v)D_y
\end{equation}
The residuals are then given by
\begin{align}
 R_u&=Au+\mathcal C(u,v)u+D_xp-u^n/\Delta t-f_x,\\
 R_v&=Av+\mathcal C(u,v)v+D_yp-v^n/\Delta t-f_y,\\
 R_c&=D_xu+D_yv-\tau Lp
 \label{eq:ns_discrete_residual}
\end{align}
Here, \(f_x\) and \(f_y\) are also evaluated at the new time step.
Each subdomain residual is evaluated by selecting its owned rows from this global residual and using the neighboring states needed for the derivative and convection operators.

In the full-order computation, \(w=[u;v;p]\) is treated as a single unknown vector, and the Newton update in Eq.~\eqref{eq:fom_newton_update} is applied.
A sparse direct method is used for the Newton linear systems on the \(3\times3\) training domain, while GMRES is used for the \(15\times15\) linear systems, which include velocity--pressure coupling and nonsymmetric convection terms.
The residual Jacobian is given in Appendix~\ref{app:ns_discrete_jacobian}, and the ROM mean-pressure treatment and iterative solver are described in Appendices~\ref{app:pressure_penalty} and~\ref{app:ns_rom_solver}.

\subsection{Local snapshot collection}
\label{app:data_split}
For both problems, data generation begins with small-domain FOM integration for each parameter case.
Let \(N_{\mathrm{case}}\) be the number of cases, \(\boldsymbol\mu_j\) the parameters for case \(j\), and \(N_{\mathrm{snap}}\) the number of saved snapshots per trajectory.
The \(\ell\)-th global snapshot of field \(q\) is denoted by \(q^\ell(\boldsymbol\mu_j)\).
Using a restriction operator \(E_i\) to extract the state owned by subdomain \(i\), a local sample is given by
\begin{equation}
 x_{q,i}^\ell(\boldsymbol\mu_j)
 =E_iq^\ell(\boldsymbol\mu_j)\in\mathbb R^m
 \label{eq:local_snapshot}
\end{equation}
Pooling samples over all cases, saved time instances, and subdomains of the \(3\times3\) configuration gives the field-wise dataset
\begin{equation}
 \mathcal X_q=
 \left\{x_{q,i}^\ell(\boldsymbol\mu_j):
 j=1,\ldots,N_{\mathrm{case}},\quad
 \ell=1,\ldots,N_{\mathrm{snap}},\quad i=1,\ldots,9\right\}
 \label{eq:field_dataset}
\end{equation}
The total number of local samples per field is therefore
\begin{equation}
 N_{\mathrm{sample}}=N_{\mathrm{case}}N_{\mathrm{snap}}\times9
\end{equation}
This procedure extracts local states at different times and locations from each global trajectory and uses them to train a common representation.

We distinguish between the small-domain data used for representation learning and the large-domain data used for scale-up prediction evaluation.
The AE is trained on the training split, while the validation split is used to select the learning rate and stopping point.
The POD basis is constructed from the same training data.
Local representation accuracy is assessed on small-domain test cases excluded from training and model selection, whereas scale-up predictions are evaluated against separate large-domain FOM trajectories.

\subsubsection{Data split for Burgers}
\label{sec:burgers_data_split}
For Burgers, independent trajectories are generated by initial-condition case: 1600 for training, 400 for validation, and 200 for local reconstruction testing.
All time and subdomain samples from one case belong to the same split.
Each trajectory contains 51 saved states, including \(t=0\), giving 734400 training, 183600 validation, and 91800 test local samples per field.
The snapshot matrix for local representation evaluation is therefore \(X_q\in\mathbb R^{2500\times91800}\).
Validation loss is used for learning-rate scheduling and early stopping. The test split is not used for model selection.
Scale-up evaluation uses 10 separately generated \(15\times15\) initial-condition cases.

\subsubsection{Data split for Navier--Stokes}
\label{sec:ns_data_split}
For Navier--Stokes, local snapshots generated from 1000 small-domain forcing cases are randomly split into training and validation sets with a training fraction of 0.8 and random seed 2025.
Thus, unlike Burgers, the validation split is performed at the local-snapshot level rather than by entire case.
Validation loss is used for learning-rate scheduling and early stopping.
A separate set of 200 small-domain forcing cases is generated exclusively for local reconstruction testing. These cases are excluded from training, validation, and model selection.
The 10 large-domain test cases are also generated by separate FOM integrations, and the AEs and POD bases learned on the small domain are applied without retraining.

\subsection{Summary of the data-generation procedure}
The full procedure consists of parameter sampling, small-domain FOM integration, storage of field-wise time snapshots, extraction of non-overlapping subdomain states, and representation learning using the specified training and validation splits.
Initial-condition coefficients are sampled for Burgers, and forcing directions are sampled for Navier--Stokes.
After their reconstruction performance is evaluated on separate small-domain test data, the learned representations are used for ROM predictions against large-domain FOM reference trajectories.
Separating training-data construction from these two evaluations allows separate assessment of local representation capability and scale-up performance without retraining.

\clearpage
\section{Autoencoder Architecture and Training}
\label{sec:model_architecture}
\label{app:ae_details}
This appendix describes the architecture and training of the autoencoders used to construct the local nonlinear manifolds.
Each physical field uses an independent encoder--decoder pair: the encoder is a fully connected multilayer perceptron (MLP), and the decoder has one hidden layer followed by a sparse output layer that reflects spatial connectivity.
Burgers and Navier--Stokes use the same field-wise construction principle, but differ in encoder depth, normalization, and network size.

\subsection{Field-wise local autoencoders}
For a local snapshot \(x_q\in\mathbb R^m\) of field \(q\), the encoder \(h_q\) produces reduced coordinates \(z_q\in\mathbb R^{n_q}\), and the decoder \(g_q\) reconstructs the local state from them.
The full autoencoder is \(\mathcal A_q=g_q\circ h_q\), as in Eq.~\eqref{eq:experiment_field_sets} of the main text.
The encoder and decoder are optimized jointly during training, but in the online residual minimization problem the decoder parameters are fixed and \(z_q\) is treated as the unknown.

To reconstruct each subdomain state, the field-wise decoders are concatenated according to Eq.~\eqref{eq:fieldwise_block_decoder}.
The total local reduced dimension is \(n_u+n_v\) for Burgers and \(n_u+n_v+n_p\) for Navier--Stokes.
This construction separates the learning problem by physical field.
The field index \(q\) is omitted below when describing the network architecture.

\subsection{Encoder and input normalization}
\label{app:normalization}
Let \(x_{\mathrm{ref}}\) be the input reference offset and \(s\) the scale. The normalized input is \((x-x_{\mathrm{ref}})\oslash s\).
Here, \(\oslash\) denotes element-wise division.
An encoder with \(L_h\) hidden layers is written as
\begin{align}
 a_0&=(x-x_{\mathrm{ref}})\oslash s,\\
 a_\ell&=\sigma(W_{h,\ell}a_{\ell-1}+b_{h,\ell}),
 \qquad\ell=1,\ldots,L_h,\\
 z&=W_{h,L_h+1}a_{L_h}+b_{h,L_h+1}
 \label{eq:encoder_general}
\end{align}
We use the SiLU activation,
\begin{equation}
 \sigma(t)=t\operatorname{sigmoid}(t)
\end{equation}
Bias terms are included in the hidden layers, and no activation is applied at the final latent output layer.
The final encoder layer is trained without bias, so \(b_{h,L_h+1}=0\).

For Burgers, we use \(L_h=3\) with hidden widths \((300,200,100)\).
No normalization is applied, so \(x_{\mathrm{ref}}=0\) and \(s=\mathbf1\).
For Navier--Stokes, \(L_h=1\), and the inputs are standardized.
The one-hidden-layer encoder can be written as
\begin{equation}
 z=W_{h,2}\sigma\!\left(W_{h,1}\bigl((x-x_{\mathrm{ref}})\oslash s\bigr)+b_{h,1}\right)
 \label{eq:encoder}
\end{equation}
Input normalization is thus included at the start of the encoder and paired with the inverse transformation in the physical-state decoder defined below.

\subsection{Sparse shallow decoder and physical-state reconstruction}
\label{sec:decoder_jacobian}
Let \(\eta\in\mathbb R^{M_2}\) denote the decoder hidden features. The network output is
\begin{align}
 \eta&=\sigma(W_{d,1}z+b_{d,1}),\\
 d(z)&=W_{d,2}^{\mathcal M}\eta+b_{d,2}
 \label{eq:decoder}
\end{align}
The output weight matrix \(W_{d,2}^{\mathcal M}\in\mathbb R^{m\times M_2}\) has connections restricted by the binary mask \(\mathcal M\).
Specifically, we set
\begin{equation}
 (W_{d,2}^{\mathcal M})_{i\ell}=0
 \quad\text{if}\quad\mathcal M_{i\ell}=0
 \label{eq:masked_weights}
\end{equation}
so that only the prescribed connections are trained.
The hidden-layer bias is retained, while the output layer is trained without bias, so \(b_{d,2}=0\).

When the network output \(d(z)\) represents a normalized state, the inverse transformation must be included to reconstruct the physical state.
The decoder substituted into the governing residual is
\begin{equation}
 g(z)=x_{\mathrm{ref}}+s\odot d(z)
 \label{eq:physical_decoder_interface}
\end{equation}
Thus, even without a bias in the network output layer, physical-state reconstruction may include the reference offset \(x_{\mathrm{ref}}\).
For Burgers, where normalization is not used, this relation reduces to \(g(z)=d(z)\).

The decoder Jacobian required for Gauss--Newton assembly is obtained by differentiating this map.
For the network output, we have
\begin{equation}
 D(z):=\frac{\partial d}{\partial z}
 =W_{d,2}^{\mathcal M}
 \operatorname{diag}\!\left(\sigma'(W_{d,1}z+b_{d,1})\right)W_{d,1}
 \label{eq:decoder_jacobian}
\end{equation}
and for physical-state reconstruction,
\begin{equation}
 \frac{\partial g}{\partial z}=\operatorname{diag}(s)D(z)
 \label{eq:physical_decoder_jacobian}
\end{equation}
The sparse output weights allow these computations to use sparse--dense matrix products and to evaluate selected output rows efficiently \citep{kim2022fast-c42}.
However, if \(W_{d,1}\) is dense, the resulting \(D(z)\in\mathbb R^{m\times n_q}\) can also be dense.
The output-layer sparsity, the global block-diagonal structure formed by concatenating subdomain decoders, and the local block structure of the PDE residual occur at different levels.
Future integration with residual sampling would allow only the hidden units connected to selected outputs to be evaluated, but the present experiments use the full residual.

\subsection{Nearest-neighbor-based output mask}
\label{app:decoder_mask}
The output mask is constructed from the \(k\)-nearest-neighbor connectivity of the local grid.
Let the subdomain grid points be \(\{y_i\}_{i=1}^m\), with \(y_i\in\mathbb R^2\), and define the neighbor set of output point \(i\) as
\begin{equation}
 \mathcal N_k(i)=\operatorname{kNN}(y_i)
\end{equation}
For each grid point \(j\), a hidden-index block is assigned using block size \(b\) and block stride \(db\):
\begin{equation}
 B_j=\{(j-1)db+1,\ldots,(j-1)db+b\}
\end{equation}
The total decoder hidden width is then
\begin{equation}
 M_2=b+db(m-1)
 \label{eq:decoder_hidden_width}
\end{equation}

\begin{figure}[pos=htbp]
\centering
\includegraphics[width=0.72\textwidth]{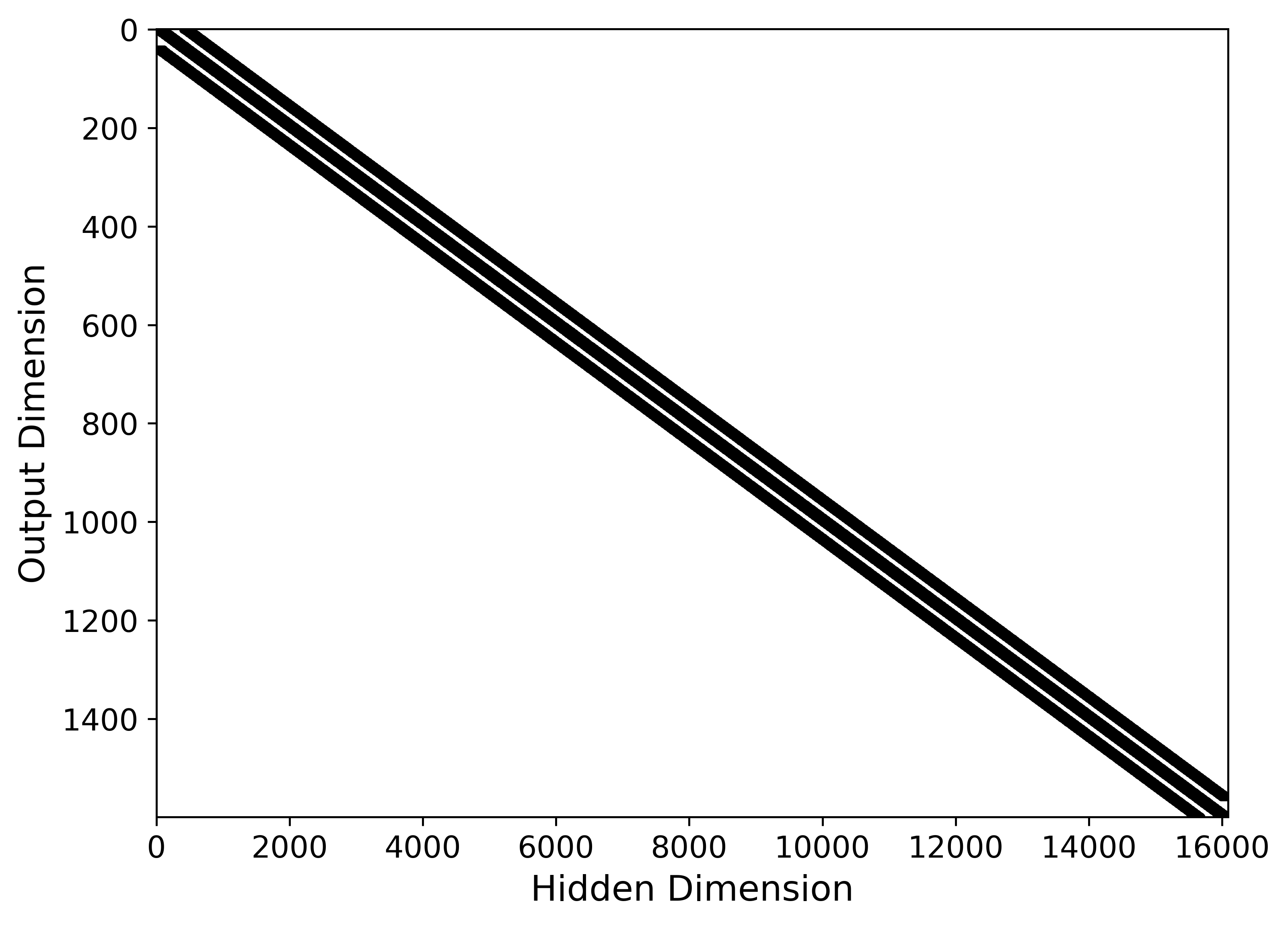}
\caption{Example of the kNN-based sparse decoder output mask with \(b=100\), \(db=10\), and \(k=5\).
The binary pattern specifies connections between output degrees of freedom and decoder hidden units.
Each output is connected only to hidden blocks associated with its selected local-grid neighbors.}
\label{fig:sparse_mask_knn}
\end{figure}

The hidden indices connected to output point \(i\) are obtained by collecting the hidden blocks of its neighboring grid points:
\begin{equation}
 \mathcal H_i=\bigcup_{j\in\mathcal N_k(i)}B_j
 \label{eq:decoder_support}
\end{equation}
The mask is given by
\begin{equation}
 \mathcal M_{i\ell}=1
 \quad\Longleftrightarrow\quad\ell\in\mathcal H_i
 \label{eq:sparse_mask}
\end{equation}
Thus, each output degree of freedom is connected only to hidden blocks associated with nearby points on the local grid.
Figure~\ref{fig:sparse_mask_knn} shows an example mask for \(b=100\), \(db=10\), and \(k=5\).

\subsection{Problem-specific architecture configurations}
\label{app:architecture_settings}
Table~\ref{tab:burgers_ae_architecture} summarizes the field-wise autoencoder architecture for Burgers.
The two velocity components are trained independently using the same architecture, with fixed encoder widths and output-mask parameters across latent dimensions.
Substituting the local dimension \(m=2500\), \(b=10\), and \(db=1\) into Eq.~\eqref{eq:decoder_hidden_width} gives a decoder hidden width of \(M_2=2509\).

\begin{table}[pos=htbp]
\centering
\small
\caption{Field-wise autoencoder architecture for the Burgers experiments.
The \(u\) and \(v\) autoencoders are trained independently with the same architecture.
Encoder widths and decoder-mask parameters are fixed across the latent dimensions.}
\label{tab:burgers_ae_architecture}
\begin{tabular}{lc}
\toprule
Setting & Value \\
\midrule
Local input/output dimension \(m\) & \(2500\) \\
Encoder hidden width & \(300,\ 200,\ 100\) \\
Latent dimension per field \(n_q\) & \(20,25,30,35,40,45,50\) \\
Decoder hidden width \(M_2\) & \(2509\) \\
Block size \(b\) & \(10\) \\
Block stride \(db\) & \(1\) \\
kNN neighbors \(k\) & \(5\) \\
\bottomrule
\end{tabular}
\end{table}

Table~\ref{tab:ae_architecture} shows the Navier--Stokes architecture configurations.
In this case, the encoder hidden width and mask parameters \(b,db\) increase together with the latent dimension.
The performance comparisons in the main text therefore include differences in network architecture as well as latent dimension.
The decoder hidden width for each configuration is determined by \(m=1600\) and the tabulated values of \(b,db\).

\begin{table}[pos=htbp]
\centering
\small
\caption{Autoencoder architecture configurations for the Navier--Stokes experiments.
The latent dimension \(r\) is specified per field.
The encoder width and decoder-mask parameters vary with \(r\). The decoder hidden width is given by Eq.~\eqref{eq:decoder_hidden_width}.}
\label{tab:ae_architecture}
\begin{tabular}{ccccc}
\toprule
Field-wise dimension \(r\) & Encoder hidden width & \(b\) & \(db\) & \(k\) \\
\midrule
5  & 400  & 100 & 10 & 5 \\
7  & 560  & 140 & 14 & 5 \\
10 & 800  & 200 & 20 & 5 \\
\bottomrule
\end{tabular}
\end{table}

\subsection{Training and model selection}
\label{app:training_hparams}
Each field-wise autoencoder is trained using Adam to minimize the reconstruction MSE in Eq.~\eqref{eq:fieldwise_ae_training_loss}.
Burgers uses case-wise training and validation splits, whereas Navier--Stokes uses a random split of small-domain local snapshots.
In both cases, validation loss is used for learning-rate scheduling and early stopping, and the checkpoint with the lowest validation loss is selected as the final model.
Neither the separate local test data nor the large-domain evaluation data are used in this selection.

The Burgers training settings are summarized in Table~\ref{tab:burgers_ae_training_hparams}.
No normalization is applied, and the batch size and initial learning rate are 918 and \(5\times10^{-4}\), respectively.
The learning-rate reduction factor is 0.5, and training runs for at most 15000 epochs, subject to the specified early-stopping criterion.
The Navier--Stokes settings are given in Table~\ref{tab:ae_training_hparams}: standardized inputs, a batch size of 450, an initial learning rate of \(10^{-3}\), and a maximum of 10000 epochs.

\begin{table}[pos=htbp]
\centering
\small
\caption{Autoencoder training settings for the Burgers experiments.
All snapshots from one initial-condition case belong to the same training, validation, or test split.}
\label{tab:burgers_ae_training_hparams}
\begin{tabular}{lc}
\toprule
Setting & Value \\
\midrule
Fields & \(u,v\) trained independently \\
Training/validation/test split unit & Initial-condition case \\
Training/validation/test cases & \(1600/400/200\) \\
Normalization & None \\
Batch size & \(918\) \\
Activation & SiLU \\
Loss & MSE \\
Optimizer & Adam \\
Maximum epochs & \(15000\) \\
Initial learning rate & \(5\times10^{-4}\) \\
LR reduction factor & \(0.5\) \\
LR scheduler patience & \(300\) \\
Minimum learning rate & \(10^{-6}\) \\
Early-stopping patience & \(600\) \\
Random seed & \(2025\) \\
\bottomrule
\end{tabular}
\end{table}

\begin{table}[pos=htbp]
\centering
\small
\caption{Autoencoder training settings for the Navier--Stokes experiments.
Local snapshots from 1000 small-domain cases are randomly split for training and validation.
The separate 200-case representation test set is not used for model selection.}
\label{tab:ae_training_hparams}
\begin{tabular}{lc}
\toprule
Setting & Value \\
\midrule
Fields & \(u,v,p\) trained independently \\
Validation split & Random local-snapshot split \\
Training fraction & \(0.8\) \\
Batch size & 450 \\
Normalization & Standardization \\
Activation & SiLU \\
Encoder layers & Fully connected \\
Decoder output layer & Sparse mask \\
Mask type & k-nearest neighbors (kNN) \\
Maximum epochs & 10000 \\
Initial learning rate & \(10^{-3}\) \\
LR scheduler patience & 300 \\
Early-stopping patience & 600 \\
Optimizer & Adam \\
\bottomrule
\end{tabular}
\end{table}

\subsection{Representative training histories}
\label{app:training_history}
\begin{figure}[pos=htbp]
\centering
\includegraphics[width=0.67\textwidth]{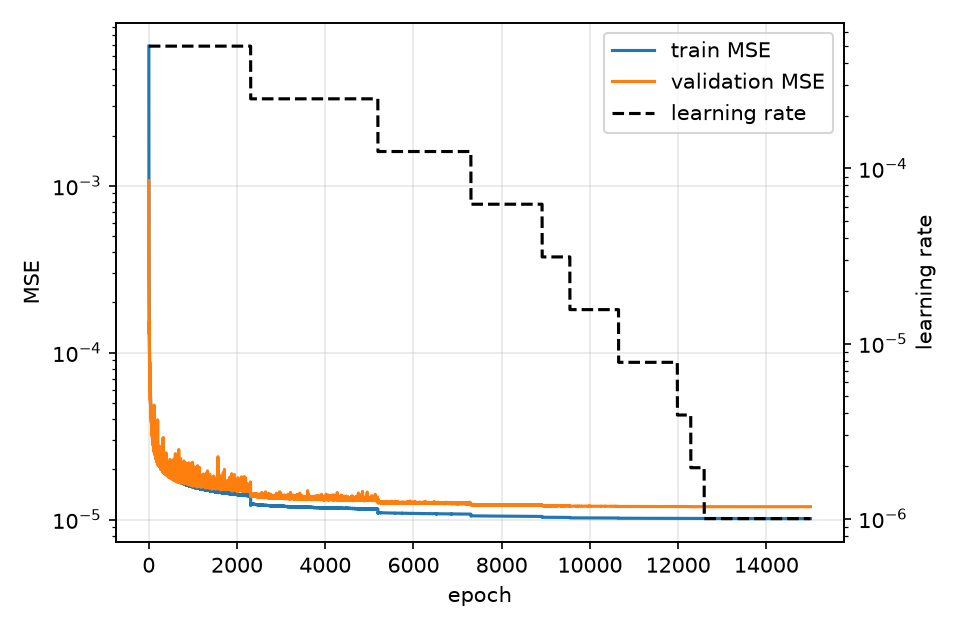}\\
\includegraphics[width=0.67\textwidth]{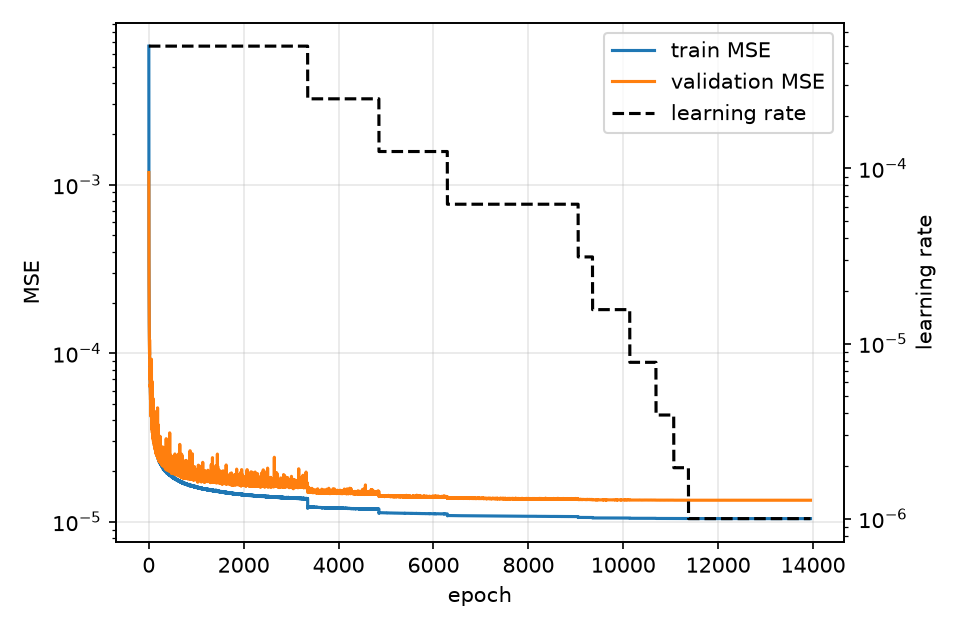}
\caption{Representative training histories for the Burgers autoencoders.
The top and bottom panels correspond to \(u\) and \(v\), respectively, and show training MSE, validation MSE, and learning rate.
The latent dimension is 40 per field, with mask parameters \(b=10\), \(db=1\), and \(k=5\).
This training example differs from the AE \((50,50)\) model used in the main-text speed-field comparison.}
\label{fig:burgers_ae_training_history}
\end{figure}

\begin{figure}[pos=p]
\centering
\includegraphics[width=0.60\textwidth]{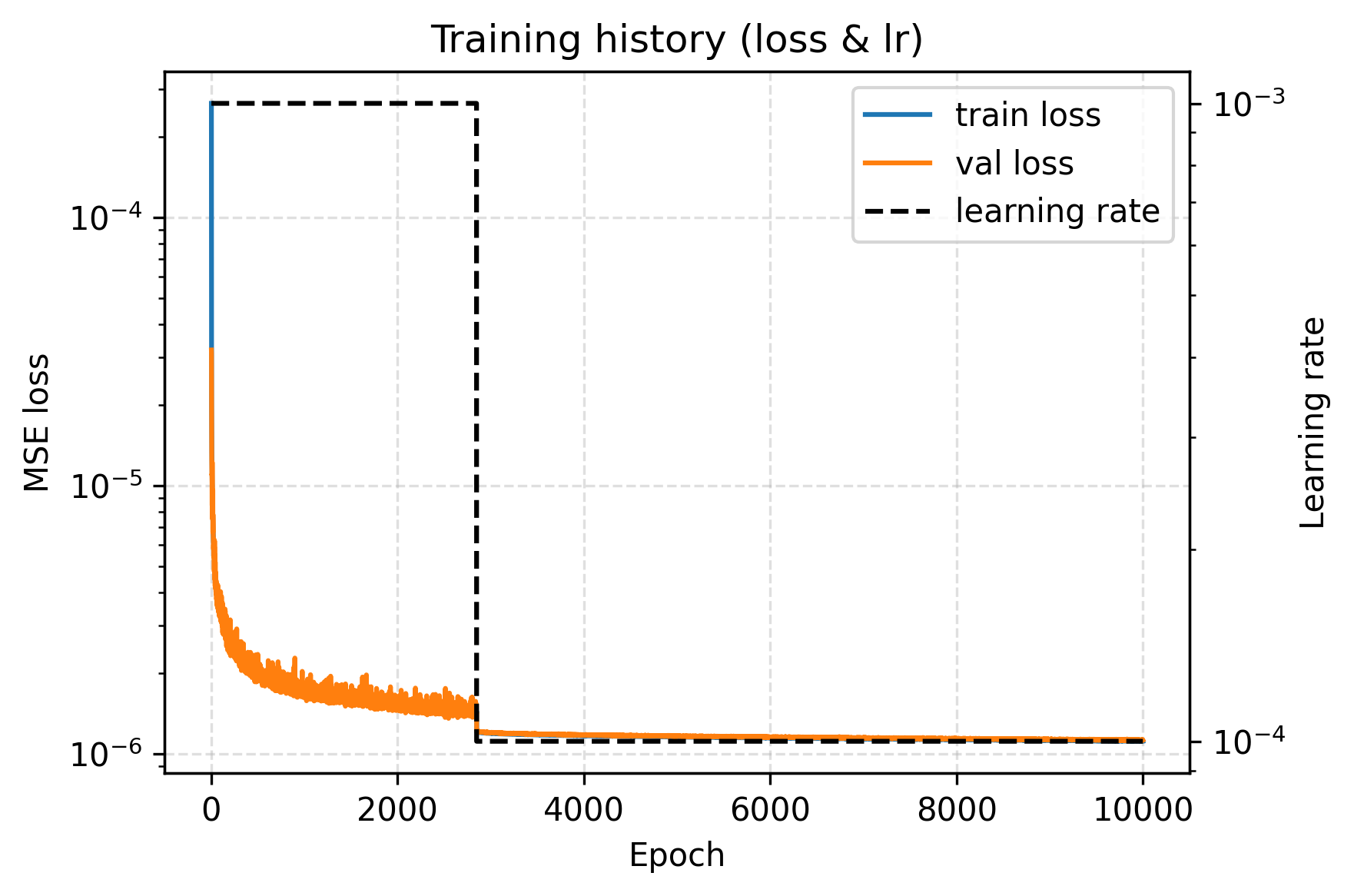}\\
\includegraphics[width=0.60\textwidth]{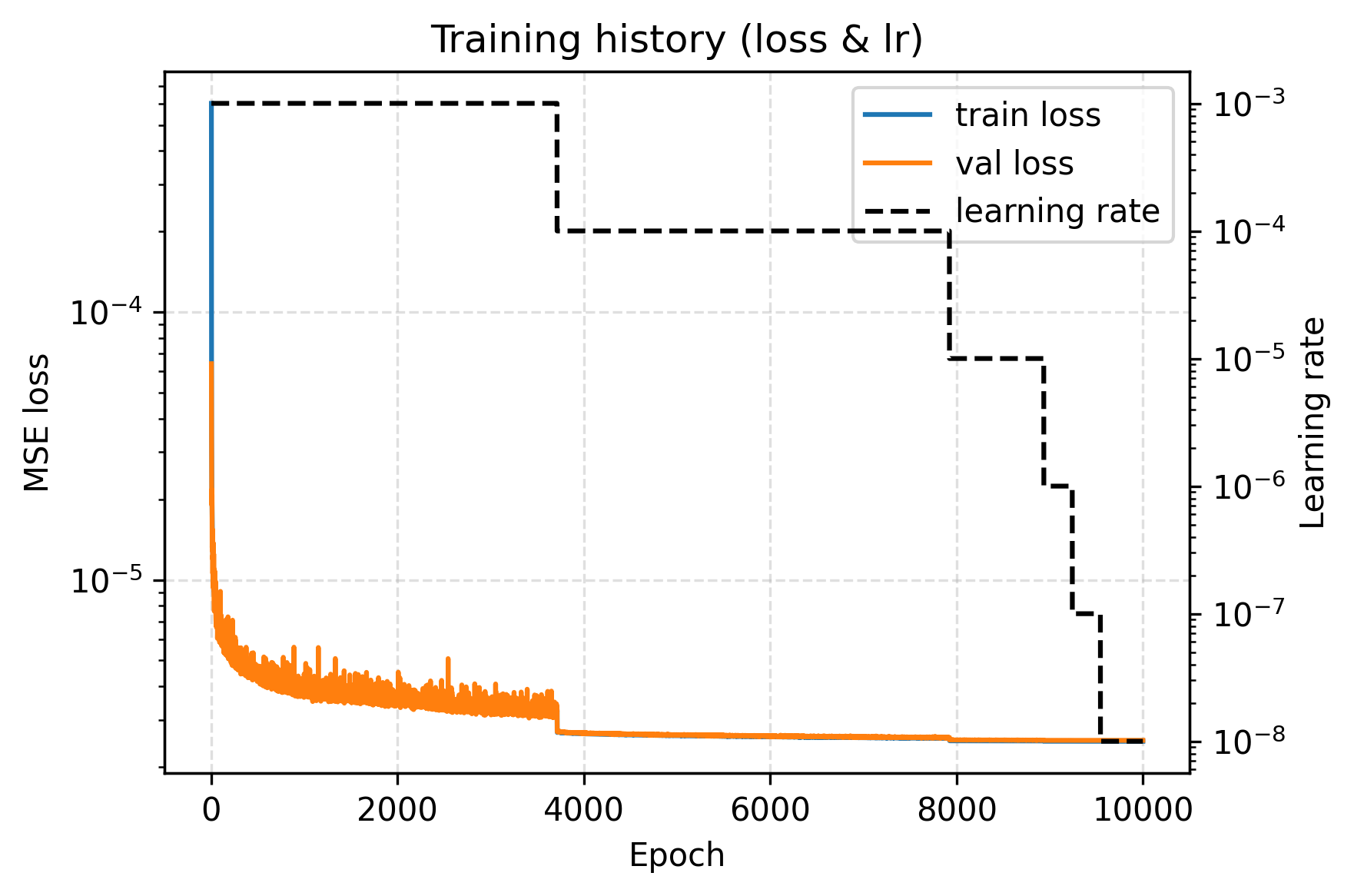}\\
\includegraphics[width=0.60\textwidth]{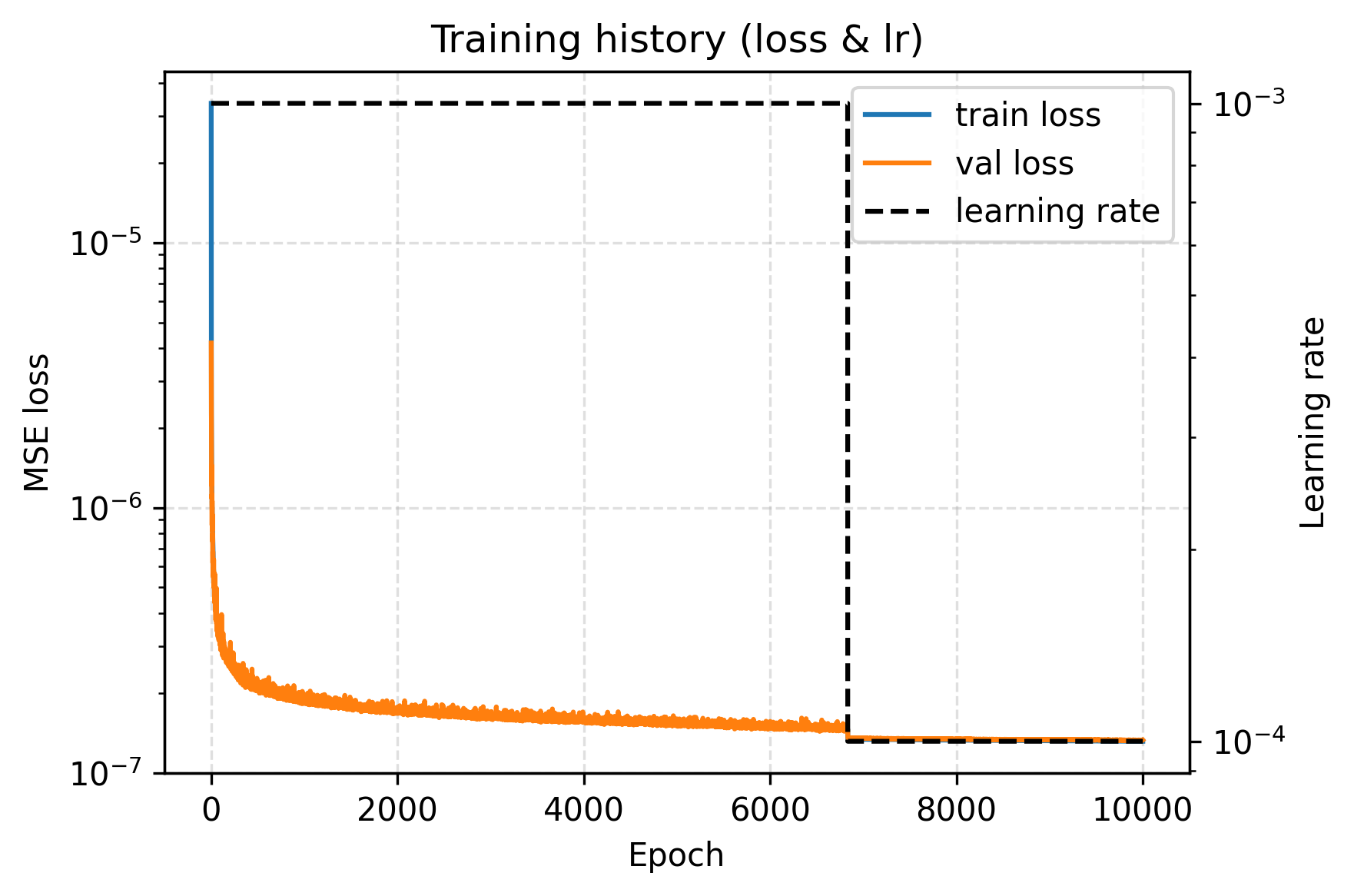}
\caption{Representative field-wise training histories for the Navier--Stokes autoencoders with latent dimension 7 per field.
The top, middle, and bottom panels correspond to \(u\), \(v\), and \(p\), respectively.
Solid curves show training and validation MSE, and dashed curves show the learning rate.}
\label{fig:ae_training_history}
\end{figure}

Figure~\ref{fig:burgers_ae_training_history} shows Burgers training histories for a field-wise dimension of 40.
The mask parameters are \(b=10\), \(db=1\), and \(k=5\), and the plots show training MSE, validation MSE, and learning-rate changes for each field.
This is a representative training configuration and differs from the AE \((50,50)\) model used in the main-text speed-field comparison.
Figure~\ref{fig:ae_training_history} shows the corresponding training histories for the Navier--Stokes models with a field-wise dimension of 7.
Separate from these learning curves, independent local test errors and large-domain time-integration performance are evaluated in the respective experiment sections of the main text.

\clearpage
\section{Numerical Solver Details}
\label{app:rom_solver}
This appendix describes the preconditioning and stopping criteria used to solve the reduced Gauss--Newton systems in the main text.
It also presents the Navier--Stokes discrete Jacobian and the contributions of the global mean-pressure term to the reduced system.
Here, \(H\) and \(q\) denote the Gauss--Newton matrix and gradient of the objective actually minimized, including the mean-pressure contribution for Navier--Stokes.

\subsection{Subdomain-wise preconditioning and linear systems}
\label{app:rom_preconditioning}
At the current Gauss--Newton iteration, the linear system is
\begin{equation}
 H_k\Delta\hat x^{(k)}=-q_k
 \label{eq:numerical_reduced_pcg_system}
\end{equation}
Following the subdomain-wise ordering of the reduced coordinates, let \(H_{k,ij}\) denote the matrix block coupling subdomains \(i,j\).
The block-Jacobi preconditioner is defined as
\begin{equation}
 M_k=\operatorname{blkdiag}
 (H_{k,11},H_{k,22},\ldots,H_{k,n_\Omega n_\Omega})
 \label{eq:block_jacobi_preconditioner}
\end{equation}
Each diagonal block includes the reduced coordinates of all fields on that subdomain.
Preconditioning therefore consists of small subdomain-wise block solves, while global coupling is handled through the PCG iterations.

If the columns of the respective Jacobian are linearly independent, the corresponding Gauss--Newton matrix is positive definite.
No separate diagonal regularization term is added to the global Gauss--Newton matrix in either experiment.
Preconditioner-block treatment and PCG tolerances follow the problem-specific settings.

\subsection{Reduced nonlinear solver for Burgers}
\label{app:burgers_rom_solver}
For Burgers, each time step uses at most 20 Gauss--Newton iterations.
Iterations stop when the gradient norm is at or below \(10^{-8}\) or the relative change in residual norm between consecutive iterations is at or below \(10^{-2}\).
The decoder of the local representation remains fixed during these iterations. The updated unknowns are the subdomain-specific reduced coordinates.

Each Gauss--Newton linear system is solved using block-Jacobi-preconditioned PCG.
The PCG relative tolerance is \(10^{-4}\), the absolute tolerance is 0, and the maximum iteration count is 1000.
Neither diagonal regularization of the global matrix nor diagonal shifts of the preconditioner blocks are used, and no adaptive diagonal shifting is applied.
These settings are used in all Burgers AE and POD scale-up experiments in the main text.

\subsection{Discrete Jacobian for Navier--Stokes}
\label{app:ns_discrete_jacobian}
Differentiating the residual in Appendix~\ref{app:ns_discretization} with respect to velocity and pressure gives the block Jacobian
\begin{equation}
 \frac{\partial(R_u,R_v,R_c)}{\partial(u,v,p)}=
 \begin{bmatrix}
 A+\mathcal C+\operatorname{diag}(D_xu)&\operatorname{diag}(D_yu)&D_x\\
 \operatorname{diag}(D_xv)&A+\mathcal C+\operatorname{diag}(D_yv)&D_y\\
 D_x&D_y&-\tau L
 \end{bmatrix}.
 \label{eq:ns_discrete_jacobian}
\end{equation}
Here, \(A=I/\Delta t-\nu L\) and \(\mathcal C=\operatorname{diag}(u)D_x+\operatorname{diag}(v)D_y\), with all variables evaluated at the current iterate.
The diagonal momentum blocks include how changes in the corresponding velocity component affect the advecting velocity, while the off-diagonal velocity blocks include changes in convection due to the other velocity component.
The last row is the derivative of the stabilized continuity equation.

The local Jacobian for subdomain \(i\) consists of the owned residual rows and required neighboring state columns of this global operator.
Multiplication by the decoder Jacobian according to Eq.~\eqref{eq:local_chain_rule} of the main text gives the local Jacobian with respect to reduced coordinates.

\subsection{Treatment of the global mean pressure}
\label{app:pressure_penalty}
The Navier--Stokes mean pressure is treated as a single global condition.
Let \(\mathcal I_i^p\) be the set of pressure degrees of freedom owned by subdomain \(i\), and let \(p_i=g_p(z_{p,i})\) be its reconstructed pressure vector.
The ownership sets are non-overlapping, and the global number of pressure degrees of freedom is
\begin{equation}
 N_p=\sum_{i=1}^{n_\Omega}|\mathcal I_i^p|=n_xn_y
\end{equation}
Counting each owned pressure value exactly once, we define the global mean and its difference from the target mean as
\begin{equation}
 \bar p(\hat x)=\frac1{N_p}\sum_{i=1}^{n_\Omega}\mathbf1_i^Tg_p(z_{p,i}),
 \qquad
 \psi(\hat x)=\bar p(\hat x)-p_\star
 \label{eq:pmean_definition}
\end{equation}
Pressure values needed by neighboring residuals also refer to the same owned state and are therefore not counted twice in the mean.
This definition controls the global mean while allowing different subdomain-wise means.

With \(\beta=\gamma N_p\), the objective in Eq.~\eqref{eq:pressure_penalty_main} can be written as
\begin{equation}
 \mathcal F(\hat x)
 =\frac12\sum_{i=1}^{n_\Omega}\|\tilde r_i(\hat x)\|_2^2
 +\frac\beta2\psi(\hat x)^2
 =\frac12\left\|
 \begin{bmatrix}\tilde r(\hat x)\\\sqrt\beta\,\psi(\hat x)\end{bmatrix}
 \right\|_2^2
 \label{eq:pmean_augmented_objective}
\end{equation}
Here, \(\tilde r=\operatorname{col}_i\tilde r_i\), and the experiments use \(\gamma=1\) and \(p_\star=0\).
The coefficient includes the global number of pressure degrees of freedom, \(N_p\).
This term is part of the residual minimization problem that determines the reduced coordinates, not a postprocessing correction that subtracts the mean from an already computed pressure field.

Define the mean-pressure Jacobian as
\begin{equation}
 C=\frac{\partial\bar p}{\partial\hat x}
 \in\mathbb R^{1\times n_\Omega n_{\mathrm{loc}}}
\end{equation}
Its block with respect to the pressure reduced coordinates of subdomain \(i\) is
\begin{equation}
 C_i^p=\frac1{N_p}\mathbf1_i^TG_{p,i},
 \qquad
 G_{p,i}=\frac{\partial g_p(z_{p,i})}{\partial z_{p,i}}
 \label{eq:pmean_jacobian}
\end{equation}
while the blocks corresponding to velocity reduced coordinates are zero.
Linearizing at the current iterate as \(\psi(\hat x+\Delta\hat x)\approx\psi(\hat x)+C\Delta\hat x\) gives the Gauss--Newton matrix and gradient
\begin{equation}
 H=H_{\mathrm{PDE}}+\beta C^TC,
 \qquad
 q=q_{\mathrm{PDE}}+\beta\psi C^T
 \label{eq:pmean_gn_correction}
\end{equation}

The PDE contributions \(H_{\mathrm{PDE}}\) and \(q_{\mathrm{PDE}}\) are computed by the local assembly in the main text, and the mean-pressure contributions are then added.
The term \(\beta C^TC\) is a rank-one correction coupling the global pressure reduced coordinates and is distinct from the direct local-neighbor coupling of the PDE residual.
Using the augmented Jacobian
\begin{equation}
 J_{\mathrm{aug}}=
 \begin{bmatrix}J\\\sqrt\beta\,C\end{bmatrix}
\end{equation}
we have \(H=J_{\mathrm{aug}}^TJ_{\mathrm{aug}}\), which is positive definite when the augmented Jacobian has linearly independent columns.

Pressure stabilization through \(-\tau Lp\) and the mean-pressure term also play different roles.
The periodic Laplacian vanishes for a constant pressure, so the stabilization term alone does not set a global pressure reference.
The mean-pressure term controls this constant mode. It does not directly enforce subdomain interface compatibility or the local incompressibility residual.
Because a finite penalty coefficient suppresses violations of the mean condition, this treatment is distinct from imposing \(\bar p=0\) as a strong equality constraint.

\subsection{Reduced nonlinear solver for Navier--Stokes}
\label{app:ns_rom_solver}
At most 7 Gauss--Newton iterations are performed at each time step.
Let \(q_k\) denote the gradient, \(\Delta\hat x_k\) the increment, and \(\rho_k=\|\tilde r(\hat x_k)\|_2\) the PDE residual norm at the current iterate. Iteration stops when any of the following conditions is satisfied:
\begin{equation}
 \|q_k\|_2<10^{-8},\qquad
 \|\Delta\hat x_k\|_2<10^{-9},\qquad
 \left|\frac{\rho_{k-1}-\rho_k}{\rho_{k-1}}\right|\le10^{-3}.
 \label{eq:ns_rom_stopping}
\end{equation}
Here, \(q_k\) includes the mean-pressure term, whereas \(\rho_k\) is the norm of the PDE residual alone.
The stopping quantities are evaluated at the state at which the current increment is computed, and the termination decision is made after updating the reduced coordinates.
The full increment is applied without a separate line search; if the maximum iteration count is reached, the current state is used to proceed to the next time step.

The linear system is solved by PCG with a relative tolerance of \(10^{-2}\), an absolute tolerance of 0, and a maximum of 150 iterations.
No diagonal regularization is applied to the global Gauss--Newton matrix, and the preconditioner blocks are unshifted by default.
If Cholesky factorization of a diagonal block \(D_B\) fails, however, the following diagonal shift is added as a factorization safeguard:
\begin{equation}
 \epsilon_B=
 \max\left\{10^{-14},\,
 10^{-12}\frac{\operatorname{tr}(D_B)}{\dim(D_B)}\right\}
 \label{eq:preconditioner_block_shift}
\end{equation}
The shifted block is then factorized again.
This shift applies only to the preconditioner block and does not change the original global linear-system matrix.
If PCG does not reach the tolerance within the maximum iteration count, the available approximate increment is used.

These stopping conditions solve the reduced residual minimization problem with a finite iteration budget.
In particular, a small relative change in residual norm is not equivalent to the exact stationarity condition \(\nabla\mathcal F(\hat x)=0\). Sensitivity of the accuracy--efficiency trade-off to stopping tolerances is outside the scope of these experiments.

\subsection{Scope of runtime measurements}
\label{app:timing_details}
Runtimes cover the time-integration loop and exclude model loading, initial reduced-coordinate computation, and file output.
We compare the wall-clock times of the serial FOM and ROM.
These timings evaluate online speed-up according to Eq.~\eqref{eq:speedup} in the main text. Savings in small-domain data-generation and autoencoder-training costs are not separately quantified.

Errors and speed-ups are computed for each of the 10 test cases for each problem.
The speed-up reported in the tables is therefore the arithmetic mean of the case-wise time ratios, not the ratio of the mean FOM runtime to the mean ROM runtime.
This aggregation is also expressed mathematically in Appendix~\ref{app:metric_interpretation}.

\clearpage
\section{Additional Results and State Visualizations}
\label{app:supporting_results}
This appendix supplements the interpretation of the main-text evaluation metrics and the primitive-variable results for representative cases.
The main text shows representative states using speed for Burgers and vorticity for Navier--Stokes; here, we compare the individual fields used directly in learning and residual evaluation.

\subsection{Interpretation of the all-field metrics}
\label{app:metric_interpretation}
Let \(X_q\) denote the FOM snapshot matrix for each field in a single test case. The space--time relative error computed from the concatenated full state satisfies
\begin{equation}
 \bigl(E_{\mathrm{all}}^{L_2}\bigr)^2
 =\frac{\sum_{q\in\mathcal Q}\|X_q\|_F^2\bigl(E_q^{L_2}\bigr)^2}
 {\sum_{q\in\mathcal Q}\|X_q\|_F^2}
 \label{eq:combined_error_weights}
\end{equation}
This relation follows directly from the Frobenius norm of the concatenated state.
The all-field error therefore reflects a state-magnitude-weighted average of squared field-wise errors. It does not imply that all fields have comparable relative errors.
In particular, the Navier--Stokes metric combines velocity and pressure and is affected by their scales and units, so it is interpreted together with field-wise results.

The metric \(E^{\mathrm{Gmax}}\) is the maximum deviation normalized by the maximum norm over the full FOM trajectory of one case.
Its average over 10 cases is not the worst error across those 10 cases.
Because its normalization denominator differs from that of \(E^{L_2}\), no fixed ordering between the two metrics is assumed.
The two accuracy--efficiency panels in the main text jointly assess accumulated trajectory error and the maximum deviation relative to a global scale.

If \(S_j=T_{\mathrm{FOM},j}/T_{\mathrm{ROM},j}\) denotes the speed-up for case \(j\), the mean shown in the tables and figures is
\begin{equation}
 \overline S=\frac1{10}\sum_{j=1}^{10}S_j
 \label{eq:mean_speedup}
\end{equation}
Errors are also computed case-wise and averaged in the same way.
The range bars in the figures indicate the minimum and maximum observed case-wise values.

\subsection{Pointwise comparisons of primitive variables}
\label{app:scaleup_fields}
The field-wise figures show the same representative test cases as the speed and vorticity figures in the main text, not mean fields over the 10 cases.
The comparisons use AE \((50,50)\) and POD \((65,65)\) for Burgers, and AE \((10,10,10)\) and POD \((46,42,37)\) for Navier--Stokes.
The pointwise absolute error for a primitive variable \(q\) is defined as
\begin{equation}
 e_q^{\mathrm{ROM}}(\mathbf x,t)
 =|q_{\mathrm{ROM}}(\mathbf x,t)-q_{\mathrm{FOM}}(\mathbf x,t)|,
 \qquad q\in\mathcal Q
 \label{eq:primitive_visual_error}
\end{equation}
The top row of each figure compares the FOM, DD-NM-ROM, and DD-LS-ROM states, while the bottom row compares the two ROM errors.
State panels for the same physical quantity share a color scale, and the two error panels share an error color scale.

\subsubsection{Velocity components for Burgers}
\label{app:burgers_scaleup_fields}
\begin{figure}[pos=p]
\centering
\includegraphics[width=0.98\textwidth]{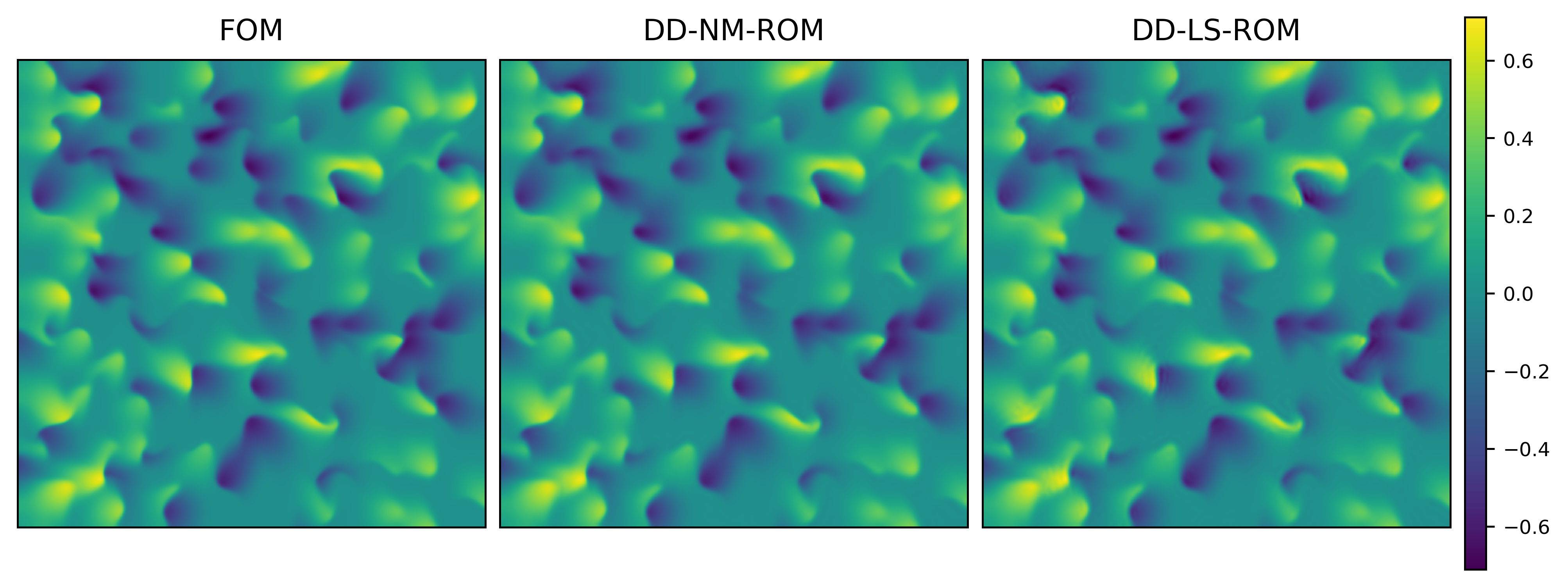}\\
\includegraphics[width=0.86\textwidth]{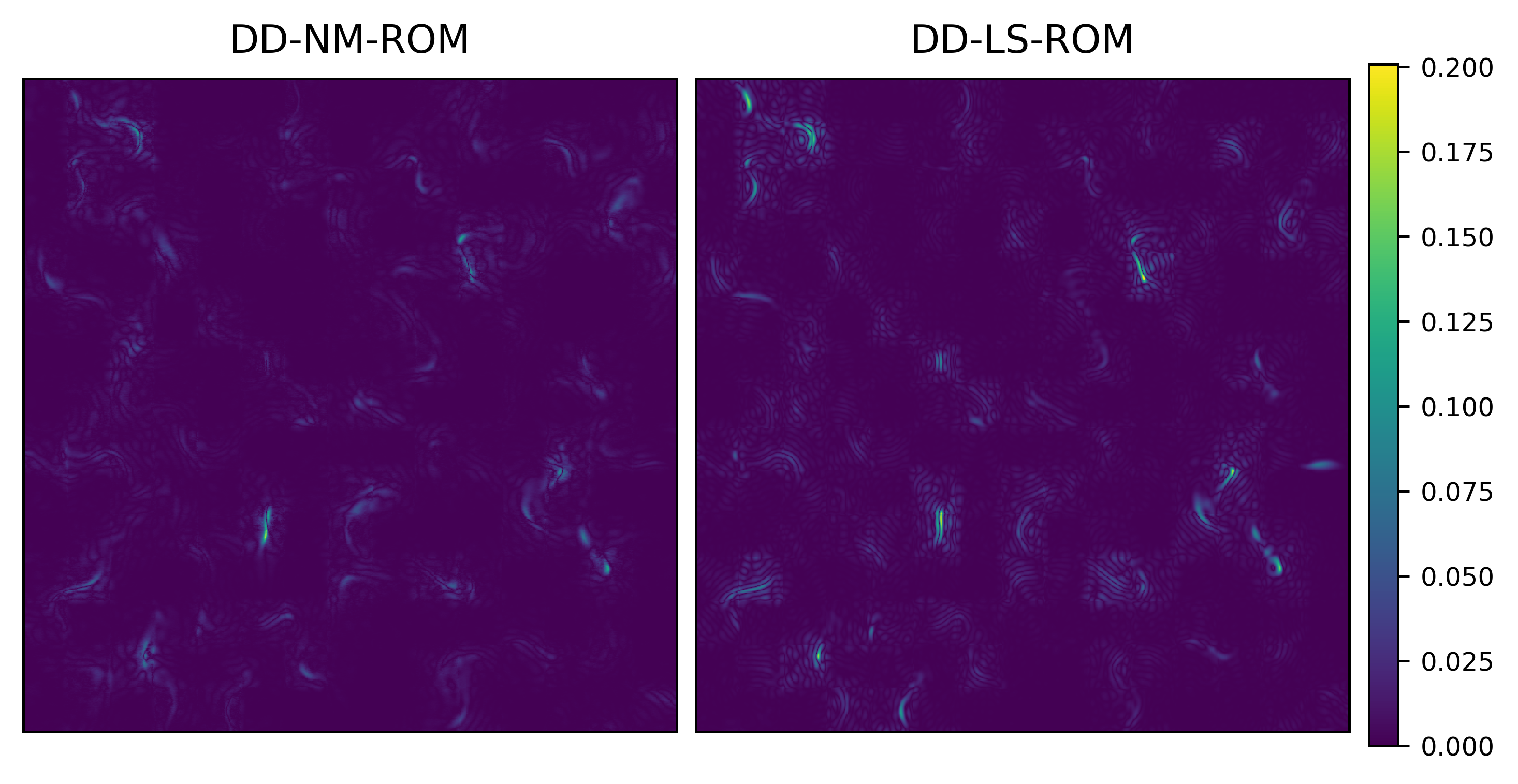}
\caption{Burgers \(u\)-velocity field and pointwise absolute error at \(t=1\) for the representative case used in Figure~\ref{fig:burgers_speed_solution}.
The top row shows the FOM, DD-NM-ROM with AE \((50,50)\), and DD-LS-ROM with POD \((65,65)\), from left to right.
The bottom row shows the corresponding DD-NM-ROM and DD-LS-ROM errors defined in Eq.~\eqref{eq:primitive_visual_error}.
A common color scale is used within each row.}
\label{fig:burgers_app_u_solution}
\label{fig:burgers_app_u_error}
\end{figure}

\begin{figure}[pos=p]
\centering
\includegraphics[width=0.98\textwidth]{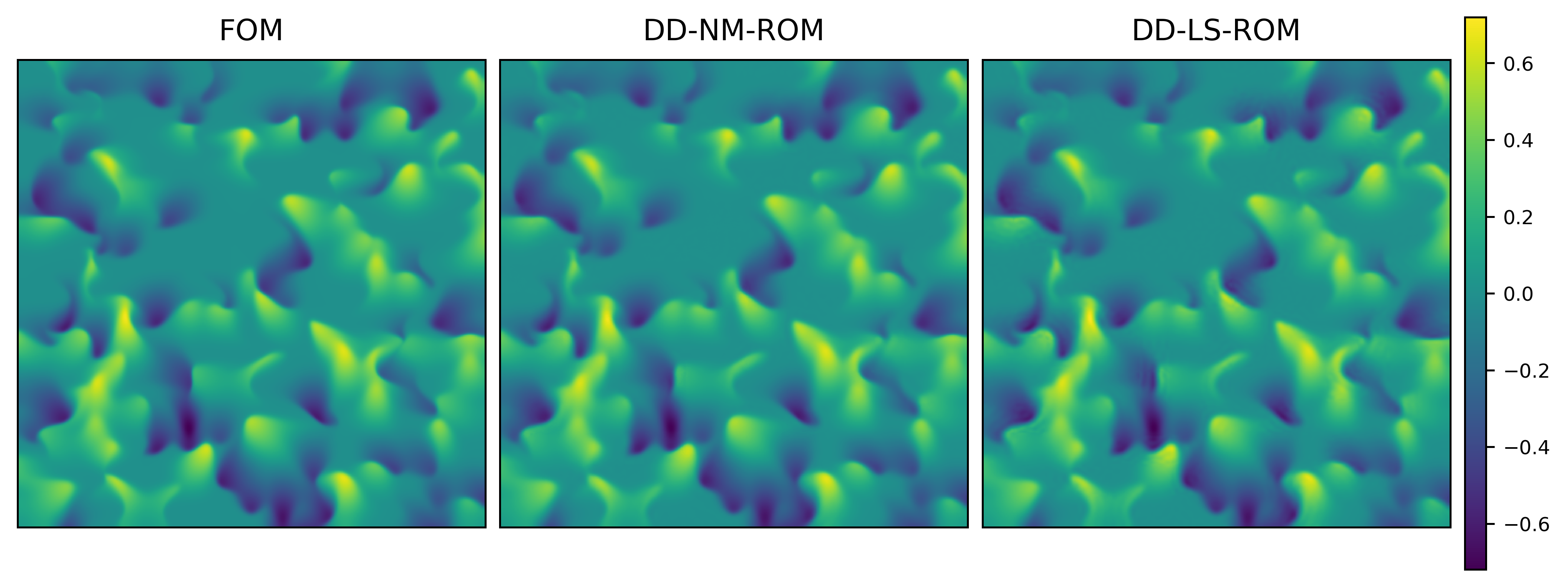}\\
\includegraphics[width=0.86\textwidth]{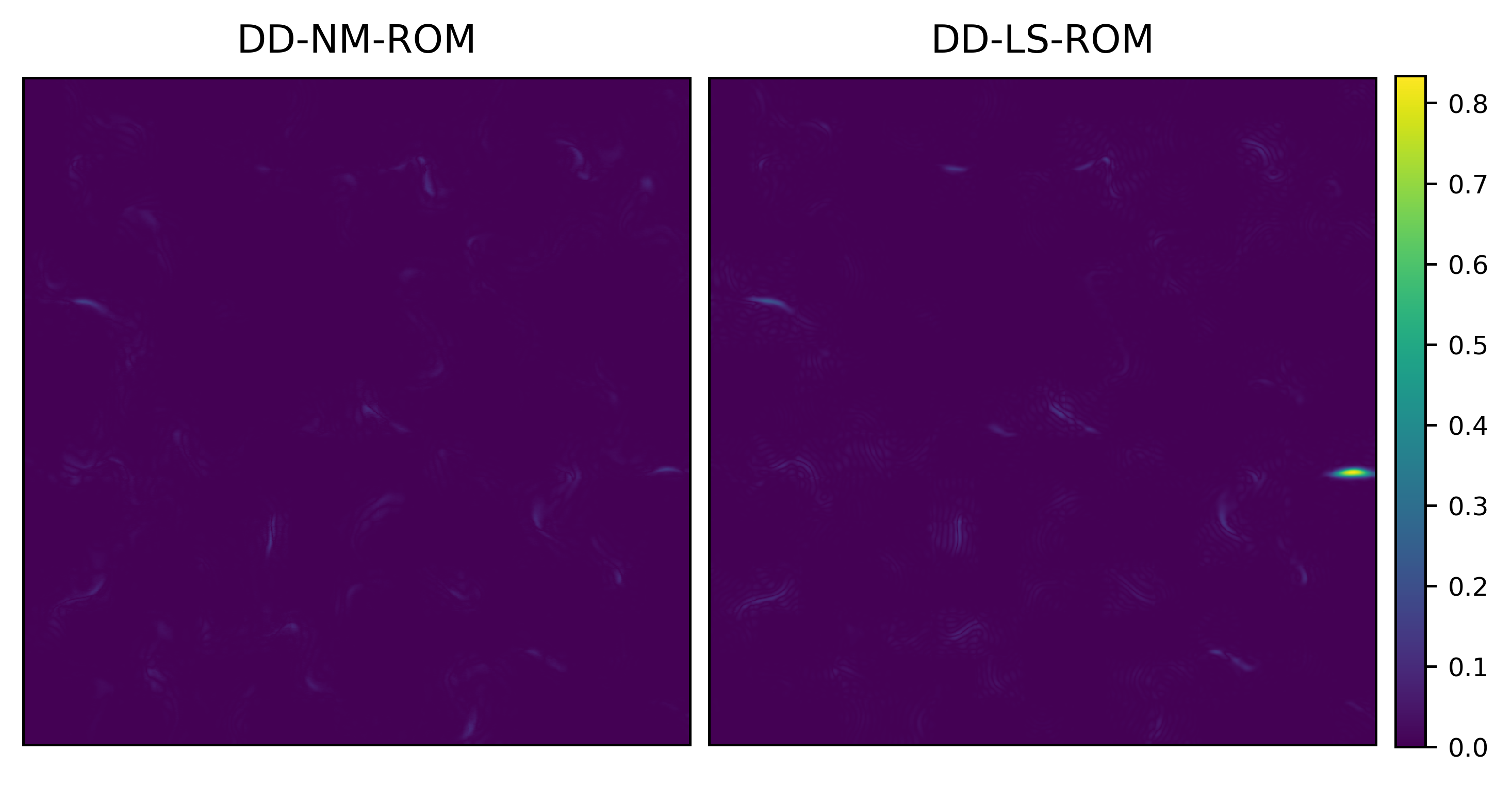}
\caption{Burgers \(v\)-velocity field and pointwise absolute error at \(t=1\) for the same representative case as Figure~\ref{fig:burgers_app_u_solution}.
The top row shows the FOM, DD-NM-ROM with AE \((50,50)\), and DD-LS-ROM with POD \((65,65)\), from left to right.
The bottom row shows the corresponding DD-NM-ROM and DD-LS-ROM errors.
A common color scale is used within each row.}
\label{fig:burgers_app_v_solution}
\label{fig:burgers_app_v_error}
\end{figure}

Figures~\ref{fig:burgers_app_u_solution} and~\ref{fig:burgers_app_v_solution} show the \(u\) and \(v\) components, respectively, at \(t=1\), corresponding to Figure~\ref{fig:burgers_speed_solution} in the main text.
Whereas speed summarizes both components in a single scalar, the field-wise comparisons allow the state and spatial error in each direction to be examined separately.
The values in Table~\ref{tab:burgers_scaleup_results} concern the full time trajectories of 10 cases and should be distinguished from these pointwise final-time errors.

\subsubsection{Velocity and pressure for Navier--Stokes}
\label{app:ns_scaleup_fields}
\begin{figure}[pos=p]
\centering
\includegraphics[width=0.98\textwidth]{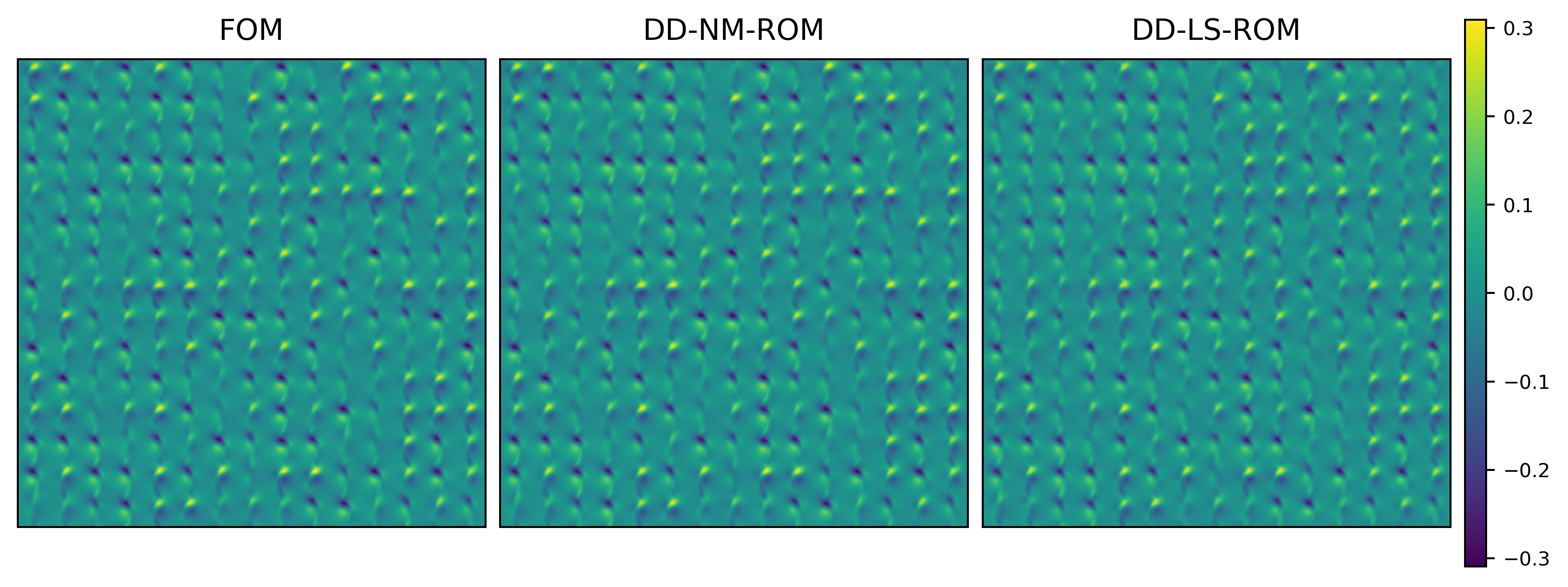}\\
\includegraphics[width=0.86\textwidth]{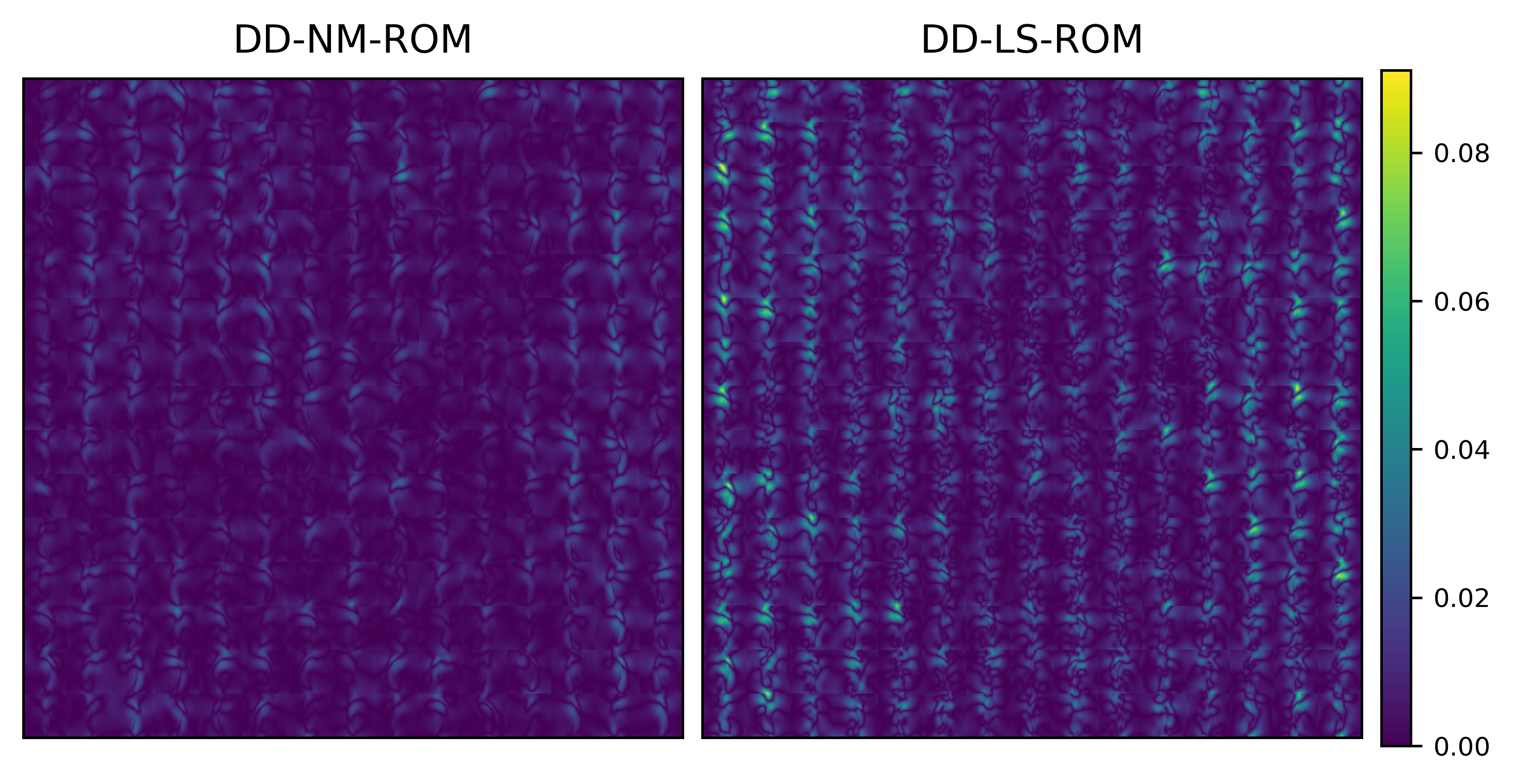}
\caption{Navier--Stokes \(u\)-velocity field and pointwise absolute error at \(t=5\) for the representative case used in Figure~\ref{fig:scaleup_vorticity_solution}.
The top row shows the FOM, DD-NM-ROM with AE \((10,10,10)\), and DD-LS-ROM with POD \((46,42,37)\), from left to right.
The bottom row shows the corresponding DD-NM-ROM and DD-LS-ROM errors defined in Eq.~\eqref{eq:primitive_visual_error}.
A common color scale is used within each row.}
\label{fig:app_scaleup_u_solution}
\label{fig:app_scaleup_u_error}
\end{figure}

\begin{figure}[pos=p]
\centering
\includegraphics[width=0.98\textwidth]{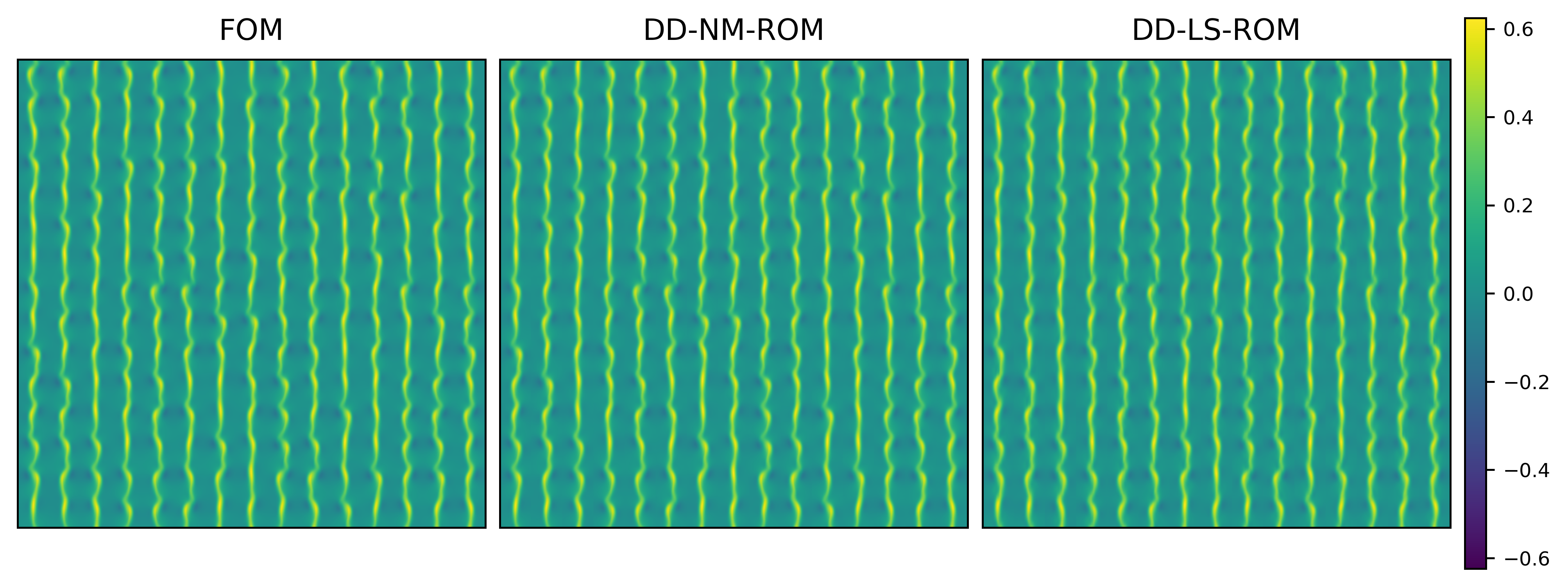}\\
\includegraphics[width=0.86\textwidth]{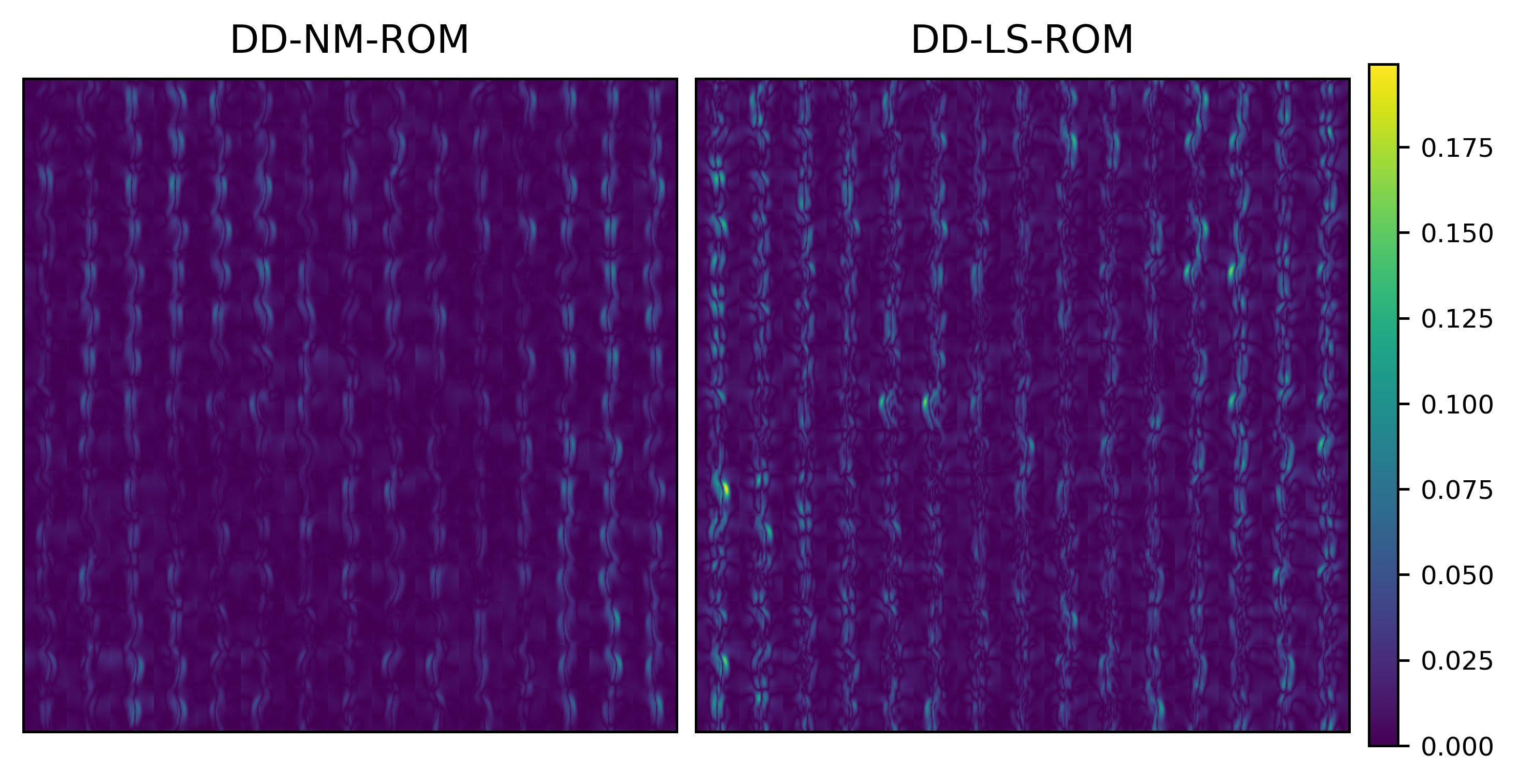}
\caption{Navier--Stokes \(v\)-velocity field and pointwise absolute error at \(t=5\) for the same representative case as Figure~\ref{fig:app_scaleup_u_solution}.
The top row shows the FOM, DD-NM-ROM with AE \((10,10,10)\), and DD-LS-ROM with POD \((46,42,37)\), from left to right.
The bottom row shows the corresponding DD-NM-ROM and DD-LS-ROM errors.
A common color scale is used within each row.}
\label{fig:app_scaleup_v_solution}
\label{fig:app_scaleup_v_error}
\end{figure}

\begin{figure}[pos=p]
\centering
\includegraphics[width=0.98\textwidth]{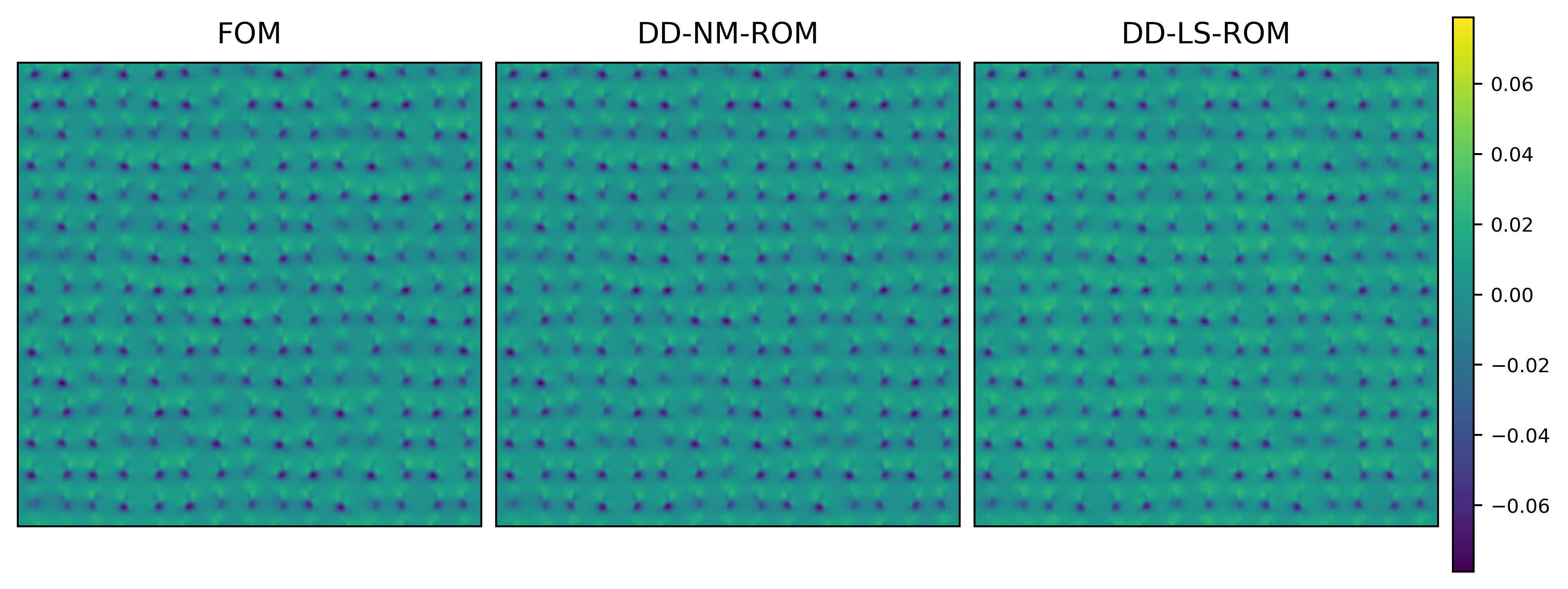}\\
\includegraphics[width=0.86\textwidth]{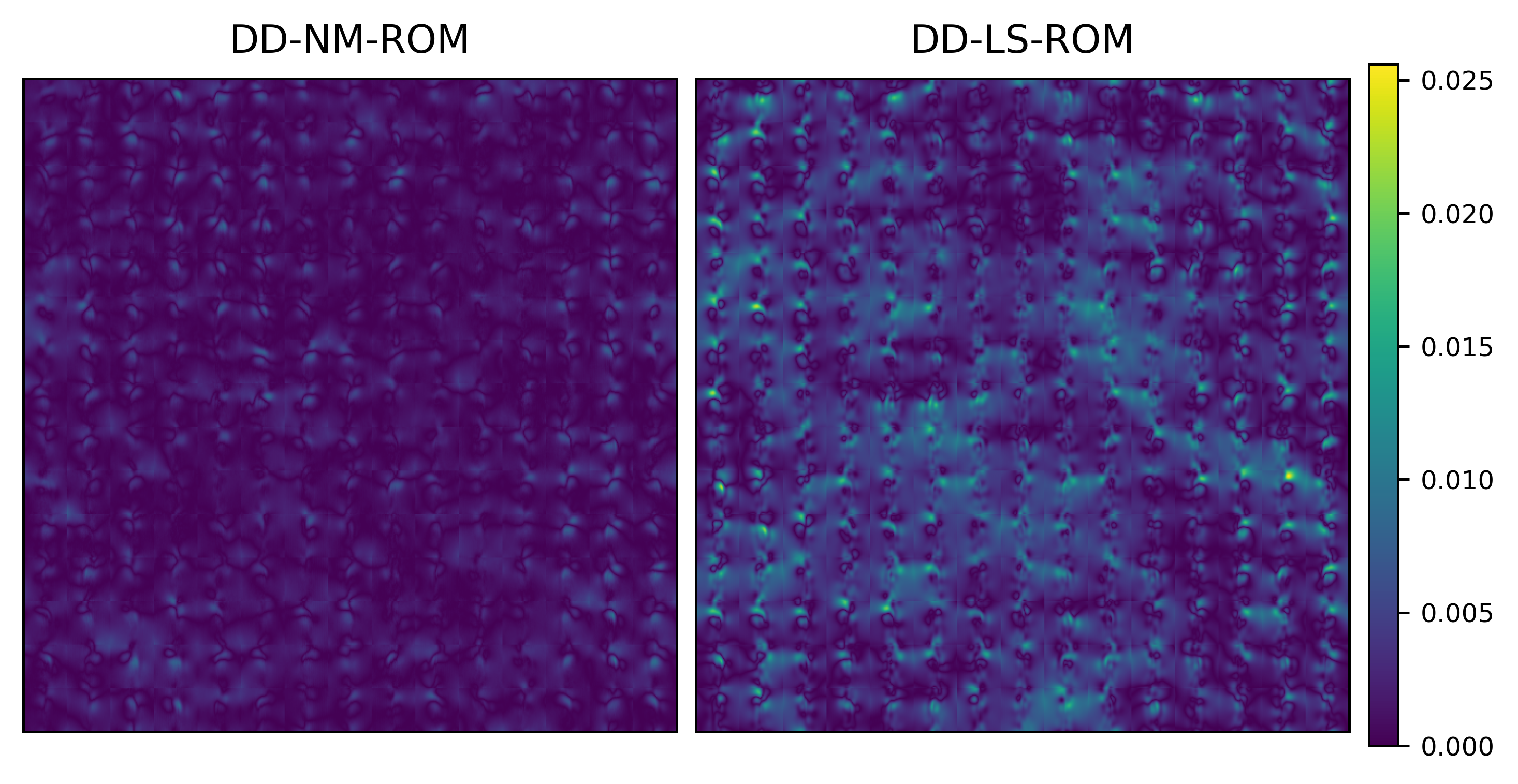}
\caption{Navier--Stokes pressure field and pointwise absolute error at \(t=5\) for the same representative case as Figure~\ref{fig:app_scaleup_u_solution}.
The top row shows the FOM, DD-NM-ROM with AE \((10,10,10)\), and DD-LS-ROM with POD \((46,42,37)\), from left to right.
The bottom row shows the corresponding DD-NM-ROM and DD-LS-ROM errors.
A common color scale is used within each row.}
\label{fig:app_scaleup_p_solution}
\label{fig:app_scaleup_p_error}
\end{figure}

Figures~\ref{fig:app_scaleup_u_solution}--\ref{fig:app_scaleup_p_solution} compare \(u,v,p\) at \(t=5\) for the same representative case as Figure~\ref{fig:scaleup_vorticity_solution} in the main text.
Unlike the main-text vorticity, which is a flow diagnostic obtained from spatial derivatives of velocity, these figures show the primitive variables directly reconstructed by the decoder and used in residual evaluation.
As for the other fields, the pressure plots show the FOM and ROM states and pointwise absolute errors.

\clearpage
\printcredits

\section*{Declaration of competing interest}
The authors declare that they have no known competing financial
interests or personal relationships that could have appeared to
influence the work reported in this paper.

\section*{Acknowledgments}
This work was performed at Lawrence Livermore National Laboratory.
This research was supported in part by the U.S. Department of Energy, Office of Science, Advanced Scientific Computing Research program, through the Scientific Discovery through Advanced Computing (SciDAC) Institute ``LEADS: LEarning-Accelerated Domain Science'' under Project No.~SCW1933.
Lawrence Livermore National Laboratory is operated by Lawrence Livermore National Security, LLC, for the U.S. Department of Energy, National Nuclear Security Administration, under Contract DE-AC52-07NA27344.
IM release number: LLNL-JRNL-2024535.

\section*{Data availability}
The code used in this study, including the data generation scripts, is
currently undergoing internal review at Lawrence Livermore National
Laboratory (LLNL) for public release. Upon approval, the code will be
made available in a single GitHub repository. The simulation datasets
will not be distributed separately because of their large size.
They can be generated using the code and instructions that will be
provided in the repository.
Interested readers may contact the corresponding author by email
to be notified when the code is released.

\section*{Declaration of generative AI and AI-assisted technologies in the manuscript preparation process}
In this work, we used ChatGPT (GPT-6 Pro) for English translation, grammar
correction, and language editing of the manuscript drafted by the authors.

We reviewed the AI-assisted manuscript to ensure that the original technical
content and intended meaning were preserved. 

\bibliographystyle{cas-model2-names}

\bibliography{refs/ref.bib}

\end{document}